\documentclass[leqno,draft]{article}

\usepackage{makeidx}
\makeindex

\def\supp{\mathop{\rm supp}}

\newtheorem{theorem}{Theorem}
\newtheorem{lemma}[theorem]{Lemma}
\newtheorem{proposition}[theorem]{Proposition}
\newtheorem{definition}[theorem]{Definition}
\newtheorem{corollary}[theorem]{Corollary}

\newcommand{\begintheorem}{\addtocounter{equation}{1}\begin{theorem}}
\newcommand{\beginlemma}{\addtocounter{equation}{1}\begin{lemma}}
\newcommand{\beginproposition}{\addtocounter{equation}{1}\begin{proposition}}
\newcommand{\begindefinition}{\addtocounter{equation}{1}\begin{definition}}
\newcommand{\begincorollary}{\addtocounter{equation}{1}\begin{corollary}}

\begin{document}

\title{Some topics related to metrics and norms, including
  ultrametrics and ultranorms}

\author{Stephen Semmes \\
        Rice University}

\date{}

\maketitle

\begin{abstract}
Here we look at some geometric properties related to connectedness and
topological dimension $0$, especially in connection with norms on
vector spaces over fields with absolute value functions, which may be
non-archimedian.
\end{abstract}

\tableofcontents

\part{Metrics and norms}
\label{metrics, norms}

\section{$q$-Metrics}
\label{q-metrics}
\setcounter{equation}{0}

        Let $M$ be a set, and let $q$ be a positive real number.
A nonnegative real-valued function $d(x, y)$ defined for $x, y \in M$
is said to be a \emph{$q$-metric}\index{qmetrics@$q$-metrics} on $M$
if it satisfies the following three conditions.  First,
\begin{equation}
\label{d(x, y) = 0 if and only if x = y}
        d(x, y) = 0 \quad\hbox{if and only if}\quad x = y.
\end{equation}
Second,
\begin{equation}
\label{d(x, y) = d(y, x) for every x, y in M}
        d(x, y) = d(y, x) \quad\hbox{for every } x, y \in M.
\end{equation}
Third,
\begin{equation}
\label{d(x, z)^q le d(x, y)^q + d(y, z)^q for every x, y, z in M}
 d(x, z)^q \le d(x, y)^q + d(y, z)^q \quad\hbox{for every } x, y, z \in M.
\end{equation}
Of course, (\ref{d(x, z)^q le d(x, y)^q + d(y, z)^q for every x, y, z
  in M}) is the version of the triangle inequality associated to $q$.
If this holds with $q = 1$, then we simply say that $d(x, y)$ is a
\emph{metric}\index{metrics} on $M$.  Thus $d(x, y)$ is a $q$-metric
on $M$ if and only if $d(x, y)^q$ is a metric on $M$.

        Similarly, a nonnegative real-valued function $d(x, y)$
defined for $x, y \in M$ is said to be an
\emph{ultrametric}\index{ultrametrics} on $M$ if it satisfies
(\ref{d(x, y) = 0 if and only if x = y}), (\ref{d(x, y) = d(y, x) for
  every x, y in M}), and
\begin{equation}
\label{d(x, z) le max(d(x, y), d(y, z)) for every x, y, z in M}
 d(x, z) \le \max(d(x, y), d(y, z)) \quad\hbox{for every } x, y, z \in M,
\end{equation}
instead of (\ref{d(x, z)^q le d(x, y)^q + d(y, z)^q for every x, y, z
  in M}).  Clearly, for each $q > 0$, (\ref{d(x, z) le max(d(x, y),
  d(y, z)) for every x, y, z in M}) is equivalent to asking that
\begin{equation}
\label{d(x, z)^q le max(d(x, y)^q, d(y, z)^q) for every x, y, z in M}
 d(x, z)^q \le \max(d(x, y)^q, d(y, z)^q) \quad\hbox{for every } x, y, z \in M.
\end{equation}
If $d(x, y)$ is an ultrametric on $M$, then it follows that $d(x, y)$
is a $q$-metric on $M$ for every $q > 0$, since (\ref{d(x, z)^q le
  max(d(x, y)^q, d(y, z)^q) for every x, y, z in M}) implies
(\ref{d(x, z)^q le d(x, y)^q + d(y, z)^q for every x, y, z in M}).  In
this case, we also get that $d(x, y)^q$ is an ultrametric on $M$ for
every $q > 0$.  Note that the discrete metric\index{discrete metric}
on any set $M$ is an ultrametric, which is defined by putting $d(x,
y)$ equal to $1$ when $x \ne y$, and to $0$ when $x = y$.

        It is sometimes convenient to reformulate (\ref{d(x, z)^q le 
d(x, y)^q + d(y, z)^q for every x, y, z in M}) as saying that
\begin{equation}
\label{d(x, z) le (d(x, y)^q + d(y, z)^q)^{1/q} for every x, y, z in M}
 d(x, z) \le (d(x, y)^q + d(y, z)^q)^{1/q} \quad\hbox{for every } x, y, z \in M.
\end{equation}
Observe that
\begin{equation}
\label{max(a, b) le (a^q + b^q)^{1/q} le 2^{1/q} max(a, b)}
        \max(a, b) \le (a^q + b^q)^{1/q} \le 2^{1/q} \, \max(a, b)
\end{equation}
for any pair $a$, $b$ of nonnegative real numbers, which implies that
\begin{equation}
\label{lim_{q to infty} (a^q + b^q)^{1/q} = max(a, b)}
        \lim_{q \to \infty} (a^q + b^q)^{1/q} = \max(a, b).
\end{equation}
Thus (\ref{d(x, z) le max(d(x, y), d(y, z)) for every x, y, z in M})
corresponds to taking the limit as $q \to \infty$ in (\ref{d(x, z) le
  (d(x, y)^q + d(y, z)^q)^{1/q} for every x, y, z in M}), so that one
might think of an ultrametric as being a $q$-metric with $q = \infty$.

        If $0 < q_1 < q_2 < \infty$, then
\begin{equation}
\label{a^{q_2} + b^{q_2} le max(a, b)^{q_2 - q_1} (a^{q_1} + b^{q_1})}
        a^{q_2} + b^{q_2} \le \max(a, b)^{q_2 - q_1} \, (a^{q_1} + b^{q_1})
\end{equation}
for every $a, b \ge 0$.  We also have that
\begin{equation}
\label{max(a, b) le (a^{q_1} + b^{q_1})^{1/q_1}}
        \max(a, b) \le (a^{q_1} + b^{q_1})^{1/q_1},
\end{equation}
as in (\ref{max(a, b) le (a^q + b^q)^{1/q} le 2^{1/q} max(a, b)}), so that
\begin{equation}
\label{a^{q_2} + b^{q_2} le ... = (a^{q_1} + b^{q_1})^{q_2/q_1}}
        a^{q_2} + b^{q_2} \le (a^{q_1} + b^{q_1})^{((q_2 - q_1)/ q_1) + 1}
                            = (a^{q_1} + b^{q_1})^{q_2/q_1}
\end{equation}
for every $a, b \ge 0$.  Equivalently,
\begin{equation}
\label{(a^{q_2} + b^{q_2})^{1/q_2} le (a^{q_1} + b^{q_1})^{1/q_1}}
        (a^{q_2} + b^{q_2})^{1/q_2} \le (a^{q_1} + b^{q_1})^{1/q_1}
\end{equation}
for every $a, b \ge 0$ when $0 < q_1 < q_2 < \infty$.  If $d(x, y)$ is
a $q_2$-metric on $M$ for some $q_2 > 0$, then it follows that $d(x,
y)$ is a $q_1$-metric on $M$ as well when $0 < q_1 < q_2$.  Of course,
the topology on $M$ determined by the metric $d(x, y)^{q_1}$ is the
same as the topology on $M$ determined by the metric $d(x, y)^{q_2}$
in this case.

\section{Open and closed balls}
\label{open, closed balls}
\setcounter{equation}{0}

        Let $M$ be a set, and suppose that $d(x, y)$ is a $q$-metric
on $M$ for some positive real number $q$.  If $x \in M$ and $r$
is a positive real number, then the open ball\index{open balls} centered
at $x$ with radius $r$ is defined as usual by
\begin{equation}
\label{B(x, r) = {z in M : d(x, z) < r}}
        B(x, r) = \{z \in M : d(x, z) < r\}.
\end{equation}
Equivalently,
\begin{equation}
\label{B(x, r) = {z in M : d(x, z)^q < r^q}}
        B(x, r) = \{z \in M : d(x, z)^q < r^q\},
\end{equation}
which is the open ball in $M$ centered at $x$ with radius $r^q$ with
respect to the metric $d(\cdot, \cdot)^q$.  If $y \in B(x, r)$, so
that $d(x, y)^q < r^q$, then let $t$ be the positive real number
determined by
\begin{equation}
\label{t^q = r^q - d(x, y)^q}
        t^q = r^q - d(x, y)^q.
\end{equation}
It is easy to see that
\begin{equation}
\label{B(y, t) subseteq B(x, r)}
        B(y, t) \subseteq B(x, r)
\end{equation}
under these conditions, because $d(\cdot, \cdot)^q$ is a metric
on $M$.

        Let us say that $U \subseteq M$ is an open set\index{open sets}
in $M$ if for each $x \in U$ there is an $r > 0$ such that
\begin{equation}
\label{B(x, r) subseteq U}
        B(x, r) \subseteq U.
\end{equation}
This is analogous to the standard definition for metric spaces, and it
is equivalent to saying that $U$ is an open set in $M$ with respect to
the metric $d(\cdot, \cdot)^q$, because of (\ref{B(x, r) = {z in M :
    d(x, z)^q < r^q}}).  In particular, this defines a topology on
$M$, which is the same as the topology on $M$ determined by the metric
$d(\cdot, \cdot)^q$.  Open balls in $M$ are open sets with respect to
this topology, by (\ref{B(y, t) subseteq B(x, r)}), which corresponds
to the standard argument for metric spaces.

        Similarly, if $d(\cdot, \cdot)$ is an ultrametric on $M$, then
\begin{equation}
\label{B(y, r) subseteq B(x, r)}
        B(y, r) \subseteq B(x, r)
\end{equation}
for every $y \in B(x, r)$, which is the same as (\ref{B(y, t) subseteq
  B(x, r)}) with $t = r$.  More precisely, this holds when $d(x, y) <
r$, which is symmetric in $x$ and $y$.  Thus we also have that
\begin{equation}
\label{B(x, r) subseteq B(y, r)}
        B(x, r) \subseteq B(y, r)
\end{equation}
when $d(x, y) < r$, so that
\begin{equation}
\label{B(x, r) = B(y, r)}
        B(x, r) = B(y, r)
\end{equation}
in this situation.

        Let $d(\cdot, \cdot)$ be a $q$-metric on $M$ for some $q > 0$
again.  The closed ball\index{closed balls} in $M$ centered at $x \in
M$ with radius $r \ge 0$ with respect to $d(\cdot, \cdot)$ is defined
by
\begin{equation}
\label{overline{B}(x, r) = {z in M : d(x, z) le r}}
        \overline{B}(x, r) = \{z \in M : d(x, z) \le r\}.
\end{equation}
As before, this is the same as
\begin{equation}
\label{overline{B}(x, r) = {z in M : d(x, z)^q le r^q}}
        \overline{B}(x, r) = \{z \in M : d(x, z)^q \le r^q\},
\end{equation}
which is the closed ball in $M$ centered at $x$ with radius $r^q$ with
respect to the metric $d(\cdot, \cdot)^q$.  If $y \in \overline{B}(x,
r)$, so that $d(x, y)^q \le r^q$, then let $t$ be the nonnegative real
number determined by (\ref{t^q = r^q - d(x, y)^q}).  In analogy with
(\ref{B(y, t) subseteq B(x, r)}), we have that
\begin{equation}
\label{overline{B}(y, t) subseteq overline{B}(x, r)}
        \overline{B}(y, t) \subseteq \overline{B}(x, r),
\end{equation}
since $d(\cdot, \cdot)^q$ is a metric on $M$.  If $d(\cdot, \cdot)$ is
an ultrametric on $M$, then
\begin{equation}
\label{overline{B}(y, r) subseteq overline{B}(x, r)}
        \overline{B}(y, r) \subseteq \overline{B}(x, r)
\end{equation}
for every $y \in \overline{B}(x, r)$, which is the same as
(\ref{overline{B}(y, t) subseteq overline{B}(x, r)}) with $t = r$.
Equivalently, (\ref{overline{B}(y, r) subseteq overline{B}(x, r)})
holds when $d(x, y) \le r$, which is symmetric in $x$ and $y$.  Thus
the opposite inclusion also holds in this case, so that
\begin{equation}
\label{overline{B}(x, r) = overline{B}(y, r)}
        \overline{B}(x, r) = \overline{B}(y, r)
\end{equation}
for every $x, y \in M$ with $d(x, y) \le r$.  This implies that
closed balls in $M$ are open sets when $d(\cdot, \cdot)$ is an
ultrametric on $M$.

\section{Some related facts}
\label{some related facts}
\setcounter{equation}{0}

        If $d(x, y)$ is a $q$-metric on a set $M$ for some positive
real number $q$, then we can reexpress (\ref{d(x, z)^q le d(x, y)^q +
  d(y, z)^q for every x, y, z in M}) as
\begin{equation}
\label{d(x, z)^q - d(y, z)^q le d(x, y)^q for every x, y, z in M}
 d(x, z)^q - d(y, z)^q \le d(x, y)^q \quad\hbox{for every } x, y, z \in M.
\end{equation}
Of course, this is nontrivial only when $d(y, z)^q < d(x, z)^q$, which
is to say that
\begin{equation}
\label{d(y, z) < d(x, z)}
        d(y, z) < d(x, z).
\end{equation}
If $d(\cdot, \cdot)$ is an ultrametric on $M$, then (\ref{d(x, z) le
  max(d(x, y), d(y, z)) for every x, y, z in M}) and (\ref{d(y, z) <
  d(x, z)}) imply that
\begin{equation}
\label{d(x, z) le d(x, y)}
        d(x, z) \le d(x, y).
\end{equation}
In this case, we also have that
\begin{equation}
\label{d(x, y) le max(d(x, z), d(z, y)) le d(x, z)}
        d(x, y) \le \max(d(x, z), d(z, y)) \le d(x, z)
\end{equation}
when $d(y, z) \le d(x, z)$.  It follows that
\begin{equation}
\label{d(x, y) = d(x, z)}
        d(x, y) = d(x, z)
\end{equation}
when $d(\cdot, \cdot)$ is an ultrametric on $M$ and $x, y, z \in M$
satisfy (\ref{d(y, z) < d(x, z)}), by combining (\ref{d(x, z) le d(x,
  y)}) and (\ref{d(x, y) le max(d(x, z), d(z, y)) le d(x, z)}).

        Let $d(\cdot, \cdot)$ be a $q$-metric on $M$ for some
$q > 0$ again, and put
\begin{equation}
\label{V(x, r) = {z in M : d(x, z) > r} = {z in M : d(x, z)^q > r^q}}
 V(x, r) = \{z \in M : d(x, z) > r\} = \{z \in M : d(x, z)^q > r^q\}
\end{equation}
for every $x \in M$ and $r \ge 0$, which is the same as the complement
of $\overline{B}(x, r)$ in $M$.  If $z \in V(x, r)$, then let $t$ be
the positive real number determined by
\begin{equation}
\label{t^q = d(x, z)^q - r^q}
        t^q = d(x, z)^q - r^q.
\end{equation}
If $y \in B(z, t)$, so that $d(y, z)^q < t^q$, then (\ref{d(x, z)^q -
  d(y, z)^q le d(x, y)^q for every x, y, z in M}) implies that $d(x,
y)^q > r^q$, which means that $y \in V(x, r)$.  This shows that
\begin{equation}
\label{B(z, t) subseteq V(x, r),}
        B(z, t) \subseteq V(x, r),
\end{equation}
which implies that $V(x, r)$ is an open set in $M$, and hence that
$\overline{B}(x, r)$ is a closed set in $M$.  If $d(\cdot, \cdot)$ is
an ultrametric on $M$, then (\ref{B(z, t) subseteq V(x, r),}) holds
with $t = d(x, z)$, because of (\ref{d(x, y) = d(x, z)}).

       Let $d(\cdot, \cdot)$ be any $q$-metric on $M$ again, and put
\begin{equation}
\label{W(x, r) = {z in M : d(x, z) ge r} = {z in M : d(x, z)^q ge r^q}}
 W(x, r) = \{z \in M : d(x, z) \ge r\} = \{z \in M : d(x, z)^q \ge r^q\}
\end{equation}
for every $x \in M$ and $r > 0$, which is the same as the complement
of $B(x, r)$ in $M$.  If $z \in W(x, r)$, then let $t$ be the
nonnegative real number determined by (\ref{t^q = d(x, z)^q - r^q}).
If $y \in \overline{B}(z, t)$, so that $d(y, z)^q \le t^q$, then
(\ref{d(x, z)^q - d(y, z)^q le d(x, y)^q for every x, y, z in M})
implies that $d(x, y)^q \ge r^q$, and thus $y \in W(x, r)$.  It
follows that
\begin{equation}
\label{overline{B}(z, t) subseteq W(x, r)}
        \overline{B}(z, t) \subseteq W(x, r)
\end{equation}
under these conditions, which is trivial when $d(x, z) = r$, so that
$t = 0$.  Note that $W(x, r)$ is a closed set in $M$ for every $x \in
M$ and $r > 0$, since it is the complement of an open set.

        If $d(\cdot, \cdot)$ is an ultrametric on $M$, then
\begin{equation}
\label{B(z, d(x, z)) subseteq W(x, r)}
        B(z, d(x, z)) \subseteq W(x, r)
\end{equation}
for every $z \in W(x, r)$.  More precisely, if $y \in B(z, d(x, z))$,
then (\ref{d(y, z) < d(x, z)}) holds, which implies that (\ref{d(x, y)
  = d(x, z)}) holds as well.  If we also have $z \in W(x, r)$, then it
follows that
\begin{equation}
\label{d(x, y) = d(x, z) ge r}
        d(x, y) = d(x, z) \ge r
\end{equation}
for every $y \in B(z, d(x, z))$, so that $y \in W(x, r)$, as desired.
In particular, this shows that $W(x, r)$ is an open set in $M$ for
every $x \in M$ and $r > 0$ when $d(\cdot, \cdot)$ is an ultrametric
on $M$.  Thus $B(x, r)$ is a closed set in $M$ for every $x \in M$ and
$r > 0$ in this case.

\section{Absolute value functions}
\label{absolute value functions}
\setcounter{equation}{0}

        Let $k$ be a field, and let $q$ be a positive real number
again.  A nonnegative real-valued function $|\cdot|$ defined on $k$
is said to be a \emph{$q$-absolute value function}\index{qabsolute
value functions@$q$-absolute value functions} if it satisfies
the following three conditions.  First, for each $x \in k$,
\begin{equation}
\label{|x| = 0 if and only if x = 0}
        |x| = 0 \quad\hbox{if and only if}\quad x = 0.
\end{equation}
Second,
\begin{equation}
\label{|x y| = |x| |y| for every x, y in k}
        |x \, y| = |x| \, |y| \quad\hbox{for every } x, y \in k.
\end{equation}
Third,
\begin{equation}
\label{|x + y|^q le |x|^q + |y|^q for every x, y in k}
        |x + y|^q \le |x|^q + |y|^q \quad\hbox{for every } x, y \in k.
\end{equation}
If (\ref{|x + y|^q le |x|^q + |y|^q for every x, y in k}) holds with
$q = 1$, then we may simply say that $|\cdot|$ is an \emph{absolute
  value function}\index{absolute value functions} on $k$.
Equivalently, $|x|$ is a $q$-absolute value function on $k$ if and
only if $|x|^q$ is an absolute value function on $k$.

        Suppose for the moment that $|\cdot|$ is a nonnegative
real-valued function on $k$ that satisfies (\ref{|x| = 0 if and only
  if x = 0}) and (\ref{|x y| = |x| |y| for every x, y in k}).
Let us use $1$ to denote the multiplicative identity element in $k$,
as well as usual positive integer, depending on the context.  Thus
$|1| > 0$, by (\ref{|x| = 0 if and only if x = 0}), since $1 \ne 0$
in $k$, by definition of a field.  We also have that
\begin{equation}
\label{|1| = |1^2| = |1|^2}
        |1| = |1^2| = |1|^2,
\end{equation}
which implies that
\begin{equation}
\label{|1| = 1}
        |1| = 1.
\end{equation}
Similarly, if $x \in k$ satisfies $x^n = 1$ for some positive integer
$n$, then we get that
\begin{equation}
\label{|x^n| = |x|^n = 1}
        |x^n| = |x|^n = 1,
\end{equation}
and hence
\begin{equation}
\label{|x| = 1}
        |x| = 1.
\end{equation}

        Let $-x$ be the additive inverse of $x \in k$, which is equal to
$(-1) \, x$, where $-1$ is the additive inverse of $1$ in $k$.  In
particular, $(-1)^2 = 1$ in $k$, which implies that
\begin{equation}
\label{|-1| = 1}
        |-1| = 1,
\end{equation}
as before.  It follows that
\begin{equation}
\label{|-x| = x}
        |-x| = x
\end{equation}
for every $x \in k$.  If $|\cdot|$ is a $q$-absolute value function
on $k$, then we get that
\begin{equation}
\label{d(x, y) = |x - y|}
        d(x, y) = |x - y|
\end{equation}
defines a $q$-metric on $k$.  More precisely, this uses (\ref{|-x| =
  x}) to get the symmetry condition (\ref{d(x, y) = d(y, x) for every
  x, y in M}).

        A nonnegative real-valued function $|\cdot|$ on $k$ is said
to be an \emph{ultrametric absolute value function}\index{ultrametric
absolute value functions} on $k$ if
\begin{equation}
\label{|x + y| le max(|x|, |y|) for every x, y in k}
        |x + y| \le \max(|x|, |y|) \quad\hbox{for every } x, y \in k.
\end{equation}
This implies that (\ref{d(x, y) = |x - y|}) defines an ultrametric on
$k$.  As before, for each $q > 0$, (\ref{|x + y| le max(|x|, |y|) for
  every x, y in k}) is equivalent to asking that
\begin{equation}
\label{|x + y|^q le max(|x|^q, |y|^q) for every x, y in k}
        |x + y|^q \le \max(|x|^q, |y|^q) \quad\hbox{for every } x, y \in k.
\end{equation}
If $|\cdot|$ is an ultrametric absolute value function on $k$, then it
follows that $|\cdot|$ is a $q$-absolute value function on $k$ for
every $q > 0$, because (\ref{|x + y|^q le max(|x|^q, |y|^q) for every
  x, y in k}) implies (\ref{|x + y|^q le |x|^q + |y|^q for every x, y
  in k}).  We also get that $|x|^q$ is an ultrametric absolute value
function on $k$ for every $q > 0$ in this case.

        As in Section \ref{q-metrics}, we can reformulate
(\ref{|x + y|^q le |x|^q + |y|^q for every x, y in k}) as saying that
\begin{equation}
\label{|x + y| le (|x|^q + |y|^q)^{1/q} for every x, y in k}
 |x + y| \le (|x|^q + |y|^q)^{1/q} \quad\hbox{for every } x, y \in k.
\end{equation}
Using (\ref{lim_{q to infty} (a^q + b^q)^{1/q} = max(a, b)}), (\ref{|x
  + y| le max(|x|, |y|) for every x, y in k}) corresponds to taking
the limit as $q \to \infty$ in (\ref{|x + y| le (|x|^q + |y|^q)^{1/q}
  for every x, y in k}), so that an ultrametric absolute value
function may be considered as a $q$-absolute value function with $q =
\infty$.  We have also seen that the right side of the inequality in
(\ref{|x + y| le (|x|^q + |y|^q)^{1/q} for every x, y in k}) decreases
monotonically in $q$, by (\ref{(a^{q_2} + b^{q_2})^{1/q_2} le (a^{q_1}
  + b^{q_1})^{1/q_1}}).  If $0 < q_1 < q_2 < \infty$, and if $|\cdot|$
is a $q_2$-absolute value function on $k$, then it follows that
$|\cdot|$ is a $q_1$-absolute value function on $k$ too.

        The \emph{trivial absolute value function}\index{trivial absolute
value function} may be defined on any field $k$ by putting $|x| = 1$ when
$x \ne 0$, and $|0| = 0$.  It is easy to see that this defines
an ultrametric absolute value function on $k$, for which the
corresponding ultrametric (\ref{d(x, y) = |x - y|}) is the same as the
discrete metric.  Suppose for the moment that $|\cdot|$ is a
nonnegative real-valued function on a field $k$ that satisfies
(\ref{|x| = 0 if and only if x = 0}) and (\ref{|x y| = |x| |y| for
  every x, y in k}), and which is not the trivial absolute
value function on $k$.  This means that there is an $x \in k$ such that
$x \ne 0$ and $|x| \ne 1$, and we may as well suppose that
\begin{equation}
\label{0 < |x| < 1}
        0 < |x| < 1,
\end{equation}
since otherwise we can replace $x$ with $1/x$.  Of course,
we also have that
\begin{equation}
\label{|1/x| = 1/|x| > 1}
        |1/x| = 1/|x| > 1
\end{equation}
in this case.

        It is well known that the standard absolute value functions
on the fields ${\bf R}$,\index{R@${\bf R}$} ${\bf C}$\index{C@${\bf
    C}$} of real and complex numbers are absolute value functions in
the sense described in this section.  Hence they are also $q$-absolute
value functions when $0 < q < 1$, as before.  However, it is easy to
see that they are not $q$-absolute value functions when $q > 1$, even
when restricted to the field ${\bf Q}$\index{Q@${\bf Q}$} of rational
numbers.

\section{Some additional properties}
\label{some additional properties}
\setcounter{equation}{0}

        Let $k$ be a field, and let ${\bf Z}_+$\index{Z_+@${\bf Z}_+$}
be the set of positive integers.  If $x \in k$ and $n \in {\bf Z}_+$,
then let $n \cdot x$ be the sum of $n$ $x$'s in $k$.  Note that
\begin{equation}
\label{n_1 cdot (n_2 cdot x) = (n_1 n_2) cdot x}
        n_1 \cdot (n_2 \cdot x) = (n_1 \, n_2) \cdot x
\end{equation}
for every $x \in k$ and $n_1, n_2 \in {\bf Z}_+$, and
\begin{equation}
\label{n cdot (x y) = (n cdot x) y = x (n cdot y)}
        n \cdot (x \, y) = (n \cdot x) \, y = x \, (n \cdot y)
\end{equation}
for every $x, y \in k$ and $n \in {\bf Z}_+$.  In particular,
\begin{equation}
\label{n^j cdot 1 = (n cdot 1)^j}
        n^j \cdot 1 = (n \cdot 1)^j
\end{equation}
for every $j, n \in {\bf Z}_+$.

        An absolute value function $|\cdot|$ on $k$ is said to be
\emph{archimedian}\index{archimedian absolute value functions}
if there are positive integers $n$ such that $|n \cdot 1|$ is as large
as one wants.  Equivalently, $|\cdot|$ is archimedian when
\begin{equation}
\label{|n cdot 1| > 1}
        |n \cdot 1| > 1
\end{equation}
for some $n \in {\bf Z}_+$, since this implies that
\begin{equation}
\label{|n^j cdot 1| = |n cdot 1|^j to infty as j to infty}
        |n^j \cdot 1| = |n \cdot 1|^j \to \infty \quad\hbox{as } j \to \infty,
\end{equation}
by (\ref{n^j cdot 1 = (n cdot 1)^j}).  Otherwise, $|\cdot|$ is
\emph{non-archimedian}\index{non-archimedian absolute value functions}
when
\begin{equation}
\label{|n cdot 1| le 1}
        |n \cdot 1| \le 1
\end{equation}
for every $n \in {\bf Z}_+$.  The previous argument shows that it is
enough to check that $|n \cdot 1|$ is bounded for $n \in {\bf Z}_+$,
to get that $|\cdot|$ is non-archimedian.  It is easy to see that
ultrametric absolute value functions are non-archimedian, using the
ultrametric version (\ref{|x + y| le max(|x|, |y|) for every x, y in
  k}) of the triangle inequality.  Conversely, it can be shown that
non-archimedian absolute value functions satisfy the ultrametric
version of the triangle inequality, as in Lemma 1.5 on p16 of
\cite{c}, and Theorem 2.2.2 on p28 of \cite{fg}.  There is an
analogous statement for a $q$-absolute value function $|\cdot|$ on $k$
for any $q > 0$, which can be derived from the previous statement for
absolute value functions applied to $|x|^q$.

        A pair of absolute value functions $|\cdot|_1$, $|\cdot|_2$
on $k$ are said to be \emph{equivalent}\index{equivalent absolute
value functions} if there is a positive real number $a$ such that
\begin{equation}
\label{|x|_2 = |x|_1^a}
        |x|_2 = |x|_1^a
\end{equation}
for every $x \in k$.  This implies that the topologies on $k$
determined by the metrics associated to $|\cdot|_1$, $|\cdot|_2$ as in
(\ref{d(x, y) = |x - y|}) are the same.  Conversely, if the topologies
on $k$ determined by the metrics associated to $|\cdot|_1$ and
$|\cdot|_2$ are the same, then one can show that $|\cdot|_1$ and
$|\cdot|_2$ are equivalent on $k$, as in Lemma 3.2 on p20 of \cite{c},
and Lemma 3.1.2 on p42 of \cite{fg}.  Similarly, if $|\cdot|_1$ and
$|\cdot|$ are $q_1$ and $q_2$-absolute value functions on $k$ for some
$q_1, q_2 > 0$, then let us say that $|\cdot|_1$ and $|\cdot|_2$ are
equivalent when (\ref{|x|_2 = |x|_1^a}) holds for some $a > 0$.  This
is the same as saying that
\begin{equation}
\label{|x|_2^{q_2} = (|x|_1^{q_1})^{a q_2/q_1}}
        |x|_2^{q_2} = (|x|_1^{q_1})^{a \, q_2/q_1}
\end{equation}
for every $x \in k$, so that $|x|_1^{q_1}$ and $|x|_2^{q_2}$ are
equivalent as absolute value functions on $k$.

        Let $|\cdot|$ be an absolute value function on $k$, which
leads to a metric on $k$ as in (\ref{d(x, y) = |x - y|}), and hence
a topology on $k$.  Using standard arguments, one can check that
addition and multiplication on $k$ are continuous as mappings from
$k \times k$ into $k$, where $k \times k$ is equipped with the
corresponding product topology.  Similarly,
\begin{equation}
\label{x mapsto 1/x}
        x \mapsto 1/x
\end{equation}
is continuous as a mapping from $k \setminus \{0\}$ into itself.

        If $k$ is not already complete as a metric space with respect
to the metric associated to $|\cdot|$, then one can obtain a
completion of $k$ in the usual way.  The field operations on $k$ can
be extended to the completion in a natural way, so that the completion
of $k$ is also a field.  The absolute value function on $k$ can be
extended to an absolute value function on the completion of $k$ as
well, in such a way that the metric associated to the extension of the
absolute value function to the completion of $k$ is the same as the
metric already given on the completion of $k$.  If $|\cdot|$ is an
ultrametric absolute value function on $k$, then the extension of
$|\cdot|$ to the completion of $k$ is an ultrametric absolute value
function too.  Of course, $k$ is a dense subset of its completion,
by construction.

        Let $|\cdot|$ be an ultrametric absolute value function
on any field $k$.  If $x, y \in k$ and $|x - y| \le |x|$, then
\begin{equation}
\label{|y| le max(|x|, |x - y|) le |x|}
        |y| \le \max(|x|, |x - y|) \le |x|.
\end{equation}
If $|x - y| < |x|$, then
\begin{equation}
\label{|x| le max(|y|, |x - y|)}
        |x| \le \max(|y|, |x - y|)
\end{equation}
implies that $|x| \le |y|$.  Combining this with (\ref{|y| le max(|x|,
  |x - y|) le |x|}), we get that
\begin{equation}
\label{|x| = |y|}
        |x| = |y|
\end{equation}
when $|x - y| < |x|$.

        If $x \in k$ and $n$ is a nonnegative integer, then
\begin{equation}
\label{(1 - x) sum_{j = 0}^n x^j = 1 - x^{n + 1}}
        (1 - x) \, \sum_{j = 0}^n x^j = 1 - x^{n + 1},
\end{equation}
where $x^j$ is interpreted as being equal to $1$ when $j = 0$,
as usual.  It follows that
\begin{equation}
\label{sum_{j = 0}^n x^j = frac{1 - x^{n + 1}}{1 - x}}
        \sum_{j = 0}^n x^j = \frac{1 - x^{n + 1}}{1 - x}
\end{equation}
for every $n \ge 0$ when $x \ne 1$.  Let $|\cdot|$ be an absolute
value function on $k$, and suppose that $|x| < 1$, so that
\begin{equation}
\label{|x^{n + 1}| = |x|^{n + 1} to 0 as n to infty}
        |x^{n + 1}| = |x|^{n + 1} \to 0 \quad\hbox{as } n \to \infty.
\end{equation}
This implies that
\begin{equation}
\label{lim_{n to infty} sum_{j = 0}^n x^j = frac{1}{1 - x}}
        \lim_{n \to \infty} \sum_{j = 0}^n x^j = \frac{1}{1 - x}
\end{equation}
when $|x| < 1$, where the limit is taken with respect to the metric
associated to $|\cdot|$ on $k$.

\section{$p$-Adic numbers}
\label{p-adic numbers}
\setcounter{equation}{0}

        If $p$ is a prime number, then the \emph{$p$-adic absolute
value}\index{p-adic absolute value@$p$-adic absolute value} $|x|_p$
of a rational number $x$ is defined as follows.  Of course, $|0|_p = 0$.
Otherwise, if $x \ne 0$, then $x$ can be expressed as
\begin{equation}
\label{x = p^j (a / b)}
        x = p^j \, (a / b),
\end{equation}
where $a$, $b$, and $j$ are integers, $a, b \ne 0$, and neither $a$
nor $b$ is divisible by $p$.  In this case, we put
\begin{equation}
\label{|x|_p = p^{-j}}
        |x|_p = p^{-j}.
\end{equation}
One can check that this defines an ultrametric absolute value function
on ${\bf Q}$, so that the corresponding \emph{$p$-adic
  metric}\index{p-adic metric@$p$-adic metric}
\begin{equation}
\label{d_p(x, y) = |x - y|_p}
        d_p(x, y) = |x - y|_p
\end{equation}
is an ultrametric on ${\bf Q}$.

        The field ${\bf Q}_p$\index{Q_p@${\bf Q}_p$} of
\emph{$p$-adic numbers}\index{p-adic numbers@$p$-adic numbers}
is obtained by completing ${\bf Q}$ with respect to the $p$-adic
metric, as in the previous section.  The natural extension of the
$p$-adic absolute value function to ${\bf Q}_p$ is also called the
$p$-adic absolute value,\index{p-adic absolute value@$p$-adic absolute
value} and denoted $|\cdot|_p$.  Similarly, the natural extension of
the $p$-adic metric to ${\bf Q}_p$ is called the $p$-adic
metric\index{p-adic metric@$p$-adic metric} too, and is denoted
$d_p(\cdot, \cdot)$.  By construction, these extensions of the
$p$-adic absolute value and metric to ${\bf Q}_p$ are related as
in (\ref{d_p(x, y) = |x - y|_p}).  Note that $|\cdot|_p$
is an ultrametric absolute value function on ${\bf Q}_p$, and
that $d_p(\cdot, \cdot)$ is an ultrametric on ${\bf Q}_p$, because of
the corresponding properties on ${\bf Q}$.  The possible values of the
$p$-adic absolute value and metric on ${\bf Q}_p$ are $0$ and integer
powers of $p$, as on ${\bf Q}$.  This can be obtained from the
construction of the completion, or from the fact that ${\bf Q}$
is dense in ${\bf Q}_p$.

        The set of \emph{$p$-adic integers}\index{p-adic integers@$p$-adic
integers} is defined by\index{Z_p@${\bf Z}_p$}
\begin{equation}
\label{{bf Z}_p = {x in {bf Q}_p : |x|_p le 1}}
        {\bf Z}_p = \{x \in {\bf Q}_p : |x|_p \le 1\},
\end{equation}
which is a closed set in ${\bf Q}_p$ with respect to the topology
determined by the $p$-adic metric.  Note that the set ${\bf
  Z}$\index{Z@${\bf Z}$} of ordinary integers is contained in ${\bf
  Z}_p$, by definition of the $p$-adic metric on ${\bf Q}$.  It
follows that the closure of ${\bf Z}$ in ${\bf Q}_p$ is contained in
${\bf Z}_p$, and in fact ${\bf Z}_p$ is equal to the closure of ${\bf
  Z}$ in ${\bf Q}_p$.  To see this, let $y \in {\bf Z}_p$ be given,
and remember that $y$ can be approximated by elements of ${\bf Q}$
with respect to the $p$-adic metric, since ${\bf Q}$ is dense in ${\bf
  Q}_p$.  If $w \in {\bf Q}$ satisfies $|y - w|_p \le 1$, then $|w|_p
\le 1$, by the ultrametric version of the triangle inequality.  This
implies that $w$ can be expressed as $a/b$ for some $a, b \in {\bf
  Z}$, where $b \ne 0$ and $b$ is not divisible by $p$, by the
definition of the $p$-adic absolute value on ${\bf Q}$.  Because the
integers modulo $p$ form a field, there is a $c \in {\bf Z}$ such that
$b \, c = 1 - p \, z$ for some $z \in {\bf Z}$.  Thus
\begin{equation}
\label{w = frac{a}{b} = frac{a c}{b c} = frac{a c}{1 - p z}}
 w = \frac{a}{b} = \frac{a \, c}{b \, c} = \frac{a \, c}{1 - p \, z}.
\end{equation}
Of course, $|p \, z|_p = (1/p) \, |z|_p \le 1/p < 1$, and so we can
apply (\ref{lim_{n to infty} sum_{j = 0}^n x^j = frac{1}{1 - x}}) with
$x = p \, z$ to get that
\begin{equation}
\label{w = lim_{n to infty} a c sum_{j = 0}^n p^j z^j}
        w = \lim_{n \to \infty} a \, c \, \sum_{j = 0}^n p^j \, z^j,
\end{equation}
where the limit is taken with respect to the $p$-adic metric.  This
shows that $w$ can be approximated by integers with respect to the
$p$-adic metric when $w \in {\bf Q}$ and $|w|_p \le 1$.  It follows
that every $y \in {\bf Z}_p$ can be approximated by integers with
respect to the $p$-adic metric, since $y$ can be approximated by $w
\in {\bf Q}$ with $|w|_p \le 1$, as before.

        It is easy to see that ${\bf Z}_p$ is a subgroup of ${\bf Q}_p$
with respect to addition, because of the ultrametric version of the
triangle inequality.  Similarly,
\begin{equation}
\label{p^j {bf Z}_p = {p^j x : x in {bf Z}_p} = ...}
        p^j \, {\bf Z}_p = \{p^j \, x : x \in {\bf Z}_p\}
                         = \{y \in {\bf Q}_p : |y| \le p^{-j}\}
\end{equation}
is a subgroup of ${\bf Q}_p$ with respect to addition for every $j \in
{\bf Z}$.  One can also check that ${\bf Z}_p$ is a subring of ${\bf
  Q}_p$, and that $p^j \, {\bf Z}_p$ is an ideal in ${\bf Z}_p$ when
$j \ge 0$.

        Thus the quotient
\begin{equation}
\label{{bf Z}_p / p^j {bf Z}_p}
        {\bf Z}_p / p^j \, {\bf Z}_p
\end{equation}
is defined as a commutative ring for every nonnegative integer $j$.
The natural inclusion of ${\bf Z}$ into ${\bf Z}_p$ may be considered
as a ring homomorphism, which leads to a ring homomorphism from ${\bf
  Z}$ into (\ref{{bf Z}_p / p^j {bf Z}_p}), by composition with the
quotient homomorphism from ${\bf Z}_p$ onto (\ref{{bf Z}_p / p^j {bf
    Z}_p}).  The kernel of this homomorphism from ${\bf Z}$ into
(\ref{{bf Z}_p / p^j {bf Z}_p}) is equal to
\begin{equation}
\label{{bf Z} cap (p^j {bf Z}_p) = p^j {bf Z}}
        {\bf Z} \cap (p^j \, {\bf Z}_p) = p^j \, {\bf Z},
\end{equation}
using the definition of the $p$-adic absolute value on ${\bf Z}$ in
the second step.  Hence the homomorphism from ${\bf Z}$ into (\ref{{bf
    Z}_p / p^j {bf Z}_p}) leads to an injective ring homomorphism from
\begin{equation}
\label{{bf Z} / p^j {bf Z}}
        {\bf Z} / p^j \, {\bf Z}
\end{equation}
into (\ref{{bf Z}_p / p^j {bf Z}_p}).  The usual homomorphism from
${\bf Z}$ into (\ref{{bf Z}_p / p^j {bf Z}_p}) is actually surjective,
because ${\bf Z}$ is dense in ${\bf Z}_p$ with respect to the $p$-adic
metric.  This implies that the we get an ring isomorphism from
(\ref{{bf Z} / p^j {bf Z}}) onto (\ref{{bf Z}_p / p^j {bf Z}_p})
for each $j \ge 0$.  In particular, (\ref{{bf Z}_p / p^j {bf Z}_p})
has exactly $p^j$ elements for each $j \ge 0$.

        It follows that for each nonnegative integer $j$, ${\bf Z}_p$
can be expressed as the union of $p^j$ pairwise-disjoint translates of
$p^j \, {\bf Z}_p$.  Of course, the translates of $p^j \, {\bf Z}_p$
in ${\bf Q}_p$ are the same as closed balls of radius $p^{-j}$ with
respect to the $p$-adic metric.  This implies that ${\bf Z}_p$ is totally
bounded in ${\bf Q}_p$, since ${\bf Z}_p$ can be covered by finitely
many ball of arbitrarily small radius.  It is well known that a subset
of a complete metric space is compact if and only if it is closed and
totally bounded.  This shows that ${\bf Z}_p$ is compact in ${\bf Q}_p$,
because ${\bf Z}_p$ is closed and totally bounded in ${\bf Q}_p$, and
${\bf Q}_p$ is complete by construction.

        An analogous argument implies that $p^l \, {\bf Z}_p$ is
compact in ${\bf Q}_p$ for every integer $l$.  This can also be
obtained from the compactness of ${\bf Z}_p$ and continuity of
multiplication on ${\bf Q}_p$.  Similarly, one can use continuity of
translations on ${\bf Q}_p$ to get that every closed ball in ${\bf
  Q}_p$ is compact.  It follows that closed and bounded subsets of
${\bf Q}_p$ are compact, since closed subsets of compact sets are
compact.  More precisely, it suffices to use the compactness of
$p^l \, {\bf Z}_p$ for each $l \in {\bf Z}$, because every bounded
subset of ${\bf Q}_p$ is contained in $p^l \, {\bf Z}_p$ for some $l$.

\section{$q$-Norms}
\label{q-norms}
\setcounter{equation}{0}

        Let $k$ be a field, and let $V$ be vector space over $k$.
Also let $|\cdot|$ be a $q$-absolute value function on $k$ for some
positive real number $q$.  A nonnegative real-valued function $N$ on
$V$ is said to be a \emph{$q$-norm}\index{qnorms@$q$-norms} on $V$ if
it satisfies the following three conditions.  First, for every $v \in V$,
\begin{equation}
\label{N(v) = 0 if and only if v = 0}
        N(v) = 0 \quad\hbox{if and only if}\quad v = 0.
\end{equation}
Second,
\begin{equation}
\label{N(t v) = |t| N(v) for every t in k and v in V}
 N(t \, v) = |t| \, N(v) \quad\hbox{for every } t \in k \hbox{ and } v \in V.
\end{equation}
Third,
\begin{equation}
\label{N(v + w)^q le N(v)^q + N(w)^q for every v, w in V}
        N(v + w)^q \le N(v)^q + N(w)^q \quad\hbox{for every } v, w \in V.
\end{equation}
If $q = 1$, then we may simply say that $N$ is a
\emph{norm}\index{norms} on $V$.

        Remember that $|\cdot|$ is a $q$-absolute value function on $k$
if and only if $|x|^q$ is an absolute value function on $k$.  In this
case, $N(v)$ is a $q$-norm on $V$ with respect to $|x|$ on $k$ if and
only if $N(v)^q$ is a norm on $V$ with respect to $|x|^q$ on $k$.

        As usual, (\ref{N(v + w)^q le N(v)^q + N(w)^q for every v, w in V})
can be reformulated as saying that
\begin{equation}
\label{N(v + w) le (N(v)^q + N(w)^q)^{1/q} for every v, w in V}
        N(v + w) \le (N(v)^q + N(w)^q)^{1/q} \quad\hbox{for every } v, w \in V.
\end{equation}
We have seen that the right side of this inequality decreases
monotonically in $q$, as in (\ref{(a^{q_2} + b^{q_2})^{1/q_2} le
  (a^{q_1} + b^{q_1})^{1/q_1}}).  If $0 < q_1 < q_2 < \infty$ and
$|\cdot|$ is a $q_2$-absolute value function on $k$, then $|\cdot|$ is
a $q_1$-absolute value function on $k$ too, as in Section
\ref{absolute value functions}.  If we suppose in addition that $N$ is
a $q_2$-norm on $V$, then it follows that $N$ is a $q_1$-norm on $V$
as well.

        Suppose for the moment that $|\cdot|$ is a nonnegative real-valued
function on $k$, and that $N$ is a nonnegative real-valued function on
$V$ that satisfies (\ref{N(v) = 0 if and only if v = 0}), (\ref{N(t v)
  = |t| N(v) for every t in k and v in V}), and (\ref{N(v + w)^q le
  N(v)^q + N(w)^q for every v, w in V}) for some $q > 0$.  If $V \ne
\{0\}$, then one can check that $|\cdot|$ has to be a $q$-absolute
value function on $k$ under these conditions.  Of course, if $|\cdot|$
is a $q$-absolute value function on $k$, then $|\cdot|$ may also be
considered as a $q$-norm on $k$, where $k$ is considered as a
one-dimensional vector space over itself.

        Suppose now that $|\cdot|$ is an ultrametric absolute value function
on $k$.  A nonnegative real-valued function $N$ on $V$ is said to be
an \emph{ultranorm}\index{ultranorms} if it satisfies (\ref{N(v) = 0
  if and only if v = 0}), (\ref{N(t v) = |t| N(v) for every t in k and
  v in V}), and
\begin{equation}
\label{N(v + w) le max(N(v), N(w)) for every v, w in V}
        N(v + w) \le \max(N(v), N(w)) \quad\hbox{for every } v, w \in V.
\end{equation}
As usual, for each $q > 0$, (\ref{N(v + w) le max(N(v), N(w)) for
  every v, w in V}) is equivalent to asking that
\begin{equation}
\label{N(v + w)^q le max(N(v)^q, N(w)^q) for every v, w in V}
 N(v + w)^q \le \max(N(v)^q, N(w)^q) \quad\hbox{for every } v, w \in V.
\end{equation}
If $N$ is an ultranorm on $V$, then it follows that $N$ is a $q$-norm
on $V$ for every $q > 0$, because (\ref{N(v + w)^q le max(N(v)^q,
  N(w)^q) for every v, w in V}) implies (\ref{N(v + w)^q le N(v)^q +
  N(w)^q for every v, w in V}).  This also uses the fact that
$|\cdot|$ is a $q$-absolute value function on $k$ for every $q > 0$
when $|\cdot|$ is an ultrametric absolute value function on $k$,
as in Section \ref{absolute value functions}.

        Similarly, if $|x|$ is an ultrametric absolute value function
on $k$, then $|x|^q$ is an ultrametric absolute value function on $k$
for every $q > 0$, as in Section \ref{absolute value functions}.  If
$N(v)$ is an ultranorm on $V$ with respect to $|x|$ on $k$, then
$N(v)^q$ is an ultranorm on $V$ with respect to $|x|^q$ on $k$ for
every $q > 0$, by (\ref{N(v + w)^q le max(N(v)^q, N(w)^q) for every v,
  w in V}).

        Suppose for the moment again that $|\cdot|$ is a nonnegative
real-valued function on $k$, and that $N$ is a nonnegative real-valued
function on $V$ that satisfies (\ref{N(v) = 0 if and only if v = 0}),
(\ref{N(t v) = |t| N(v) for every t in k and v in V}), and (\ref{N(v +
  w) le max(N(v), N(w)) for every v, w in V}).  If $V \ne \{0\}$, then
one can check that $|\cdot|$ has to be an ultrametric absolute value
function on $k$, as before.  If $|\cdot|$ is an ultrametric absolute
value function on $k$, then $|\cdot|$ may also be considered as an
ultranorm on $k$, as a one-dimensional vector space over itself.

        As in previous situations, (\ref{N(v + w) le max(N(v), N(w)) 
for every v, w in V}) corresponds to taking the limit as $q \to \infty$
in (\ref{N(v + w)^q le N(v)^q + N(w)^q for every v, w in V}), because
of (\ref{lim_{q to infty} (a^q + b^q)^{1/q} = max(a, b)}).  Thus an
ultranorm may be considered as a $q$-norm with $q = \infty$.

        If $|\cdot|$ is a $q$-absolute value function on $k$,
and if $N$ is a $q$-norm on $V$ with respect to $|\cdot|$, then
\begin{equation}
\label{d(v, w) = N(v - w)}
        d(v, w) = N(v - w)
\end{equation}
defines a $q$-metric on $V$.  Similarly, if $|\cdot|$ is an
ultrametric absolute value function on $k$, and if $N$ is an ultranorm
on $V$, then (\ref{d(v, w) = N(v - w)}) is an ultrametric on $V$.

        Consider the function $N$ defined on $V$ by $N(v) = 1$ when
$v \ne 0$, and $N(0) = 0$.  This is an ultranorm on $V$ with respect to
the trivial absolute value function on $k$, which is known as the
\emph{trivial ultranorm}\index{trivial ultranorm} on $V$.  The
ultrametric on $V$ associated to the trivial ultranorm as in
(\ref{d(v, w) = N(v - w)}) is the same as the discrete metric on $V$.

\section{Supremum metrics and norms}
\label{supremum metrics, norms}
\setcounter{equation}{0}

        Let $X$ and $M$ be nonempty sets, and let $d(\cdot, \cdot)$ be
a $q$-metric on $M$ for some $q > 0$.  As usual, a subset of $M$ is
said to be \emph{bounded}\index{bounded sets} with respect to
$d(\cdot, \cdot)$ if it is contained in a ball of finite radius in
$M$.  Similarly, a function $f$ on $X$ with values in $M$ is said to
be bounded\index{bounded mappings} if $f(X)$ is a bounded set in $M$.
Let $B(X, M)$\index{B(X, M)@$B(X, M)$} be the space of bounded
functions on $X$ with values in $M$.  If $f, g \in B(X, M)$, then
$d(f(x), g(x))$ is a bounded nonnegative real-valued function on $X$,
so that
\begin{equation}
\label{sup_{x in X} d(f(x), g(x))}
        \sup_{x \in X} d(f(x), g(x))
\end{equation}
is defined as a nonnegative real number.  It is easy to see that
(\ref{sup_{x in X} d(f(x), g(x))}) defines a $q$-metric on $B(X, M)$,
which may be described as the \emph{supremum
  $q$-metric}.\index{supremum q-metric@supremum $q$-metric} If
$d(\cdot, \cdot)$ is an ultrametric on $M$, then (\ref{sup_{x in X}
  d(f(x), g(x))}) is an ultrametric on $B(X, M)$, which corresponds to
the previous statement with $q = \infty$.

        One can define Cauchy sequences\index{Cauchy sequences} in $M$
with respect to $d(\cdot, \cdot)$ in the same way as for a metric.  If
every Cauchy sequence of elements of $M$ converges to an element of
$M$ with respect to the topology determined by $d(\cdot, \cdot)$, then
we say that $M$ is complete\index{completeness} with respect to
$d(\cdot, \cdot)$, as usual.  Any positive power of $d(\cdot, \cdot)$
determines the same collection of Cauchy sequences in $M$, and leads
to an equivalent version of completeness.  In particular, this permits
one to reduce to the case of ordinary metrics, using suitable powers
of $d(\cdot, \cdot)$.  If $M$ is complete with respect to $d(\cdot,
\cdot)$, then $B(X, M)$ is complete with respect to (\ref{sup_{x in X}
  d(f(x), g(x))}), by standard arguments.

        Suppose for the moment that $X$ is a topological space, and
let $C(X, M)$\index{C(X, M)@$C(X, M)$} be the space of continuous
mappings from $X$ into $M$.  Also let
\begin{equation}
\label{C_b(X, M) = B(X, M) cap C(X, M)}
        C_b(X, M) = B(X, M) \cap C(X, M)
\end{equation}
be the space of bounded continuous mappings from $X$ into
$M$.\index{C_b(X, M)@$C_b(X, M)$} It is easy to see that $C_b(X, M)$
is a closed set in $B(X, M)$ with respect to the supremum $q$-metric,
by standard arguments.  If $M$ is complete with respect to $d(\cdot,
\cdot)$, then it follows that $C_b(X, M)$ is complete with respect to
the supremum $q$-metric.  Note that compact subsets of $M$ are
bounded, and hence that continuous mappings from $X$ into $M$ are
bounded when $X$ is compact.

        Let $k$ be a field, and let $|\cdot|$ be a $q$-absolute
value function on $k$ for some $q > 0$.  Also let $V$ be a vector
space over $k$, and let $N$ be a $q$-norm on $V$ with respect to
$|\cdot|$ on $k$.  Thus (\ref{d(v, w) = N(v - w)}) defines a
$q$-metric on $V$, as in the previous section.  If $X$ is a nonempty
set again, then we shall also use the notation $\ell^\infty(X,
V)$\index{l^\infty(X, V)@$\ell^\infty(X, V)$} for the space of bounded
$V$-valued functions on $X$.  It is easy to see that this is a vector
space over $k$ with respect to pointwise addition and scalar
multiplication.  Put
\begin{equation}
\label{||f||_infty = ||f||_{ell^infty(X, V)} = sup_{x in X} N(f(x))}
        \|f\|_\infty = \|f\|_{\ell^\infty(X, V)} = \sup_{x \in X} N(f(x))
\end{equation}
for each $f \in \ell^\infty(X, V)$, which defines a $q$-norm on
$\ell^\infty(X, V)$ with respect to $|\cdot|$ on $k$.  This is the
\emph{supremum $q$-norm}\index{supremum q-norm@supremum $q$-norm} on
$\ell^\infty(X, V)$ corresponding to $N$ on $V$.  By construction, the
$q$-metric on $\ell^\infty(X, V)$ associated to (\ref{||f||_infty =
  ||f||_{ell^infty(X, V)} = sup_{x in X} N(f(x))}) is the supremum
$q$-metric that corresponds to the $q$-metric (\ref{d(v, w) = N(v -
  w)}) on $V$ associated to $N$.  If $|\cdot|$ is an ultrametric
absolute value function on $k$, and if $N$ is an ultranorm on $V$,
then (\ref{||f||_infty = ||f||_{ell^infty(X, V)} = sup_{x in X}
  N(f(x))}) is an ultranorm on $\ell^\infty(X, V)$ as well.  If $X$ is
a topological space, then $C(X, V)$ is a vector space over $k$ with
respect to pointwise addition and scalar multiplication too, and
$C_b(X, V)$ is a linear subspace of $\ell^\infty(X, V)$.  Of course,
if $X$ is equipped with the discrete topology, then every function
on $X$ is continuous, so that $C_b(X, V)$ is the same as $\ell^\infty(X, V)$.

\section{Summable functions}
\label{summable functions}
\setcounter{equation}{0}

        Let $X$ be a nonempty set, and let $f$ be a nonnegative
real-valued function on $X$.  The sum
\begin{equation}
\label{sum_{x in X} f(x)}
        \sum_{x \in X} f(x)
\end{equation}
is defined as a nonnegative extended real number to be the supremum
of the sums
\begin{equation}
\label{sum_{x in A} f(x)}
        \sum_{x \in A} f(x)
\end{equation}
over all nonempty finite subsets $A$ of $X$.  If $g$ is another
nonnegative real-valued function on $X$ and $a$ is a positive real
number, then one can check that
\begin{equation}
\label{sum_{x in X} (f(x) + g(x)) = sum_{x in X} f(x) + sum_{x in X} g(x)}
        \sum_{x \in X} (f(x) + g(x)) = \sum_{x \in X} f(x) + \sum_{x \in X} g(x)
\end{equation}
and
\begin{equation}
\label{sum_{x in X} a f(x) = a sum_{x in X} f(x)}
        \sum_{x \in X} a \, f(x) = a \, \sum_{x \in X} f(x),
\end{equation}
with the usual interpretations for nonnegative extended real numbers.
If (\ref{sum_{x in X} f(x)}) is finite, then $f$ is said to be
\emph{summable}\index{summable functions} on $X$.  If $f$ and $g$ are
summable on $X$, then it follows that $f + g$ is summable on $X$, and
that $a \, f$ is summable on $X$ for every $a \ge 0$.

        Similarly, $f$ is said to be \emph{$r$-summable}\index{r-summable
functions@$r$-summable functions} on $X$ for some positive real number
$r$ if $f(x)^r$ is summable on $X$.  Put
\begin{equation}
\label{||f||_r = (sum_{x in X} f(x)^r)^{1/r}}
        \|f\|_r = \Big(\sum_{x \in X} f(x)^r\Big)^{1/r}
\end{equation}
when $0 < r < \infty$, and
\begin{equation}
\label{||f||_infty = sup_{x in X} f(x)}
        \|f\|_\infty = \sup_{x \in X} f(x).
\end{equation}
Thus (\ref{||f||_r = (sum_{x in X} f(x)^r)^{1/r}}) is finite exactly
when $f$ is $r$-summable on $X$, and (\ref{||f||_infty = sup_{x in X}
  f(x)}) is finite exactly when $f$ is bounded on $X$.  If $f$ is
bounded on $X$, then (\ref{||f||_infty = sup_{x in X} f(x)}) is the
same as the supremum norm of $f$, with respect to the standard
absolute value function on ${\bf R}$.  Note that
\begin{equation}
\label{||a f||_r = a ||f||_r}
        \|a \, f\|_r = a \, \|f\|_r
\end{equation}
for every $a, r > 0$, and in particular that $a \, f$ is $r$-summable
on $X$ for every $a \ge 0$ when $f$ is $r$-summable on $X$.

        If $f$ is $r$-summable on $X$ for some $r > 0$, then it is easy
to see that $f$ is bounded on $X$, and that
\begin{equation}
\label{||f||_infty le ||f||_r}
        \|f\|_\infty \le \|f\|_r.
\end{equation}
This implies that for each $t > r$ and $x \in X$, we have that
\begin{equation}
\label{f(x)^t le ||f||_infty^{t - r} f(x)^r le ||f||_r^{t - r} f(x)^r}
 f(x)^t \le \|f\|_\infty^{t - r} \, f(x)^r \le \|f\|_r^{t - r} \, f(x)^r.
\end{equation}
Summing over $x \in X$, we get that
\begin{equation}
\label{||f||_t^t = sum_{x in X} f(x)^t le ... = ||f||_r^t}
 \|f\|_t^t = \sum_{x \in X} f(x)^t \le \|f\|_r^{t - r} \, \sum_{x \in X} f(x)^r
                  = \|f\|_r^{t - r} \, \|f\|_r^r = \|f\|_r^t,
\end{equation}
and hence
\begin{equation}
\label{||f||_t le ||f||_r}
        \|f\|_t \le \|f\|_r.
\end{equation}
In particular, $f$ is $t$-summable on $X$ for every $t > r$ when $f$
is $r$-summable on $X$.

        Let $g$ be another nonnegative real-valued function on $X$ again,
and observe that
\begin{equation}
\label{(f(x) + g(x))^r le (2 max(f(x), g(x)))^r le 2^r (f(x)^r + g(x)^r)}
 (f(x) + g(x))^r \le (2 \, \max(f(x), g(x)))^r \le 2^r \, (f(x)^r + g(x)^r)
\end{equation}
for every $x \in X$ and $r > 0$.  If $f$ and $g$ are both $r$-summable
on $X$, then it follows that $f + g$ is $r$-summable on $X$ too, by
summing over $x \in X$.  More precisely, if $0 < r \le 1$, then we
have that
\begin{equation}
\label{(f(x) + g(x))^r le f(x)^r + g(x)^r}
        (f(x) + g(x))^r \le f(x)^r + g(x)^r
\end{equation}
for every $x \in X$.  This follows from (\ref{a^{q_2} + b^{q_2} le
  ... = (a^{q_1} + b^{q_1})^{q_2/q_1}}), with $q_1 = r$ and $q_2 = 1$,
and it can also be derived from (\ref{||f||_t le ||f||_r}), with $t =
1$.  Summing both sides of (\ref{(f(x) + g(x))^r le f(x)^r + g(x)^r})
over $x \in X$, we get that
\begin{equation}
\label{||f + g||_r^r le ||f||_r^r + ||g||_r^r}
        \|f + g\|_r^r \le \|f\|_r^r + \|g\|_r^r
\end{equation}
when $0 < r \le 1$.  If $r \ge 1$, then we have that
\begin{equation}
\label{||f + g||_r le ||f||_r + ||g||_r}
        \|f + g\|_r \le \|f\|_r + \|g\|_r,
\end{equation}
by Minkowski's inequality\index{Minkowski's inequality} for sums.  Of
course, (\ref{||f + g||_r^r le ||f||_r^r + ||g||_r^r}) and (\ref{||f +
  g||_r le ||f||_r + ||g||_r}) reduce to (\ref{sum_{x in X} (f(x) +
  g(x)) = sum_{x in X} f(x) + sum_{x in X} g(x)}) when $r = 1$, and it
is easy to verify (\ref{||f + g||_r le ||f||_r + ||g||_r}) directly
when $r = \infty$.

\section{$\ell^r$ Norms}
\label{l^r norms}
\setcounter{equation}{0}

        Let $k$ be a field, and let $|\cdot|$ be a $q$-absolute
value function on $k$ for some $q > 0$.  Also let $V$ be a vector
space over $k$, and let $N$ be a $q$-norm on $V$ with respect to
$|\cdot|$ on $k$.  A $V$-valued function $f$ on $X$ is said to be
\emph{$r$-summable}\index{r-summable functions@$r$-summable functions}
on a nonempty set $X$ for some positive real number $r$ if $N(f(x))$
is $r$-summable as a nonnegative real-valued function on $X$, as in
the previous section.  If $f$ is $r$-summable with $r = 1$, then we
may simply say that $f$ is summable on $X$.\index{summable functions}
The space of $V$-valued $r$-summable functions on $X$ is denoted
$\ell^r(X, V)$.\index{l^r(X, V)@$\ell^r(X, V)$}

        Let $f$ and $g$ be $V$-valued functions on $X$, and observe that
\begin{equation}
\label{N(f(x) + g(x))^r le (N(f(x))^q + N(g(x))^q)^{r/q}}
        N(f(x) + g(x))^r \le (N(f(x))^q + N(g(x))^q)^{r/q}
\end{equation}
for every $x \in X$, by the $q$-norm version of the triangle
inequality.  Thus
\begin{equation}
\label{N(f(x) + g(x))^r le 2^{r/q} (N(f(x))^r + N(g(x))^r)}
        N(f(x) + g(x))^r \le 2^{r/q} \, (N(f(x))^r + N(g(x))^r)
\end{equation}
for every $x \in X$, as in (\ref{(f(x) + g(x))^r le (2 max(f(x),
  g(x)))^r le 2^r (f(x)^r + g(x)^r)}), but with $r$ replaced by $r /
q$.  If $f$ and $g$ are both $r$-summable on $X$, then it follows that
$f + g$ is $r$-summable too, by summing over $x \in X$.  It is easy to
see that $r$-summability is also preserved by scalar multiplication,
so that $\ell^r(X, V)$ is a vector space with respect to pointwise
addition and scalar multiplication.

        Put
\begin{equation}
\label{||f||_r = ||f||_{ell^r(X, V)} = (sum_{x in X} N(f(x))^r)^{1/r}}
 \|f\|_r = \|f\|_{\ell^r(X, V)} = \Big(\sum_{x \in X} N(f(x))^r\Big)^{1/r}
\end{equation}
for each $f \in \ell^r(X, V)$, which clearly satisfies the usual
positivity and homogeneity requirements of a norm.  If $r \le q$, then
(\ref{N(f(x) + g(x))^r le (N(f(x))^q + N(g(x))^q)^{r/q}}) implies that
\begin{equation}
\label{N(f(x) + g(x))^r le N(f(x))^r + N(g(x))^r}
        N(f(x) + g(x))^r \le N(f(x))^r + N(g(x))^r
\end{equation}
for every $f, g \in \ell^r(X, V)$ and $x \in X$, as in (\ref{(f(x) +
  g(x))^r le f(x)^r + g(x)^r}), with $r$ replaced by $r / q$.  Summing
over $x \in X$, we get that
\begin{equation}
\label{||f + g||_r^r le ||f||_r^r + ||g||_r^r, 2}
        \|f + g\|_r^r \le \|f\|_r^r + \|g\|_r^r
\end{equation}
for every $f, g \in \ell^r(X, V)$, so that $\|f\|_r$ defines an
$r$-norm on $\ell^r(X, V)$ when $r \le q$.  If $q \le r$, then
(\ref{N(f(x) + g(x))^r le (N(f(x))^q + N(g(x))^q)^{r/q}}) implies that
\begin{equation}
\label{||f + g||_r^q le ||f||_r^q + ||g||_r^q}
        \|f + g\|_r^q \le \|f\|_r^q + \|g\|_r^q
\end{equation}
for every $f, g \in \ell^r(X, V)$, using (\ref{||f + g||_r le ||f||_r
  + ||g||_r}) with $r$ replaced by $r/q$.  This shows that $\|f\|_r$
is a $q$-norm on $\ell^r(X, V)$ when $q \le r$.

        If $N$ is an ultranorm on $V$, then we have that
\begin{equation}
\label{N(f(x) + g(x))^r le max(N(f(x)), N(g(x)))^r le N(f(x))^r + N(g(x))^r}
 \qquad N(f(x) + g(x))^r \le \max(N(f(x)), N(g(x)))^r \le N(f(x))^r + N(g(x))^r
\end{equation}
for all $V$-valued functions $f$, $g$ on $X$, $r > 0$, and $x \in X$.
This implies that (\ref{||f + g||_r^r le ||f||_r^r + ||g||_r^r, 2})
holds for every $f, g \in \ell^r(X, V)$ and $r > 0$, by summing over
$x \in X$.  Thus $\|f\|_r$ is an $r$-norm on $\ell^r(X, V)$ for every
$r > 0$ in this case, which corresponds to $q = \infty$ in the
previous discussion.

\section{Infinite series}
\label{infinite series}
\setcounter{equation}{0}

        Let $k$ be a field with a $q$-absolute value function $|\cdot|$
for some $q > 0$.  Also let $V$ be a vector space over $k$ again, and
let $N$ be a $q$-norm on $V$ with respect to $|\cdot|$ on $k$.  This
leads to a $q$-metric $d(v, w)$ on $V$ associated to $N$ as in
(\ref{d(v, w) = N(v - w)}), and hence to a topology on $V$, as in
Section \ref{open, closed balls}.  As usual, an infinite series
$\sum_{j = 1}^\infty a_j$ with terms in $V$ is said to converge in $V$
if the corresponding sequence of partial sums
\begin{equation}
\label{s_n = sum_{j = 1}^n a_j}
        s_n = \sum_{j = 1}^n a_j
\end{equation}
converges to an element of $V$ with respect to this topology, in which
case the value of the sum $\sum_{j = 1}^\infty a_j$ is defined to be the
limit of the sequence $\{s_n\}_{n = 1}^\infty$.

        As in Section \ref{supremum metrics, norms}, one can define
Cauchy sequences and completeness with respect to $d(v, w)$ in
the same way as for ordinary metrics, and this is equivalent to
defining Cauchy sequences with respect to the metric $d(v, w)^q$.
It is easy to see that the sequence (\ref{s_n = sum_{j = 1}^n a_j})
of partial sums is Cauchy sequence in $V$ for each $\epsilon > 0$
there is a positive integer $L$ such that
\begin{equation}
\label{N(sum_{j = l}^n a_j) < epsilon}
        N\Big(\sum_{j = l}^n a_j\Big) < \epsilon
\end{equation}
for every $n \ge l \ge L$.  In particular, this implies that
$\{a_j\}_{j = 1}^\infty$ converges to $0$ in $V$, by taking $l = n$.
Note that
\begin{equation}
\label{N(sum_{j = l}^n a_j)^q le sum_{j = l}^n N(a_j)^q}
        N\Big(\sum_{j = l}^n a_j\Big)^q \le \sum_{j = l}^n N(a_j)^q
\end{equation}
for every $n \ge j \ge 1$, by the $q$-norm version of the triangle
inequality.  Similarly, if $N$ is an ultranorm on $V$, then
\begin{equation}
\label{N(sum_{j = l}^n a_j) le max_{l le j le n} N(a_j)}
        N\Big(\sum_{j = l}^n a_j\Big) \le \max_{l \le j \le n} N(a_j)
\end{equation}
for every $n \ge l \ge 1$.

        Let us say that $\sum_{j = 1}^\infty a_j$ converges
\emph{$q$-absolutely}\index{qabsolute convergence@$q$-absolute convergence}
if
\begin{equation}
\label{sum_{j = 1}^infty N(a_j)^q}
        \sum_{j = 1}^\infty N(a_j)^q
\end{equation}
converges as an infinite series of nonnegative real numbers.  Of
course, this reduces to the usual notion of absolute
convergence\index{absolute convergence} when $q = 1$.  If (\ref{sum_{j
    = 1}^infty N(a_j)^q}) converges, then one can use (\ref{N(sum_{j =
    l}^n a_j)^q le sum_{j = l}^n N(a_j)^q}) to check that $\sum_{j =
  1}^\infty a_j$ satisfies the Cauchy criterion described in the
preceding paragraph, as in the $q = 1$ case.  If $V$ is complete, then
$\sum_{j = 1}^\infty a_j$ converges in $V$, and we have that
\begin{equation}
\label{N(sum_{j = 1}^infty a_j)^q le sum_{j = 1}^infty N(a_j)^q}
 N\Big(\sum_{j = 1}^\infty a_j\Big)^q \le \sum_{j = 1}^\infty N(a_j)^q,
\end{equation}
by standard arguments.  Similarly, if $N$ is an ultranorm on $V$, and
if $\{a_j\}_{j = 1}^\infty$ converges to $0$ in $V$, then
(\ref{N(sum_{j = l}^n a_j) le max_{l le j le n} N(a_j)}) implies that
$\sum_{j = 1}^\infty a_j$ satisfies the Cauchy criterion.  If $V$ is
complete, then it follows that $\sum_{j = 1}^\infty a_j$ converges in
$V$, and that
\begin{equation}
\label{N(sum_{j = 1}^infty a_j) le max_{j ge 1} N(a_j)}
        N\Big(\sum_{j = 1}^\infty a_j\Big) \le \max_{j \ge 1} N(a_j).
\end{equation}
Note that the maximum on the right side of (\ref{N(sum_{j = 1}^infty
  a_j) le max_{j ge 1} N(a_j)}) is attained under these conditions,
because $N(a_j) \to 0$ as $j \to \infty$.

        If every $q$-absolutely convergent series in $V$ converges to
an element of $V$, then a well-known argument imples that $V$ has to
be complete.  To see this, let $\{v_j\}_{j = 1}^\infty$ be any Cauchy
sequence of elements of $V$.  It is easy to see that there is a
subsequence $\{v_{j_l}\}_{l = 1}^\infty$ of $\{v_j\}_{j = 1}^\infty$
such that
\begin{equation}
\label{N(v_{j_l} - v_{j_{l + 1}}) < 2^{-l}}
        N(v_{j_l} - v_{j_{l + 1}}) < 2^{-l}
\end{equation}
for each $l \ge 1$.  This implies that $\sum_{l = 1}^\infty (v_{j_l} -
v_{j_{l + 1}})$ converges $q$-absolutely, and hence that this series
converges in $V$, by hypothesis.  Of course,
\begin{equation}
\label{sum_{l = 1}^n (v_{j_l} - v_{j_{l + 1}}) = v_{j_1} - v_{j_{n + 1}}}
        \sum_{l = 1}^n (v_{j_l} - v_{j_{l + 1}}) = v_{j_1} - v_{j_{n + 1}}
\end{equation}
for each positive integer $n$, so that $\sum_{l = 1}^\infty (v_{j_l} -
v_{j_{l + 1}})$ coverges in $V$ if and only if $\{v_{j_n}\}_{n =
  1}^\infty$ converges as a sequence in $V$.  Because $\{v_j\}_{j =
  1}^\infty$ is a Cauchy sequence in $V$, the convergence of a
subsequence $\{v_{j_n}\}_{n = 1}^\infty$ in $V$ implies that
$\{v_j\}_{j = 1}^\infty$ converges to the same limit, as desired.  If
$N$ is an ultranorm on $V$, then one can consider infinite series
$\sum_{j = 1}^\infty a_j$ with terms in $V$ such that $N(a_j) \to 0$
as $j \to \infty$, as the analogue of $q$-absolute convergence with $q
= \infty$.  If every such series converges in $V$, then $V$ has to be
complete, as before.  In this case, if $\{v_j\}_{j = 1}^\infty$ is any
Cauchy sequence in $V$, then one can apply the hypothesis on infinite
series directly to $\sum_{j = 1}^\infty (v_j - v_{j + 1})$.

\section{Vanishing at infinity}
\label{vanishing at infinity}
\setcounter{equation}{0}

        Let $k$ be a field with a $q$-absolute value function $|\cdot|$
for some $q > 0$ again, let $V$ be a vector space over $k$, and let
$N$ be a $q$-norm on $V$ with respect to $|\cdot|$.  Also let $X$ be a
(nonempty) locally compact Hausdorff topological space.  A continuous
$V$-valued function $f$ on $X$ is said to \emph{vanish at
  infinity}\index{vanishing at infinity} if for each $\epsilon > 0$
there is a compact set $K(\epsilon) \subseteq X$ such that
\begin{equation}
\label{N(f(x)) < epsilon}
        N(f(x)) < \epsilon
\end{equation}
for every $x \in X \setminus K(\epsilon)$.  This is equivalent to saying
that
\begin{equation}
\label{{x in X : N(f(x)) ge epsilon}}
        \{x \in X : N(f(x)) \ge \epsilon\}
\end{equation}
is a compact subset of $X$ for each $\epsilon > 0$.  More precisely,
if (\ref{{x in X : N(f(x)) ge epsilon}}) is compact for some $\epsilon
> 0$, then one can simply take $K(\epsilon)$ to be (\ref{{x in X :
    N(f(x)) ge epsilon}}).  Conversely, if $K(\epsilon)$ is a compact
subset of $X$ such that (\ref{N(f(x)) < epsilon}) holds for every $x
\in X \setminus K(\epsilon)$, then (\ref{{x in X : N(f(x)) ge
    epsilon}}) is contained in $K(\epsilon)$.  If $f$ is continuous,
then (\ref{{x in X : N(f(x)) ge epsilon}}) is a closed set in $X$,
since it is the inverse image of a closed ball in $V$.  This implies
that (\ref{{x in X : N(f(x)) ge epsilon}}) is compact, because closed
subsets of compact sets are compact.

        The space of continuous $V$-valued functions on $X$ that vanish
at infinity is denoted $C_0(X, V)$.\index{C_0(X, V)@$C_0(X, V)$} It is
easy to see that
\begin{equation}
\label{C_0(X, V) subseteq C_b(X, V)}
        C_0(X, V) \subseteq C_b(X, V),
\end{equation}
by taking $\epsilon = 1$ in the previous definition, and using the
fact that continuous functions are bounded on compact sets.  Moreover,
$C_0(X, V)$ is a linear subspace of $C_b(X, V)$, as a vector space
with respect to pointwise addition and scalar multiplication.  One can
also check that $C_0(X, V)$ is a closed set in $C_b(X, V)$, with
respect to the topology determined by the supremum $q$-norm.

        If $f$ is any $V$-valued function on $X$, then the
\emph{support}\index{support of a function} is denoted $\supp f$,
and defined to be the closure in $X$ of the set of $x \in X$ such that
$f(x) \ne 0$.  The space of continuous $V$-valued functions with
compact support in $X$ may be denoted
$C_{com}(X, V)$\index{C_com(X, V)@$C_{com}(X, V)$} or
$C_{00}(X, V)$,\index{C_00(X, V)@$C_{00}(X, V)$} and is a linear
subspace of $C_0(X, V)$.  If $X$ is equipped with the discrete topology,
so that every function on $X$ is continuous, then $C_0(X, V)$ may also
be denoted $c_0(X, V)$,\index{c_0(X, V)@$c_0(X, V)$} and $C_{00}(X, V)$
may be denoted $c_{00}(X, V)$.\index{c_00(X, V)@$c_{00}(X, V)$}
In this case, the support of a $V$-valued function $f$ on $X$ is
simply the set of $x \in X$ such that $f(x) \ne 0$, and the only
compact subsets of $X$ are those with only finitely many elements.
Thus $c_{00}(X, V)$ consists of the $V$-valued functions $f$ on $X$
such that $f(x) = 0$ for all but finitely many $x \in X$, and $c_0(X, V)$
consists of the $V$-valued functions $f$ on $V$ such that for each
$\epsilon > 0$, (\ref{N(f(x)) < epsilon}) holds for all but finitely many
$x \in X$.

        Let $X$ be a locally compact Hausdorff topological space again.
If $K \subseteq X$ is compact, $U \subseteq X$ is an open set, and $K
\subseteq U$, then it is well known that there is a continuous
real-valued function on $X$ with compact support contained in $U$
which is equal to $1$ on $K$, and which takes values between $0$ and
$1$ on all of $X$, by Urysohn's lemma.  If $k = {\bf R}$ or ${\bf C}$
equipped with the standard absolute value function, then one can use
this to show that $C_{com}(X, V)$ is dense in $C_0(X, V)$ with respect
to the supremum $q$-norm.  Of course, the same argument can be used
when $k = {\bf R}$ or ${\bf C}$ is equipped with a $q$-absolute value
function which is a power of the standard absolute value function.

        If $K_1$ is a compact open subset of $X$, then the function on
$X$ equal to $1$ on $K_1$ and to $0$ on $X \setminus K_1$ is continuous
and has compact support equal to $K_1$.  If every compact subset of
$X$ is contained in a compact open set, then one can use these
functions to show that $C_{com}(X, V)$ is dense in $C_0(X, V)$ with
respect to the supremum $q$-norm.  More precisely, this works for any
field $k$ with a $q$-absolute value function, and for any vector space
$V$ over $k$ with a $q$-norm $N$.  In particular, this condition holds
when $X$ is equipped with the discrete topology.  Note that this
condition also holds when $X$ is locally compact and has topological
dimension $0$, as in Section \ref{dimension 0}.

        Suppose for the moment that $|\cdot|$ is an ultrametric absolute
value function on $k$, and that $N$ is an ultranorm on $V$.  Thus
\begin{equation}
\label{{v in V : N(v) ge epsilon}}
        \{v \in V : N(v) \ge \epsilon\}
\end{equation}
is an open set in $V$ with respect to the topology determined by the
ultrametric associated to $N$ for every $\epsilon > 0$, as in Section
\ref{open, closed balls}.  If $f : X \to V$ is continuous, then it
follows that (\ref{{x in X : N(f(x)) ge epsilon}}) is an open set in
$X$ for every $\epsilon > 0$, since (\ref{{x in X : N(f(x)) ge
    epsilon}}) is the same as the inverse image of (\ref{{v in V :
    N(v) ge epsilon}}) under $f$.  If $f$ also vanishes at infinity on
$X$, then we have seen that (\ref{{x in X : N(f(x)) ge epsilon}})
is a compact subset of $X$ for every $\epsilon > 0$ too.  Using this
and the remarks at the beginning of the previous paragraph, one
can check that $C_{com}(X, V)$ is dense in $C_0(X, V)$ with respect
to the supremum norm in this case as well.

        Let $k$ be any field with a $q$-absolute value function $|\cdot|$
again, and let $V$ be a vector space with a $q$-norm $N$.  Also let
$X$ be a nonempty set, which may be considered as being equipped with
the discrete topology, and let $r$ be a positive real number.  Observe
that every $V$-valued function on $X$ with finite support is
$r$-summable, so that
\begin{equation}
\label{c_{00}(X, V) subseteq ell^r(X, V)}
        c_{00}(X, V) \subseteq \ell^r(X, V).
\end{equation}
If $f \in \ell^r(X, V)$, then
\begin{equation}
\label{sum_{x in X} N(f(x))^r < infty}
        \sum_{x \in X} N(f(x))^r < \infty,
\end{equation}
where the sum is defined as the supremum of the corresponding finite
subsums, as in Section \ref{summable functions}.  Thus for each
$\epsilon > 0$ there should be a finite set $A(\epsilon) \subseteq X$
such that
\begin{equation}
\label{sum_{x in X} N(f(x))^r < sum_{x in A(epsilon)} N(f(x))^r + epsilon}
 \sum_{x \in X} N(f(x))^r < \sum_{x \in A(\epsilon)} N(f(x))^r + \epsilon,
\end{equation}
which implies that
\begin{equation}
\label{sum_{x in X setminus A(epsilon)} N(f(x))^r < epsilon}
        \sum_{x \in X \setminus A(\epsilon)} N(f(x))^r < \epsilon.
\end{equation}
It follows from this that $f$ can be approximated by $V$-valued
functions with finite support in $X$ with respect to the $\ell^r$
norm, so that $c_{00}(X, V)$ is dense in $\ell^r(X, V)$.  In
particular, this argument shows that $f$ vanishes at infinity on $X$,
which implies that
\begin{equation}
\label{ell^r(X, V) subseteq c_0(X, V)}
        \ell^r(X, V) \subseteq c_0(X, V).
\end{equation}

\part{Topological dimension}
\label{topological dimension}

\section{Separation conditions}
\label{separation conditions}
\index{separation conditions}
\setcounter{equation}{0}

        Remember that a topological space $X$ satisfies the
\emph{first separation condition}\index{first separation condition}
if for each $x, y \in X$ with $x \ne y$ there is an open subset of $X$
that contains $x$ and not $y$.  This implies that there is also an
open subset of $X$ that contains $y$ and not $x$, by interchanging the
roles of $x$ and $y$.  Equivalently, $X$ satisfies the first
separation condition if and only if every subset of $X$ with exactly
one element is a closed set, which implies that finite subsets of $X$
are closed sets.  Similarly, $X$ satisfies the \emph{second separation
  condition}\index{second separation condition} if every pair of
distinct elements of $X$ is contained in a pair of disjoint open
subsets of $X$.  This obviously implies that $X$ satisfies the first
separation condition, and topological spaces that satisfy the second
separation condition are said to be \emph{Hausdorff}.\index{Hausdorff
  topological spaces}

        The \emph{$0$th separation condition} asks that for each
pair of distinct elements of $X$ there be an open subset of $X$ that
contains one of the two points and not the other, but without
specifying which of the two points is contained in the open set.  Thus
the first separation condition automatically implies the $0$th
separation condition.  Equivalently, $X$ satisfies the $0$th
separation condition if for every pair of distinct elements of $X$
there is a closed set in $X$ that contains one of the two points and
not the other, without specifying which of the two points is contained
in the closed set.

        If $X$ satisfies the $0$th separation condition, and for each
$x \in X$ and closed set $E \subseteq X$ there are disjoint open sets
$U, V \subseteq X$ such that $p \in U$ and $E \subseteq V$, then $X$
satisfies the \emph{third separation condition},\index{third separation
condition} and is said to be \emph{regular}.\index{regular topological spaces}
It is easy to see that regular topological spaces are Hausdorff, and
in particular that they satisfy the first separation condition.  More
precisely, if $x, y \in X$ and $x \ne y$, then the $0$th separation
condition implies that there is a closed set in $X$ that contains one
of $x$, $y$ and not the other, and one can use the rest of the
regularity condition to show that $x$, $y$ are contained in disjoint
open subsets of $X$.  Sometimes regularity of topological spaces is
defined by including the first separation condition in the definition
instead of the $0$th separation condition, which would be equivalent
by the previous remarks.  Regularity can also be characterized by
asking that $X$ satisfy the $0$th separation condition, and that for
every $x \in X$ and open set $W \subseteq X$ with $x \in W$ there be
an open set $U \subseteq X$ such that $x \in U$ and $\overline{U}
\subseteq W$, where $\overline{U}$ denotes the closure of $U$ in $X$.

        A topological space $X$ is said to be \emph{completely 
Hausdorff}\index{completely Hausdorff spaces} if every pair of distinct
elements of $X$ is contained in a pair of open subsets of $X$ with
disjoint closures in $X$.  This is also known as separation condition
number two and a half.\index{separation conditions} Completely
Hausdorff space are obviously Hausdorff, and regular topological
spaces are completely Hausdorff.

        If $X$ satisfies the first separation condition, and if every
pair of disjoint closed subsets of $X$ are contained in disjoint open
subsets of $X$, then $X$ satisfies the \emph{fourth separation
  condition}.\index{fourth separation condition} This implies that $X$
satisfies the second and third separation conditions, and $X$ is said
to be \emph{normal}\index{normal topological spaces} in this case.
Equivalently, $X$ is normal if $X$ satisfies the first separation
condition, and for every closed set $A \subseteq X$ and open set $W
\subseteq X$ with $A \subseteq W$ there is an open set $U \subseteq X$
such that $A \subseteq U$ and $\overline{U} \subseteq W$.

        Remember that a pair of subsets $A$, $B$ of a topological space
$X$ are said to be \emph{separated}\index{separated sets} in $X$ if
\begin{equation}
\label{overline{A} cap B = A cap overline{B} = emptyset}
        \overline{A} \cap B = A \cap \overline{B} = \emptyset.
\end{equation}
If $X$ satisfies the first separation condition, and if every pair of
separated subsets of $X$ are contained in disjoint open subsets of
$X$, then $X$ satisfies the \emph{fifth separation
  condition},\index{fifth separation condition} and $X$ is said to be
\emph{completely normal}.\index{completely normal spaces} Completely
normal topological spaces are automatically normal, because disjoint
closed sets are obviously separated.  It is well known that metric
spaces are completely normal.

        Let $Y$ be a subset of a topological space $X$, equipped with
the induced topology.  If $X$ satisfies any of the $0$th, first,
second, or third separation conditions, then $Y$ has the same
property.  This also works for completely Hausdorff and completely
normal spaces, but not for normal spaces.  In the case of completely
normal spaces, this uses the fact that a pair of subsets of $Y$ are
separated with respect to the induced topology on $Y$ if and only if
they are separated as subsets of $X$.

        Let $\tau_1$ and $\tau_2$ be topologies on a set $X$ with
$\tau_1 \subseteq \tau_2$.  If $X$ satisfies any of the $0$th, first,
or second separation conditions with respect to $\tau_1$, then $X$ has
the same property with respect to $\tau_2$.  This also works for the
completely Hausdorff condition, but not for regularity.

        If $X$ is a Hausdorff topological space and $K \subseteq X$
is compact, then it is well known that $K$ is a closed set in $X$.
More precisely, if $x \in X \setminus K$, then $x$ and $K$ are
contained in disjoint open subsets of $X$.  To see this, one can
use the Hausdorff condition to cover $K$ by open sets, each of which
is disjoint from an open set that contains $x$, and then use compactness
to reduce to a finite subcovering.  Similarly, one can show that every
pair of disjoint compact subsets of $X$ is contained in a pair of
disjoint open sets.  In particular, this implies that compact Hausdorff
spaces are normal, because closed subsets of compact spaces are compact
as well.  If $X$ is regular, $E \subseteq X$ is a closed set, $K \subseteq X$
is compact, and $E \cap K = \emptyset$, then $E$ and $K$ are contained
in disjoint open subsets of $X$.  This can be obtained by covering $K$
by open sets, each of which is disjoint from an open set that contains
$E$, and using compactness to reduce to a finite subcovering.

        A topological space $X$ is said to be \emph{locally 
compact}\index{locally compact spaces} if for each $x \in X$
there is an open set $W \subseteq X$ and a compact set $K \subseteq X$
such that $x \in W$ and $W \subseteq K$.  If $X$ is also Hausdorff,
then $K$ is a closed set in $X$, so that $\overline{W} \subseteq K$.
This implies that $\overline{W}$ is compact, since closed subsets
of compact sets are compact.  If $X$ is locally compact and
$H \subseteq X$ is compact, then it is easy to see that $H$ is
contained in an open set in $X$ that is contained in another compact
set, by covering $H$ by finitely many open sets that are contained
in compact sets.

        Suppose that $X$ is a Hausdorff topological space, $W \subseteq X$
is an open set, $x \in W$, and $\overline{W}$ is compact.  Thus the
boundary $\partial W = \overline{W} \setminus W$ of $W$ is compact,
and of course $x \not\in \partial W$.  As before, there is an open set
$U \subseteq X$ that contains $x$ and is disjoint from an open set
that contains $\partial W$, which means that
\begin{equation}
\label{overline{U} cap partial W = emptyset}
        \overline{U} \cap \partial W = \emptyset.
\end{equation}
If $U_1 = U \cap W$, then $U_1$ is an open set in $X$ that contains
$x$ and satisfies $\overline{U_1} \subseteq \overline{U} \cap
\overline{W}$, which implies that
\begin{equation}
\label{overline{U_1} subseteq W}
        \overline{U_1} \subseteq W,
\end{equation}
by (\ref{overline{U} cap partial W = emptyset}).  Using this, one can
check that locally compact Hausdorff spaces are regular, since one can
always replace $W$ with a smaller open set if necessary to get that
$\overline{W}$ is compact.

\section{Dimension $0$}
\label{dimension 0}
\setcounter{equation}{0}

        A subset $E$ of a topological space $X$ is said to be
\emph{connected}\index{connected sets} in $X$ if $E$ cannot be
expressed as the union of two nonempty separated sets in $X$.  If $E
\subseteq Y \subseteq X$, then $E$ is connected in $X$ if and only if
$E$ is connected in $Y$, with respect to the induced topology on $Y$.
This follows from the analogous statement for separated sets,
which was mentioned in the previous section.  
As before, disjoint closed sets in $X$ are automatically separated, and
disjoint open subsets of $X$ are separated too.  If $A, B \subseteq X$
are separated and $A \cup B = X$, then $A$ and $B$ are both open and
both closed.  Thus $X$ is connected if and only if it cannot be
expressed as the union of two disjoint nonempty open sets, which is
equivalent to saying that $X$ cannot be expressed as the union of two
disjoint nonempty closed sets.  A set $E \subseteq X$ is said to be
\emph{totally disconnected}\index{totally disconnected sets} if it
does not contain any connected sets with at least two elements.

        A topological space $X$ is said to be \emph{totally
separated}\index{totally separated spaces} if for every $x, y \in X$
with $x \ne y$ there are disjoint open subsets $U$, $V$ of $X$ such
that $x \in U$, $y \in V$, and $U \cup V = X$.  Note that $U$ and $V$
are also closed sets in $X$ under these conditions, so that totally
separated spaces are completely Hausdorff.  If $X$ is totally
separated and $Y \subseteq X$, then $Y$ is totally separated with
respect to the induced topology.  If $\tau_1$ and $\tau_2$ are
topologies on a set $X$ such that $\tau_1 \subseteq \tau_2$, and if
$X$ is totally separated with respect to $\tau_1$, then $X$ is totally
separated with respect to $\tau_2$ as well.  Totally separated spaces
are totally disconnected, which can be derived from the previous
statement about subspaces of totally separated spaces and the fact
that totally separated spaces with at least two elements are not
connected.

        A topological space $X$ is said to have \emph{topological
dimension $0$}\index{topological dimension 0@topological dimension $0$}
at a point $x \in X$ if for every open set $W \subseteq X$ with $x \in
W$ there is an open set $U \subseteq X$ such that $x \in U$, $U
\subseteq W$, and $U$ is also a closed set in $X$.  Of course, this is
the same as saying that there is a local base for the topology of $X$
at $x$ consisting of subsets of $X$ that are both open and closed.
Similarly, $X$ is said to have topological dimension $0$ if $X$ has
topological dimension $0$ at every point $x \in X$, which is the same
as saying that there is a base for the topology of $X$ consisting of
sets that are both open and closed.  One may also ask that $X$ be
nonempty, and define the topological dimension of the empty set to be
$-1$.  Ultrametric spaces have topological dimension $0$, because open
and closed balls of positive radius are both open and closed with
respect to the corresponding topology, as in Section \ref{open, closed
  balls}.

        If $X$ has topological dimension $0$ and $Y \subseteq X$, then
one can check that $Y$ has topological dimension $0$ with respect to
the induced topology.  More precisely, if the topological dimension of
the empty set is defined to be $-1$, then one should ask that $Y \ne
\emptyset$ too.  If $X$ satisfies the $0$th separation condition and
has topological dimension $0$, then it is easy to see that $X$ is
regular as a topological space.  In this case, $X$ is totally
separated, and in particular $X$ is totally disconnected.

        Suppose that $X$ is totally separated, $K \subseteq X$
is compact, and $x \in X \setminus K$.  If $y \in K$, then $y \ne x$,
and so there is an open set $V(y) \subseteq X$ that is also closed
such that $y \in V(y)$ and $x \not\in V(y)$, because $X$ is totally
separated.  It follows that $K$ can be covered by finitely many of
these sets $V(y)$, because $K$ is compact, which leads to an open set
$V \subseteq X$ such that $K \subseteq V$, $x \not\in V$, and $V$ is a
closed set too.  Equivalently, $U = X \setminus V$ is an open set that
is closed as well and satisfies $x \in U$ and $K \cap U = \emptyset$.
If $H$ and $K$ are disjoint compact subsets of $X$, then one can repeat
the process to get a subset of $X$ that is both open and closed, and which
contains $H$ and is disjoint from $K$.

        If $X \ne \emptyset$ is totally separated and compact, then $X$
has topological dimension $0$.  To see this, let $x \in X$ and an open
set $W \subseteq X$ be given, with $x \in W$.  Thus $X \setminus W$ is
a closed set in $X$, so that $X \setminus W$ is compact, because $X$
is compact.  As in the previous paragraph, there is a set $U \subseteq
X$ that is both open and closed, which contains $x$, and is disjoint
from $X \setminus W$.  Hence $U \subseteq W$, as desired.

        As another variant of this type of argument, suppose that $X$
has topological dimension $0$, $K \subseteq X$ is compact, and that
$W \subseteq X$ is an open set that contains $K$.  Thus each element
of $K$ is contained in a subset of $W$ that is both open and closed
in $X$.  It follows that $K$ is contained in a subset of $W$ that
is both open and closed in $X$, using compactness of $K$ to reduce
to a finite subcovering.

        If $X$ is locally compact and has topological dimension $0$,
then for each $x \in X$ and open set $W \subseteq X$ with $x \in W$
there is an open set $U \subseteq X$ such that $x \in U$, $U \subseteq
W$, and $U$ is also closed and compact.  More precisely, if $X$ is
locally compact, then we can always replace $W$ by a smaller open set
that contains $x$ and is contained in a compact set.  This implies
that $U$ is compact in this situation, since it is a closed set
contained in a compact set.  Similarly, if $X$ is locally compact and
has topological dimension $0$, $H \subseteq X$ is compact, and $W
\subseteq X$ is an open set such that $H \subseteq W$, then $H$ is
contained in a subset of $W$ that is open, closed, and compact.

        Suppose that $X$ is totally separated, $W \subseteq X$ is
an open set, $x \in W$, and $\overline{W}$ is compact.  This implies
that $\partial W$ is a compact set that does not contain $x$, so that
there is an open set $U \subseteq X$ that is also closed, contains $x$,
and satisfies
\begin{equation}
\label{U cap partial W = emptyset}
        U \cap \partial W = \emptyset,
\end{equation}
as before.  Thus $U_1 = U \cap W$ is an open set in $X$ that contains
$x$ and satisfies $\overline{U_1} \subseteq U \cap \overline{W}$,
which implies that
\begin{equation}
\label{overline{U_1} subseteq U cap W = U_1}
        \overline{U_1} \subseteq U \cap W = U_1,
\end{equation}
by (\ref{U cap partial W = emptyset}).  Of course, this means that
$U_1$ is a closed set too.  It follows that a nonempty totally
separated locally compact topological space has topological dimension
$0$, since one can replace $W$ with a smaller open set to get
$\overline{W}$ to be compact, as usual.

\section{Chain connectedness}
\label{chain connectedness}
\setcounter{equation}{0}

        Suppose that $M$ is a nonempty set equipped with a $q$-metric
$d(x, y)$ for some $q > 0$, which leads to a topology on $M$, as in
Section \ref{open, closed balls}.  Of course, one can always reduce to
the case of ordinary metrics, using $d(x, y)^q$ when $q < 1$.  Let us
say that $A, B \subseteq M$ are
\emph{$\eta$-separated}\index{separated sets} in $M$ for some $\eta >
0$ if
\begin{equation}
\label{d(x, y) ge eta}
        d(x, y) \ge \eta
\end{equation}
for every $x \in A$ and $y \in B$.  This implies that $A$ and $B$ are
separated in the usual topological sense, and in fact that the
closures of $A$ and $B$ are disjoint.  In the other direction, if $A,
B \subseteq M$ are separated in the topological sense, and if at least
one of $A$ and $B$ is compact, then $A$ and $B$ are $\eta$-separated
for some $\eta > 0$.  More precisely, if $A$ is compact, then one may
as well suppose that $B$ is a closed set, since otherwise one can
replace $B$ with its closure.  The initial statement can be shown
using standard arguments without this observation, but it is perhaps
more commonly given in this way.

        A finite sequence $w_1, \ldots, w_n$ of elements of $M$ is said
to be an \emph{$\eta$-chain}\index{chains} for some $\eta > 0$ if
\begin{equation}
\label{d(w_j, w_{j + 1}) < eta}
        d(w_j, w_{j + 1}) < \eta
\end{equation}
for each $j$ with $1 \le j < n$, which is vacuous when $n = 1$.  Put
\begin{equation}
\label{x sim_eta y}
        x \sim_\eta y
\end{equation}
when $x, y \in M$ can be connected by an $\eta$-chain in $M$, which is
to say that there is an $\eta$-chain $w_1, \ldots, w_n$ of elements of
$M$ with $x = w_1$ and $y = w_n$.  It is easy to see that this defines
an equivalence relation on $M$, which leads to a partition of $M$ into
equivalence classes.  Each of these equivalence classes is an open set
in $M$, and in fact each equivalence class associated to (\ref{x
  sim_eta y}) contains the open ball of radius $\eta$ in $M$ centered
at any element of the equivalence class.  Any two distinct equivalence
classes associated to (\ref{x sim_eta y}) are $\eta$-separated in $M$.

        Let us say that $M$ is \emph{$\eta$-connected}\index{connected sets}
if every pair of elements of $M$ can be connected by an $\eta$-chain
of elements of $M$.  If $M$ is not $\eta$-connected, then $M$ can be
expressed as the union of two nonempty $\eta$-separated subsets of
$M$.  More precisely, if $M$ is not $\eta$-connected, then there are
points $x, y \in M$ that cannot be connected by an $\eta$-chain of
elements of $M$.  Let $A$ be the set of points in $M$ that can be
connected to $x$ by an $\eta$-chain of elements of $M$, and put $B = M
\setminus A$.  Thus $x \in A$, $y \in B$, $A \cup B = M$, and one can
check that $A$ and $B$ are $\eta$-separated in $M$.  Conversely, if
$A$, $B$ are $\eta$-separated subsets of $M$ such that $A \cup B = M$,
then there is no $\eta$-chain of elements of $M$ that connects a point
in $A$ to a point in $B$.  This is because such an $\eta$-chain would
have to go directly from an element of $A$ to an element of $B$ at some
step, which is not possible if $A$ and $B$ are $\eta$-separated in $M$.
It follows that $M$ is not $\eta$-connected when $M$ can be expressed
as the union of two nonempty $\eta$-separated sets.

        Similarly, a set $E \subseteq M$ is said to be
$\eta$-connected\index{connected sets} if every pair of
elements of $E$ can be connected by an $\eta$-chain of elements of
$E$.  Equivalently, $E$ is $\eta$-connected if $E$ cannot be
expressed as the union of two nonempty $\eta$-separated sets.
This follows from the discussion in the previous paragraph when
$E = M$.  Otherwise, one can reduce to that case, because $E$
is $\eta$-connected as a subset of $M$ if and only if $E$ is
$\eta$-connected as a subset of itself, using the restriction of
$d(x, y)$ to $x, y \in E$.

        If $E$ is $\eta$-connected for every $\eta > 0$, then we say that
$E \subseteq M$ is \emph{chain connected}.\index{chain connectedness}
Thus if $E$ is not chain connected, then $E$ is not $\eta$-connected
for some $\eta > 0$, so that $E$ can be expressed as the union of two
nonempty $\eta$-separated sets.  This implies that $E$ is not
connected, since $\eta$-separated sets are separated in the usual
sense.  It follows that connected subsets of $M$ are chain connected.
In the other direction, if $E \subseteq M$ is compact and not
connected, then $E$ can be expressed as the union of two nonempty
separated sets $A$ and $B$, and one can check that $A$ and $B$ also
have to be compact in this case.  This implies that $A$ and $B$ are
$\eta$-separated for some $\eta > 0$, as mentioned earlier, so that
$E$ is not $\eta$-connected.  Hence compact chain-connected subsets of
$M$ are connected.

        Let us say that $M$ is \emph{strongly totally
separated}\index{strongly totally separated spaces} if for each $x, y \in M$
with $x \ne y$ there are an $\eta > 0$ and $\eta$-separated sets $U, V
\subseteq M$ such that $x \in U$, $y \in V$, and $U \cup V = M$.  Note
that $U$ and $V$ have to be open subsets of $M$ under these
conditions, since they are separated and their union is equal to $M$.
Thus $M$ is totally separated when $M$ is strongly totally separated.
Equivalently, $M$ is strongly totally separated if for each $x, y \in
M$ with $x \ne y$ there is an $\eta > 0$ such that $x$ and $y$ cannot
be connected by an $\eta$-chain of elements of $M$, as in the earlier
discussion of $\eta$-connectedness.  If $M$ is strongly totally
separated and $Y \subseteq M$, then it is easy to see that $Y$ is
strongly totally separated too, with respect to the restriction of
$d(x, y)$ to $x, y \in Y$.

        Similarly, let us say that $M$ is \emph{strongly
  $0$-dimensional}\index{strongly 0-dimensional spaces@strongly
$0$-dimensional spaces} if for each $x \in M$ and $r > 0$ there is
an open set $U \subseteq M$ such that $x \in U$, $U \subseteq B(x,
r)$, and $U$, $M \setminus U$ are $\eta$-separated for some $\eta >
0$.  This implies that $M$ is strongly totally separated, and that $M$
has topological dimension $0$.  As before, one may wish to require
that $M$ be nonempty in order to be strongly $0$-dimensional, in
particular to be consistent in the second part of the preceding
statement.  If $M$ is strongly $0$-dimensional and $Y \subseteq M$,
then $Y$ is strongly $0$-dimensional with respect to the restriction
of $d(x, y)$ to $x, y \in Y$.  If nonemptiness is included in the
definition of strongly $0$-dimensional spaces, then one should also
ask that $Y$ be nonempty in the previous statement.

        If $M$ has topological dimension $0$ and is locally compact,
then $M$ is strongly $0$-dimensional.  This uses the fact that if $U
\subseteq M$ is compact and open, then $U$ and $M \setminus U$ are
$\eta$-separated for some $\eta > 0$.  Note that the set ${\bf Q}$ of
rational numbers has topological dimension $0$ with respect to the
standard topology, even though ${\bf Q}$ is chain connected with
respect to the standard metric on ${\bf R}$.  If $d(x, y)$ is an
ultrametric on a set $M$, then $M$ is strongly $0$-dimensional with
respect to the corresponding topology, because $B(x, r)$ and $M
\setminus B(x, r)$ are $r$-separated for every $x \in M$ and $r > 0$.

        Suppose that $M$ is strongly totally separated, $K \subseteq M$
is compact, and $x$ is an element of $M \setminus K$.  Using a
covering argument as in the previous section, one can check that there
is an open set $V \subseteq M$ such that $K \subseteq V$, $x \not\in
V$, and $V$, $M \setminus V$ are $\eta$-separated for some $\eta > 0$.
If $H$, $K$ are disjoint compact subsets of $M$, then one can repeat
the process to get an open set $U \subseteq M$ such that $H \subseteq
U$, $U \cap K = \emptyset$, and $U$, $M \setminus U$ are
$\eta$-separated for some $\eta > 0$.  If $M$ is strongly
$0$-dimensional, $K \subseteq M$ is compact, $W \subseteq M$ is an
open set, and $K \subseteq W$, then an analogous argument implies that
there is an open set $U \subseteq M$ such that $K \subseteq U$, $U
\subseteq W$, and $U$, $M \setminus U$ are $\eta$-separated for some
$\eta > 0$.  If $M$ is also locally compact, then one can take $U$ to
be compact as well.

\section{$\ell^r$ Spaces}
\label{l^r spaces}
\setcounter{equation}{0}

        Let $k$ be a field with an ultrametric absolute value function
$|\cdot|$, and let $V$ be a vector space over $k$ with an ultranorm
$N$ with respect to $|\cdot|$ on $k$.  Also let $X$ be a nonempty set,
so that $\ell^r(X, V)$ can be defined as in Sections \ref{supremum
  metrics, norms} and \ref{l^r norms} for $0 < r \le \infty$.  Under
these conditions, $\|f\|_r$ defines an $r$-norm on $\ell^r(X, V)$,
which leads to an $r$-metric on $\ell^r(X, V)$, as in Section
\ref{q-norms}.  In particular, the supremum norm defines an ultranorm
on $\ell^\infty(X, V)$ in this situation, and the corresponding
supremum metric is an ultrametric.  Thus $\ell^\infty(X, V)$ is
strongly $0$-dimensional with respect to the supremum metric, as in
the previous section.

        Suppose from now on in this section that $0 < r < \infty$.
Remember that $r$-summable functions are bounded on $X$, so that
\begin{equation}
\label{ell^r(X, V) subseteq ell^infty(X, V)}
        \ell^r(X, V) \subseteq \ell^\infty(X, V).
\end{equation}
This implies that $\ell^r(X, V)$ is strongly $0$-dimensional with
respect to the supremum metric, because of the analogous property of
$\ell^\infty(X, V)$.  Similarly,
\begin{equation}
\label{||f||_infty le ||f||_r, 2}
        \|f\|_\infty \le \|f\|_r
\end{equation}
for every $f \in \ell^r(X, V)$, which means that the $r$-metric on
$\ell^r(X, V)$ associated to $\|f\|_r$ is greater than or equal to the
supremum metric.  If $A, B \subseteq \ell^r(X, V)$ are
$\eta$-separated with respect to the supremum metric for some $\eta >
0$, then it follows that $A$ and $B$ are also $\eta$-separated with
respect to the $r$-metric associated to $\|f\|_r$.  Of course,
$\ell^\infty(X, V)$ is strongly totally separated with respect to the
supremum metric, since it is strongly $0$-dimensional.  It follows
that $\ell^r(X, V)$ is also strongly totally separated with respect to
the supremum metric, and hence with respect to the $r$-metric
associated to $\|f\|_r$.

        Suppose for the moment that $|\cdot|$ is the trivial absolute
value function on $k$, and that $N$ is the trivial ultranorm on $V$.
This implies that the supremum norm is the trivial ultranorm on
$\ell^\infty(X, V)$, so that the supremum metric on $\ell^\infty(X,
V)$ is the same as the discrete metric.  Note that every $r$-summable
$V$-valued function on $X$ has finite support in $X$ in this case, and
hence
\begin{equation}
\label{ell^r(X, V) = c_{00}(X, V)}
        \ell^r(X, V) = c_{00}(X, V).
\end{equation}
Using (\ref{||f||_infty le ||f||_r, 2}), we get that the $r$-metric
associated to $\|f\|_r$ on $\ell^r(X, V)$ is greater than or equal to
the discrete metric, which implies that the topology on $\ell^r(X, V)$
determined by the $r$-metric associated to $\|f\|_r$ is the discrete
topology.  It is easy to see that $\ell^r(X, V)$ is strongly
$0$-dimensional with respect to the $r$-metric associated to $\|f\|_r$
in this situation.

        Let $|\cdot|$ be any ultrametric absolute value function
on a field $k$ again, and let $N$ be any ultranorm on a vector
space $V$ over $k$.  If $X$ is a finite set with $n$ elements,
then every $V$-valued function $f$ on $X$ is $r$-summable, and
satisfies
\begin{equation}
\label{||f||_r le n^{1/r} ||f||_infty}
        \|f\|_r \le n^{1/r} \, \|f\|_\infty.
\end{equation}
Of course, this leads to a similar relationship between the $r$-metric
on $\ell^r(X, V) = \ell^\infty(X, V)$ associated to $\|f\|_r$ and the
supremum metric.  It follows that $\ell^r(X, V)$ is strongly
$0$-dimensional with respect to the $r$-metric associated to $\|f\|_r$
in this situation, because of the analogous property of
$\ell^\infty(X, V)$ with respect to the supremum metric.

        Let us suppose from now on in this section that $|\cdot|$
is nontrivial on $k$, $V \ne \{0\}$, and that $X$ has infinitely
many elements.  As in Section \ref{absolute value functions},
the nontriviality of $|\cdot|$ on $k$ means that there are nonzero
elements of $k$ with absolute value strictly less than $1$.  This
implies that there are nonzero elements of $k$ with arbitrarily
small absolute value, by taking large integer powers of the previous
elements.  It follows that there are nonzero elements of $V$ with
arbitrarily small norm, because $V \ne \{0\}$.

        Let $\eta > 0$ be given, and let $v_\eta$ be a nonzero element of $V$
with
\begin{equation}
\label{N(v_eta) < eta}
        N(v_\eta) < \eta,
\end{equation}
as in the preceding paragraph.  Also let $x_1, \ldots, x_n$ be
finitely many distinct elements of $X$.  If $j$ is a positive integer
less than or equal to $n$, then let $a_j(x)$ be the $V$-valued
function on $X$ defined by putting
\begin{equation}
\label{a_j(x_j) = v_eta}
        a_j(x_j) = v_\eta
\end{equation}
and $a_j(x) = 0$ when $x \ne x_j$.  Put
\begin{equation}
\label{f_l(x) = sum_{j = 1}^l a_j(x)}
        f_l(x) = \sum_{j = 1}^l a_j(x)
\end{equation}
for each $l = 1, \ldots, n$ and $x \in X$, and $f_0(x) = 0$ for every
$x \in X$.  Thus $f_l \in c_{00}(X, V) \subseteq \ell^r(X, V)$
for each $l = 0, 1, \ldots, n$, and
\begin{equation}
\label{||f_l - f_{l - 1}||_r = ||a_l||_r = N(v_eta)}
        \|f_l - f_{l - 1}\|_r = \|a_l\|_r = N(v_\eta)
\end{equation}
when $l \ge 1$.  This shows that $f_0, f_1, \ldots, f_n$ is an
$\eta$-chain in $c_{00}(X, V)$ with respect to the $r$-metric
associated to $\|f\|_r$.  We also have that
\begin{equation}
\label{||f_l||_r = l^{1/r} N(v_eta)}
        \|f_l\|_r = l^{1/r} \, N(v_\eta)
\end{equation}
for each $l$, because the $x_j$'s are supposed to be distinct elements
of $X$.

        Suppose that $U$ is a nonempty subset of $\ell^r(X, V)$ which
is $\eta$-separated from its complement in $\ell^r(X, V)$ with respect to
the $r$-metric associated to $\|f\|_r$.  This means that if an $\eta$-chain
of elements of $\ell^r(X, V)$ with respect to this $r$-metric starts at
an element of $U$, then this $\eta$-chain should stay in $U$ at every
step.  Using the $\eta$-chains described in the previous paragraph,
one can check that this implies that $U$ is unbounded with respect to
the $r$-metric associated to $\|f\|_r$.  It is convenient to reduce to
the case where $0 \in U$, although this is not really necessary.  This
also uses the hypothesis that $X$ have infinitely many elements, so that
the $\eta$-chain can have arbitrary length.  This implies that $\ell^r(X, V)$
is not strongly $0$-dimensional with respect to the $r$-metric associated
to $\|f\|_r$ under these conditions.  Similarly, $c_{00}(X, V)$ is not
strongly $0$-dimensional with respect to this $r$-metric.

        It is easy to see that $\ell^2({\bf Z}_+, {\bf Q})$ is totally
separated with respect to the $\ell^2$ metric, using the restriction
of the standard absolute value function on ${\bf R}$ to ${\bf Q}$.  A
well-known theorem of Erd\"os implies that $\ell^2({\bf Z}_+, {\bf
  Q})$ does not have topological dimension $0$, as in Example II 11 on
p13 of \cite{h-w}.  This argument seems to carry over nicely to
$\ell^r(X, V)$, under the same conditions as before, with some
adjustments.  In both situations, it suffices to show that a bounded
open set $U$ that contains $0$ has nonempty boundary.  To do this, one
looks for a sequence of elements of $U$ for which the distance to the
complement converges to $0$, and where the sequence converges in the
space being considered.  In the classical case of $\ell^2({\bf Z}_+,
{\bf Q})$, the $n$th term of the sequence has at most $n$ nonzero
coordinates, and one modifies the next coordinate to get closer to the
complement of $U$.  Similarly, in the context of $\ell^r(X, V)$, each
term in the sequence has only finitely many nonzero coordinates, and
each successive term modifies only finitely many coordinates that have
not been changed previously.  In these finitely many new coordinates,
one can use $\eta$-chains of the same type as before.  More precisely,
one adds an $\eta$-chain without leaving $U$, but where a single
additional step in the $\eta$-chain would leave $U$.  This is
possible, because $U$ is bounded, by hypothesis, and this ensures that
the resulting element of $U$ is as close to the complement of $U$ as
one wants, by taking $\eta$ to be sufficiently small.  As in the
classical case, it is easy to see that a sequence of elements of $U$
constructed in this way converges to an element of $\ell^r(X, V)$.
This uses the hypothesis that $U$ be bounded, so that the terms in the
sequence have bounded norm, and the fact that each new term in the
sequence only changes the coordinates that were equal to $0$ before,
by construction.  Of course, the limit of the sequence is an element
of the boundary of $U$, so that the boundary of $U$ is nonempty, as
desired.

\section{Some variants}
\label{some variants}
\setcounter{equation}{0}

        Let $k$ be a field with an absolute value function $|\cdot|$,
and let $V$ be a vector space over $k$ with a norm $N$ with respect to
$|\cdot|$ on $k$.  Also let $X$ be a nonempty set, let $r$ be a
positive real number, and let $a$ be a nonnegative real-valued
function on $X$ which is $r$-summable.  If $f$ is a $V$-valued
function on $X$ that satisfies
\begin{equation}
\label{N(f(x)) le a(x) for every x in X}
        N(f(x)) \le a(x) \quad\hbox{for every } x \in X, 
\end{equation}
then $f$ is $r$-summable on $X$ too.  Let $E_a$ be the set of $f : X
\to V$ that satisfy (\ref{N(f(x)) le a(x) for every x in X}), so that
$E_a \subseteq \ell^r(X, V)$.  This is the same as the classical
Hilbert cube when $V = k = {\bf R}$ with the standard absolute value
function, $X = {\bf Z}_+$, $r = 2$, and $a(j) = 1/j$ for each $j \in
{\bf Z}_+$.  Suppose for the moment that $V = k = {\bf Q}$, equipped
with the restriction of the standard absolute value function on ${\bf
  R}$ to ${\bf Q}$.  Remember that $\|f\|_r$ defines a norm on
$\ell^r(X, {\bf Q})$ when $r \ge 1$, and an $r$-norm when $0 < r \le
1$, as in Section \ref{l^r norms}.  One can check that $E_a$ has
topological dimension $0$ with respect to the topology determined by
the metric or $r$-metric corresponding to $\|f\|_r$.  This corresponds
to Example II 9 on p12 of \cite{h-w} when $X = {\bf Z}_+$, $r = 2$,
and $a(j) = 1/j$ for each $j \in {\bf Z}_+$, and essentially the same
argument can be used otherwise.

        Suppose now that $|\cdot|$ is an ultrametric absolute value
function on a field $k$, and that $N$ is an ultranorm on a vector
space $V$ over $k$.  Thus $\|f\|_r$ defines an $r$-norm on $\ell^r(X, V)$,
as in Section \ref{l^r norms}.  In this case, one can check that $E_a$
is strongly $0$-dimensional with respect to the $r$-metric associated
to $\|f\|_r$.  This uses the fact that for each $\epsilon > 0$ there
is a finite set $A(\epsilon) \subseteq X$ such that
\begin{equation}
\label{sum_{x in X setminus A(epsilon)} a(x)^r < epsilon}
        \sum_{x \in X \setminus A(\epsilon)} a(x)^r < \epsilon,
\end{equation}
as in (\ref{sum_{x in X setminus A(epsilon)} N(f(x))^r < epsilon}).
Of course, every subset of $\ell^r(X, V)$ is strongly totally
separated with respect to the $r$-metric associated to $\|f\|_r$ under
these conditions, since $\ell^r(X, V)$ is strongly totally separated
with respect to this $r$-metric, as in the previous section.

        Let us continue to suppose that $|\cdot|$ be an ultrametric
absolute value function on $k$, and that $N$ be an ultranorm on $V$.
Also let $X$ be a nonempty set, and let $r$ be a positive real number,
as before.  Consider the vector space $c_{00}(X, V)$ of $V$-valued
functions on $X$ with finite support, equipped with the $r$-metric
associated to $\|f\|_r$.  Suppose for the moment that $N$ takes values
on $V$ in a set of finitely or countably many nonnegative real
numbers.  This implies that $N(v)^r$ also takes values in a set of
finitely or countably many nonnegative real numbers, and hence that
the collection of all finite sums of elements of this set has only
finitely or countably many elements.  It follows that $\|f\|_r$ takes
only finitely or countably many values on $c_{00}(X, V)$, which
implies that the corresponding $r$-metric only takes finitely or
countably many values on $c_{00}(X, V)$ as well.  In this case,
$c_{00}(X, V)$ has topological dimension $0$ with respect to the
topology determined by the $r$-metric associated to $\|f\|_r$,
since open and closed balls of all but finitely or countably many radii
are automatically the same.

        Otherwise, if $N$ does not take values in a set of nonnegative
real numbers with only finitely or countably many elements, then we
can basically reduce to this case by modifying $N$.  More precisely,
let $h(t)$ be a monotonically increasing real-valued function defined
on the set of nonnegative real numbers such that $h(0) = 0$ and $h(t)
> 0$ when $t > 0$.  Under these conditions, it is easy to see that
\begin{equation}
\label{h(N(v - w))}
        h(N(v - w))
\end{equation}
defines an ultrametric on $V$ which determines the same topology on
$V$ as the ultrametric $N(v - w)$ associated to $N$.  Put
\begin{equation}
\label{d_r(f, g) = (sum_{x in X} h(N(f(x) - g(x)))^r)^{1/r}}
        d_r(f, g) = \Big(\sum_{x \in X} h(N(f(x) - g(x)))^r\Big)^{1/r}
\end{equation}
for every $f, g \in c_{00}(X, V)$, which defines an $r$-metric on
$c_{00}(X, V)$, for the same reasons as in Section \ref{l^r norms}.
If $h(t)$ and $t$ are each bounded by positive constant multiples of
the other on $[0, +\infty)$, then (\ref{d_r(f, g) = (sum_{x in X} h(N(f(x) -
  g(x)))^r)^{1/r}}) and $\|f - g\|_r$ are each bounded by the same
constant multiples of the other on $c_{00}(X, V)$.  In particular,
this implies that these two $r$-metrics determine the same topology on
$c_{00}(X, V)$.  We can also choose $h$ so that it takes values in a
countable subset of ${\bf R}$, which implies that (\ref{d_r(f, g) =
  (sum_{x in X} h(N(f(x) - g(x)))^r)^{1/r}}) takes values in a
countable set of nonnegative real numbers when $f, g \in c_{00}(X,
V)$.  As before, this means that open and closed balls in $c_{00}(X,
V)$ with respect to (\ref{d_r(f, g) = (sum_{x in X} h(N(f(x) -
  g(x)))^r)^{1/r}}) of all but finitely or countably many radii are
the same, and hence that $c_{00}(X, V)$ has topological dimension $0$
with respect to the corresponding topology.

\section{$\ell^r$ Spaces, continued}
\label{l^r spaces, continued}
\setcounter{equation}{0}

        Let $k$ be a field with an absolute value function $|\cdot|$
again, and let $V$ be a vector space over $k$ with a norm $N$
with respect to $|\cdot|$ on $k$.  Also let $X$ be a nonempty set,
let $r$ and $t$ be positive real numbers, and let $f$ be an
$r$-summable $V$-valued function on $X$ with
\begin{equation}
\label{||f||_r = (sum_{x in X} N(f(x))^r)^{1/r} = t}
        \|f\|_r = \Big(\sum_{x \in X} N(f(x))^r\Big)^{1/r} = t.
\end{equation}
Thus for each $\epsilon > 0$ there is a finite set $A(\epsilon)
\subseteq X$ such that
\begin{equation}
\label{sum_{x in A(epsilon)} N(f(x))^r > ... = t^r - epsilon}
        \sum_{x \in A(\epsilon)} N(f(x))^r > \sum_{x \in X} N(f(x))^r - \epsilon
                                         = t^r - \epsilon ,
\end{equation}
as in (\ref{sum_{x in X} N(f(x))^r < sum_{x in A(epsilon)} N(f(x))^r +
  epsilon}).  As before, this implies that
\begin{equation}
\label{sum_{x in X setminus A(epsilon)} N(f(x))^r < epsilon, 2}
        \sum_{x \in X \setminus A(\epsilon)} N(f(x))^r < \epsilon.
\end{equation}
Suppose that $g$ is another $V$-valued function on $X$ that is
sufficiently close to $f$ on $A(\epsilon)$ so that
\begin{equation}
\label{sum_{x in A(epsilon)} N(g(x))^r > t^r - 2 epsilon}
        \sum_{x \in A(\epsilon)} N(g(x))^r > t^r - 2 \, \epsilon.
\end{equation}
If we also have that $g \in \ell^r(X, V)$ satisfies
\begin{equation}
\label{||g||_r = (sum_{x in X} N(g(x))^r)^{1/r} le t}
        \|g\|_r = \Big(\sum_{x \in X} N(g(x))^r\Big)^{1/r} \le t,
\end{equation}
then it follows that
\begin{eqnarray}
\label{sum_{x in X setminus A(epsilon)} N(g(x))^r = ... = 2 epsilon}
        \sum_{x \in X \setminus A(\epsilon)} N(g(x))^r
         & = & \sum_{x \in X} N(g(x))^r - \sum_{x \in A(\epsilon)} N(g(x))^r \\
         & < & t^r - (t^r - 2 \, \epsilon) = 2 \, \epsilon. \nonumber
\end{eqnarray}
This permits us to estimate
\begin{equation}
\label{||f - g||_r^r = ...}
        \|f - g\|_r^r = \sum_{x \in A(\epsilon)} N(f(x) - g(x))^r
                         + \sum_{x \in X \setminus A(\epsilon)} N(f(x) - g(x))^r
\end{equation}
in terms of how close $g$ is to $f$ on $A(\epsilon)$ under these
conditions, using (\ref{sum_{x in X setminus A(epsilon)} N(f(x))^r <
  epsilon, 2}) and (\ref{sum_{x in X setminus A(epsilon)} N(g(x))^r =
  ... = 2 epsilon}).

        Suppose for the moment that $V = k = {\bf Q}$, equipped with
the restriction of the standard absolute value function on ${\bf R}$
to ${\bf Q}$.  Using the remarks in the previous paragraph, one can
show that the sphere of radius $t$ in $\ell^r(X, {\bf Q})$ centered at
$0$ has topological dimension $0$.  Of course, the same argument shows
every sphere in $\ell^r(X, {\bf Q})$ has topological dimension $0$,
which means that $\ell^r(X, {\bf Q})$ has topological dimension $\le
1$.  This is basically the same as Example III 5 on p25f of
\cite{h-w}, in which one takes $X = {\bf Z}_+$, $r = 2$, and $t < 1$.
The argument in \cite{h-w} uses an embedding of the sphere into the
Hilbert cube, but this is just a convenience.  Basically the same type
of argument can be used for the sphere as for the Hilbert cube,
because of the remarks in the previous paragraph.  This is the other
part of Erd\"os' famous theorem that $\ell^2({\bf Z}_+, {\bf Q})$ has
topological dimension equal to $1$.

        Now let $|\cdot|$ be an ultrametric absolute value function
on any field $k$, and let $N$ be an ultranorm on $V$ with respect to
$|\cdot|$ on $k$.  In this case, one can use the earlier remarks to
show that spheres in $\ell^r(X, V)$ with respect to $\|f\|_r$ are
strongly $0$-dimensional.  In particular, this implies that $\ell^r(X,
V)$ has topological dimension $\le 1$.  If $X$ has infinitely many
elements, $|\cdot|$ is nontrivial on $k$, and $V \ne \{0\}$, then we
have already seen that $\ell^r(X, V)$ does not have topological
dimension $0$, as in Section \ref{l^r spaces}.  It follows that
$\ell^r(X, V)$ also has topological dimension $1$ under these
conditions.

\section{Uniform conditions}
\label{uniform conditions}
\setcounter{equation}{0}

        Let $M$ be a nonempty set with a $q$-metric $d(x, y)$ for some
$q > 0$.  Let us say that $M$ is
\emph{uniformly totally separated}\index{uniformly totally separated spaces}
if for each $r > 0$ there is an $\eta(r) > 0$ such that for every
$x, y \in M$ with $d(x, y) \ge r$ there are $\eta(r)$-separated sets
$U, V \subseteq M$ with $x \in U$, $y \in V$, and $U \cup V = M$.
Note that this implies that $M$ is strongly totally separated, as
in Section \ref{chain connectedness}.  If $M$ is uniformly totally
separated, then it follows that
\begin{eqnarray}
\label{uniformly totally separated condition}
 & & \hbox{for each $r > 0$ there is an $\eta(r) > 0$ such that for every }
               x, y \in M \\
 & & \hbox{with } d(x, y) \ge r,
      \hbox{ we have that $x$ and $y$ cannot be connected by} \nonumber \\
 & & \hbox{an $\eta(r)$-chain in $M$.} \nonumber
\end{eqnarray}
More precisely, if $M$ is uniformly totally separated, then
(\ref{uniformly totally separated condition}) holds with the same
choice of $\eta(r)$ as in the initial definition.

        Conversely, suppose that $M$ satisfies (\ref{uniformly totally 
separated condition}), and let $r > 0$ and $x \in M$ be given.  Put
\begin{equation}
\label{U = {z in M : x can be connected to z by an eta(r)-chain in M}}
  U = \{z \in M : \hbox{ $x$ can be connected to $z$ by an $\eta(r)$-chain in }
                                                               M\},
\end{equation}
where $\eta(r)$ is as in (\ref{uniformly totally separated
  condition}).  Thus $x \in U$ automatically, and it is easy to see
that $U$, $M \setminus U$ are $\eta(r)$-separated in $M$.  By
hypothesis, $M \setminus U$ contains every $y \in M$ with $d(x, y) \ge
r$, which is the same as saying that
\begin{equation}
\label{U subseteq B(x, r)}
        U \subseteq B(x, r).
\end{equation}
In particular, this implies that $M$ is uniformly totally separated,
with $V = M \setminus U$, and with the same choice of $\eta(r)$.

        Let us say that $M$ is \emph{uniformly $0$-dimensional}\index{uniformly
0-dimensional spaces@uniformly $0$-dimensional spaces} if for each $r > 0$
there is an $\eta(r) > 0$ such that for every $x \in M$ there is an
open set $U \subseteq M$ with $x \in U$ that satisfies (\ref{U
  subseteq B(x, r)}) and has the property that $U$, $M \setminus U$
are $\eta(r)$-separated in $M$.  This condition clearly implies that
$M$ is strongly $0$-dimensional, and that $M$ is uniformly totally
separated, with the same choice of $\eta(r)$.  In fact, the argument
in the previous paragraph shows that uniformly totally separated
spaces are uniformly $0$-dimensional, with the same choice of
$\eta(r)$.  This is because $U \subseteq M$ is automatically an open
set when $U$, $M \setminus U$ are separated in $M$.

        As a variant of this, let us say that $M$ is uniformly totally
separated at $x \in M$ if for each $r > 0$ there is an $\eta(x, r) >
0$ such that for every $y \in M$ with $d(x, y) \ge r$ there are
$\eta(x, r)$-separated sets $U, V \subseteq M$ with $x \in U$, $y \in
V$, and $U \cap V = M$.  This implies that
\begin{eqnarray}
\label{uniformly separated at x}
 & & \hbox{for each $r > 0$ there is an $\eta(x, r) > 0$ such that for every }
                                          y \in M \\
 & & \hbox{with } d(x, y) \ge r,
      \hbox{ we have that $x$ and $y$ cannot be connected by} \nonumber \\
 & & \hbox{an $\eta(x, r)$-chain in $M$,} \nonumber
\end{eqnarray}
with the same choice of $\eta(x, r)$ as in the previous definition.
Conversely, suppose that $M$ satisfies (\ref{uniformly separated at
  x}), and let $r > 0$ be given.  Put
\begin{equation}
\label{U = {z in M : x can be connected to z by an eta(x, r)-chain in M}}
  \quad U = \{z \in M : \hbox{ $x$ can be connected to } z
                               \hbox{ by an $\eta(x, r)$-chain in } M\},
\end{equation}
where $\eta(x, r)$ is as in (\ref{uniformly separated at x}).  This is
the same as (\ref{U = {z in M : x can be connected to z by an
    eta(r)-chain in M}}), but with $\eta(r)$ replaced by $\eta(x, r)$.
As before, $x \in U$ automatically, and $U$, $M \setminus U$ are
$\eta(x, r)$-separated in $M$.  Our hypothesis (\ref{uniformly
  separated at x}) says exactly that $U$ also satisfies (\ref{U
  subseteq B(x, r)}).  This implies that $M$ is uniformly totally
separated at $x$, with $V = M \setminus U$, and with the same choice
of $\eta(x, r)$.

        Let us say that $M$ is strongly $0$-dimensional at $x \in M$
if\index{strongly 0-dimensional spaces@strongly $0$-dimensional spaces}
for each $r > 0$ there is an $\eta = \eta(x, r) > 0$ and an open set
$U \subseteq M$ such that $x \in U$, $U$ satisfies (\ref{U subseteq
  B(x, r)}), and $U$, $M \setminus U$ are $\eta$-separated in $M$.
Thus $M$ is strongly $0$-dimensional as defined in Section \ref{chain
  connectedness} if and only if $M$ is strongly $0$-dimensional at
each point $x \in M$.  If $M$ is strongly $0$-dimensional at $x$, then
it is easy to see that $M$ is uniformly totally separated at $x$, with
the same choice of $\eta(x, r)$.  Conversely, if $M$ is strongly
totally separated at $x$, then $M$ is strongly $0$-dimensional at $x$,
with the same choice of $\eta(x, r)$, by the argument in the previous
paragraph.  As usual, this uses the fact that $U \subseteq M$ is an
open set when $U$, $M \setminus U$ are separated in $M$.

        Suppose that $K \subseteq M$ is compact, and that $M$ is
strongly $0$-dimensional at each $x \in K$.  Let $r > 0$ be given,
so that for each $x \in K$ there is an $\eta(x, r) > 0$ and an
open set $U(x, r)$ such that $x \in U(x, r)$,
\begin{equation}
\label{U(x, r) subseteq B(x, r)}
        U(x, r) \subseteq B(x, r),
\end{equation}
and
\begin{equation}
\label{U(x, r), M setminus U(x, r) are eta(x, r)-separated in M}
 U(x, r), \ M \setminus U(x, r) \hbox{ are $\eta(x, r)$-separated in } M.
\end{equation}
Because $K$ is compact, there are finitely many points $x_1, \ldots, x_n
\in K$ such that
\begin{equation}
\label{K subseteq bigcup_{j = 1}^n U(x_j, r)}
        K \subseteq \bigcup_{j = 1}^n U(x_j, r).
\end{equation}
Put
\begin{equation}
\label{eta = min_{1 le j le n} eta(x_j, r) > 0}
        \eta = \min_{1 \le j \le n} \eta(x_j, r) > 0,
\end{equation}
so that
\begin{equation}
\label{U(x_j, r), M setminus U(x_j, r) are eta-separated in M}
        U(x_j, r), \ M \setminus U(x_j, r) \hbox{ are $\eta$-separated in } M
\end{equation}
for each $j = 1, \ldots, n$, by (\ref{U(x, r), M setminus U(x, r) are
  eta(x, r)-separated in M}).  Let $w \in K$ be given, and let $j$ be
an integer such that $1 \le j \le n$ and $w \in U(x_j, r)$, as in
(\ref{K subseteq bigcup_{j = 1}^n U(x_j, r)}).  This implies that
\begin{equation}
\label{d(x_j, w) < r}
        d(x_j, w) < r,
\end{equation}
by (\ref{U(x, r) subseteq B(x, r)}) with $x = x_j$, and hence that
\begin{equation}
\label{B(x_j, r) subseteq B(w, 2^{1/q} r)}
        B(x_j, r) \subseteq B(w, 2^{1/q} \, r),
\end{equation}
since $d(\cdot, \cdot)$ is a $q$-metric on $M$.  It follows that
\begin{equation}
\label{U(x_j, r) subseteq B(w, 2^{1/q} r)}
        U(x_j, r) \subseteq B(w, 2^{1/q} \, r),
\end{equation}
by combining (\ref{U(x, r) subseteq B(x, r)}) with $x = x_j$ and
(\ref{B(x_j, r) subseteq B(w, 2^{1/q} r)}).  This shows that $M$
satisfies a version of being strongly $0$-dimensional at each $w \in
K$, with a choice of $\eta > 0$ that depends on the radius and not
$w$.  In particular, if $M$ is compact and strongly $0$-dimensional,
then $M$ is uniformly $0$-dimensional.

        If $d(\cdot, \cdot)$ is an ultrametric on $M$, then $M$ is
uniformly $0$-dimensional, with $\eta(r) = r$, for the same reasons as
for strong $0$-dimensionality in Section \ref{chain connectedness}.
If $d(\cdot, \cdot)$ is any $q$-metric on $M$, $M$ is uniformly
$0$-dimensional, and $Y \subseteq M$, then it is easy to see that $Y$
is also uniformly $0$-dimensional with respect to the restriction of
$d(\cdot, \cdot)$ to $Y$, and with the same choice of $\eta(r)$.

\section{Some examples and remarks}
\label{some examples, remarks}
\setcounter{equation}{0}

        Of course, a subset $E$ of the real line is totally disconnected
with respect to the standard topology on ${\bf R}$ if and only if the
interior of $E$ is empty, which is the same as saying that ${\bf R}
\setminus E$ is dense in ${\bf R}$.  In this case, $E$ has topological
dimension $0$, at least when $E \ne \emptyset$, if that is included in
the definition.

        Let $E$ be a subset of ${\bf R}$ again, and let $x$, $y$ be
distinct elements of $E$.  We may as well suppose that $x < y$, since
this can always be arranged by interchanging the roles of $x$ and $y$,
when needed.  If $E \cap (x, y)$ is dense in $(x, y)$, then for each
$\eta > 0$, $x$ and $y$ can be connected by an $\eta$-chain of
elements of $E$ with respect to the standard metric on ${\bf R}$.
Thus if there is an $\eta > 0$ such that $x$ and $y$ cannot be
connected by an $\eta$-chain of elements of $E$, then $E \cap (x, y)$
is not dense in $(x, y)$.  This implies that
\begin{equation}
\label{(x, y) setminus overline{E} ne emptyset}
        (x, y) \setminus \overline{E} \ne \emptyset.
\end{equation}

        Suppose now that $E \subseteq {\bf R}$ is strongly totally separated
with respect to the restriction of the standard metric on ${\bf R}$ to
$E$.  Under these conditions, the argument in the preceding paragraph
implies that (\ref{(x, y) setminus overline{E} ne emptyset}) holds for
every $x, y \in E$ with $x < y$.  The same conclusion holds when $x$
and $y$ are elements of the closure of $E \cap [x, y]$, by
approximating $x$ and $y$ by elements of $E \cap [x, y]$, and applying
the previous argument to those approximations.  If $x$ or $y$ is not
in the closure of $E \cap [x, y]$, then it is easy to see that
(\ref{(x, y) setminus overline{E} ne emptyset}) still holds.  It
follows that (\ref{(x, y) setminus overline{E} ne emptyset}) holds for
every $x, y \in {\bf R}$ with $x < y$ when $E$ is strongly totally
separated in ${\bf R}$, which means that ${\bf R} \setminus
\overline{E}$ is dense in ${\bf R}$.

        Conversely, if ${\bf R} \setminus \overline{E}$ is dense in
${\bf R}$, then it is easy to see that $E$ is strongly $0$-dimensional,
which implies that $E$ is strongly totally separated.  As usual, one
should also ask that $E \ne \emptyset$ in the first part of the
preceding statement, if that is included in the definition of being
strongly $0$-dimensional.  More precisely, $\overline{E}$ is strongly
$0$-dimensional when ${\bf R} \setminus \overline{E}$ is dense in
${\bf R}$.  If $E$ is also bounded, then $\overline{E}$ is compact,
and hence $\overline{E}$ is uniformly $0$-dimensional, as in the
preceding section.  Otherwise, the same argument implies that bounded
subsets of $\overline{E}$ are uniformly $0$-dimensional.

        Let us consider some rather different examples in the plane,
with respect to the standard Euclidean metric.  Suppose that $E_j$ is
a finite subset of $[0, 1] \times \{1/j\}$ for each positive integer
$j$, and put
\begin{equation}
\label{E = bigcup_{j = 1}^infty E_j}
        E = \bigcup_{j = 1}^\infty E_j.
\end{equation}
It is easy to see that $E$ is strongly $0$-dimensional in this
situation, with respect to the restriction of the standard Euclidean
metric on ${\bf R}^2$ to $E$.  However, we can choose the $E_j$'s so
that
\begin{equation}
\label{[0, 1] times {0} subseteq overline{E}}
        [0, 1] \times \{0\} \subseteq \overline{E},
\end{equation}
where $\overline{E}$ is the closure of $E$ in ${\bf R}^2$.  In
particular, this implies that $\overline{E}$ is not totally
disconnected.

        Let $M$ be a nonempty set with a $q$-metric $d(x, y)$ for some
$q > 0$.  If $A, B \subseteq M$ are $\eta$-separated for some $\eta > 0$,
then it is easy to see that their closures $\overline{A}$, $\overline{B}$
are $\eta$-separated in $M$ as well.  Put $E = A \cup B$, so that
\begin{equation}
\label{overline{E} = overline{A} cup overline{B}}
        \overline{E} = \overline{A} \cup \overline{B}.
\end{equation}
If $A, B \ne \emptyset$, then $\overline{A}, \overline{B} \ne
\emptyset$, and hence $\overline{E}$ is not connected in $M$.

        Suppose now that $E$ is any subset of $M$ which is strongly
totally separated with respect to the restriction of $d(x, y)$ to $x,
y \in E$.  This means that for every $x, y \in E$ with $x \ne y$ there
are an $\eta > 0$ and $\eta$-separated sets $A, B \subseteq E$ such
that $x \in A$, $y \in B$, and $A \cup B = E$.  As in the previous
paragraph, $\overline{A}$ and $\overline{B}$ are also $\eta$-separated
in $M$ and satisfy (\ref{overline{E} = overline{A} cup overline{B}}).
Thus $\overline{E}$ has a property analogous to being strongly totally
separated, but which only applies to distinct elements of $E$, instead
of $\overline{E}$.

        Similarly, let us suppose that $E \subseteq M$ is strongly
$0$-dimensional with respect to the restriction of $d(x, y)$ to
$x, y \in E$.  This implies that for each $x \in E$ and $r > 0$
there are an $\eta > 0$ and $\eta$-separated sets $A, B \subseteq E$
such that $x \in A$, $A \cup B = E$, and $A$ is contained in the
open ball centered at $x$ with radius $r$.  It follows that $\overline{A}$
and $\overline{B}$ are $\eta$-separated subsets of $M$ that satisfy
(\ref{overline{E} = overline{A} cup overline{B}}), and that $\overline{A}$
is contained in the closed ball in $M$ centered at $x$ with radius $r$.
This shows that $\overline{E}$ is strongly $0$-dimensional at
every element of $E$.

        Let us return to the case where $M = {\bf R}^2$ with the
standard Euclidean metric, and let $E$ be as in (\ref{E = bigcup_{j =
    1}^infty E_j}).  If $z$ is any element of ${\bf R}^2$, then it is
easy to see that $E \cup \{z\}$ is strongly totally separated, with
respect to the restriction of the standard Euclidean metric on ${\bf
  R}^2$ to $E \cup \{z\}$.  More precisely, $E \cup \{z\}$ is strongly
$0$-dimensional at every element of $E$, for essentially the same
reasons as before.  If $z \not\in [0, 1] \times \{0\}$, then $E \cup
\{z\}$ is strongly $0$-dimensional at $z$ too, for essentially the
same reasons again.  Otherwise, if $z \in [0, 1] \times \{0\}$, and if
we choose the $E_j$'s so that (\ref{[0, 1] times {0} subseteq
  overline{E}}) holds, then $E \cup \{z\}$ is not strongly
$0$-dimensional at $z$.  In this case, for each $\eta > 0$, there is
an $\eta$-chain of elements of $E \cup \{z\}$ that starts at $z$ and
can go a distance which is at least almost $1/2$.  If $w$, $z$ are
distinct elements of $[0, 1] \times \{0\}$, and if (\ref{[0, 1] times
  {0} subseteq overline{E}}) holds, then $E \cup \{w, z\}$ is not
strongly totally separated.  This is because $w$ and $z$ can be
connected by an $\eta$-chain of elements of $E \cup \{w, z\}$ for
every $\eta > 0$.

\section{Some additional remarks}
\label{some additional remarks}
\setcounter{equation}{0}

        Let $M$ be a nonempty set with a $q$-metric $d(x, y)$ for
some $q > 0$, and let $E$ be a subset of $M$.  As before, $E$ is
strongly totally separated with respect to the restriction of
$d(\cdot, \cdot)$ to $E$ if for every $x, y \in E$ with $x \ne y$
there are an $\eta > 0$ and $\eta$-separated sets $A, B \subseteq E$
such that $x \in A$, $y \in B$, and $A \cup B = E$.  This implies that
$A$, $B$ are relatively open in $E$, and hence that there are $t_1,
t_2 > 0$ such that
\begin{equation}
\label{B(x, t_1) cap E subseteq A, B(y, t_2) cap E subseteq B}
        B(x, t_1) \cap E \subseteq A, \quad B(y, t_2) \cap E \subseteq B.
\end{equation}
Here $B(w, t)$ denotes the open ball in $M$ centered at $w \in M$ with
radius $t > 0$ with respect to $d(\cdot, \cdot)$, as usual.  Note that
(\ref{B(x, t_1) cap E subseteq A, B(y, t_2) cap E subseteq B}) holds
with $t_1 = t_2 = \eta$, but in some circumstances (\ref{B(x, t_1) cap
  E subseteq A, B(y, t_2) cap E subseteq B}) may hold with larger
values of $t_1$, $t_2$ as well.

        As in the previous section, the closures $\overline{A}$,
$\overline{B}$ of $A$, $B$ are $\eta$-separated in $M$ too, and satisfy
(\ref{overline{E} = overline{A} cup overline{B}}).  Observe that
\begin{equation}
\label{B(x, t_1) cap overline{E} subseteq overline{A}, ...}
        B(x, t_1) \cap \overline{E} \subseteq \overline{A}, \quad
          B(y, t_2) \cap \overline{E} \subseteq \overline{B},
\end{equation}
by (\ref{B(x, t_1) cap E subseteq A, B(y, t_2) cap E subseteq B}).
Let $x'$, $y'$ be distinct elements of $\overline{E}$, and suppose
that $x$, $y$ are distinct elements of $E$ that are very close to
$x'$, $y'$, respectively.  If
\begin{equation}
\label{d(x, x') < t_1, d(y, y') < t_2}
        d(x, x') < t_1, \quad d(y, y') < t_2,
\end{equation}
where $t_1$, $t_2$ are as in (\ref{B(x, t_1) cap E subseteq A, B(y,
  t_2) cap E subseteq B}), then (\ref{B(x, t_1) cap overline{E}
  subseteq overline{A}, ...})  implies that
\begin{equation}
\label{x' in overline{A}, y' in overline{B}}
        x' \in \overline{A}, \quad y' \in \overline{B}.
\end{equation}
One might like to try to use an argument like this to show that
$\overline{E}$ is strongly totally separated too, but the problem is
that $t_1$, $t_2$ may depend on $x$ and $y$, so that (\ref{d(x, x') <
  t_1, d(y, y') < t_2}) does not hold.

          If $E$ is uniformly totally separated, as in Section
\ref{uniform conditions}, then one can take $\eta$ to depend only on
a positive lower bound for $d(x, y)$.  In this case, the argument
indicated in the previous paragraph can be used, and in fact
$\overline{E}$ is also uniformly totally separated.  If $E \subseteq
{\bf R}$ and ${\bf R} \setminus \overline{E}$ is dense in ${\bf R}$,
then (\ref{B(x, t_1) cap E subseteq A, B(y, t_2) cap E subseteq B})
holds with respect to the standard metric on ${\bf R}$ for any $t_1,
t_2 > 0$ such that
\begin{equation}
\label{t_1 + t_2 < |x - y|}
        t_1 + t_2 < |x - y|.
\end{equation}
Of course, we already know that $\overline{E}$ is strongly
$0$-dimensional in this situation, and hence strongly totally
separated.

        Let $M$ be any nonempty set with a $q$-metric $d(\cdot, \cdot)$
again.  A subset $E$ of $M$ is strongly $0$-dimensional with respect
to the restriction of $d(\cdot, \cdot)$ to $E$ if for each $x \in E$
and $r > 0$ there are an $\eta > 0$ and $\eta$-separated sets $A, B
\subseteq E$ such that $x \in A$, $A \cup B = E$, and
\begin{equation}
\label{A subseteq B(x, r)}
        A \subseteq B(x, r).
\end{equation}
This implies that $A$, $B$ are relatively open in $E$, and in
particular that
\begin{equation}
\label{B(x, t) cap E subseteq A}
        B(x, t) \cap E \subseteq A
\end{equation}
for some $t > 0$.  More precisely, (\ref{B(x, t) cap E subseteq A})
holds with $t = \eta$, but it may also hold with larger values of $t$.
As before, $\overline{A}$, $\overline{B}$ are $\eta$-separated in $M$,
and satisfy (\ref{overline{E} = overline{A} cup overline{B}}).  We
also have that
\begin{equation}
\label{B(x, t) cap overline{E} subseteq overline{A} subseteq overline{B}(x, r)}
 B(x, t) \cap \overline{E} \subseteq \overline{A} \subseteq \overline{B}(x, r)
\end{equation}
because of (\ref{A subseteq B(x, r)}) and (\ref{B(x, t) cap E subseteq
  A}).  Using the second inclusion in (\ref{B(x, t) cap overline{E}
  subseteq overline{A} subseteq overline{B}(x, r)}), we get that
$\overline{E}$ is strongly $0$-dimensional at every $x \in E$, as in
the previous section.  One might like to show that $\overline{E}$ is
strongly $0$-dimensional at a point $x' \in \overline{E}$ by
approximating $x'$ by $x \in E$, in such a way that
\begin{equation}
\label{d(x, x') < t}
        d(x, x') < t,
\end{equation}
where $t$ is as in (\ref{B(x, t) cap E subseteq A}).  This does not
always work, because $t$ may depend on $x$.  This does work when $E$
is uniformly $0$-dimensional, which is equivalent to $E$ being
uniformly totally separated, as in Section \ref{uniform conditions}.
This also works when $E \subseteq {\bf R}$ and ${\bf R} \setminus
\overline{E}$ is dense in ${\bf R}$, in which case we already know
that $\overline{E}$ is strongly $0$-dimensional with respect to the
restriction of the standard metric on ${\bf R}$ to $\overline{E}$.

\section{Another perspective}
\label{another perspective}
\setcounter{equation}{0}

        Let $M$ be a nonempty set with a $q$-metric $d(x, y)$ for some
$q > 0$, and let $E$ be a subset of $M$.  Put
\begin{equation}
\label{Z(r) = {(x, y) in E times E : d(x, y) ge r}}
        Z(r) = \{(x, y) \in E \times E : d(x, y) \ge r\}
\end{equation}
for each $r > 0$, and let $Z$ be a subset of $Z(r)$ for some $r > 0$.
Let us say that $E$ is \emph{uniformly totally separated along
  $Z$}\index{uniformly totally separated spaces} if there is an $\eta
> 0$ such that for each $(x, y) \in Z$ there are $\eta$-separated sets
$A, B \subseteq E$ such that $x \in A$, $y \in B$, and $A \cup B = E$.
Of course, this implies that $d(x, y) \ge \eta$ for every $(x, y) \in
Z$, so that $Z \subseteq Z(\eta)$.  Note that $E$ is uniformly totally
separated with respect to the restriction of $d(\cdot, \cdot)$ to $E$,
as in Section \ref{uniform conditions}, if and only if $E$ is
uniformly totally separated along $Z(r)$ for each $r > 0$.

        Suppose that $E$ is strongly totally separated with respect to
the restriction of $d(\cdot, \cdot)$ to $E$.  If $Z \subseteq E \times
E$ is compact with respect to the corresponding product topology on $E
\times E$, and if $x \ne y$ for every $(x, y) \in Z$, then it is easy
to see that $Z \subseteq Z(r)$ for some $r > 0$.  Let us check that
$E$ is uniformly totally separated along $Z$ under these conditions.
Because $E$ is strongly totally separated, for each $(x, y) \in Z$
there are an $\eta > 0$ and $\eta$-separated sets $A, B \subseteq E$
such that $x \in A$, $y \in B$, and $A \cup B = E$.  Remember that $A$
and $B$ are relatively open subsets of $E$, so that $A \times B$ is
relatively open in $E \times E$.  Using the compactness of $Z$, we get
that there are finitely many positive real numbers $\eta_1, \ldots,
\eta_n$ and pairs of subsets $A_1, B_1, \ldots, A_n, B_n$ of $E$ such
that $A_j$, $B_j$ are $\eta_j$-separated and $A_j \cup B_j = E$
for each $j = 1, \ldots, n$, and
\begin{equation}
\label{Z subseteq bigcup_{j = 1}^n A_j times B_j}
        Z \subseteq \bigcup_{j = 1}^n A_j \times B_j.
\end{equation}
It follows that $E$ is uniformly totally separated along $Z$, with
\begin{equation}
\label{eta = min (eta_1, ldots, eta_n)}
        \eta = \min (\eta_1, \ldots, \eta_n).
\end{equation}
More precisely, if $(x, y) \in Z$, then $(x, y) \in A_j \times B_j$
for some $j$, and $A_j$, $B_j$ satisfy the requirements for $(x, y)$
needed to verify that $E$ is uniformly totally separated along $Z$.

        Let $x'$, $y'$ be distinct elements of $\overline{E}$, and let
$\{x_j\}_{j = 1}^\infty$, $\{y_j\}_{j = 1}^\infty$ be sequences of elements of $E$
that converge to $x'$, $y'$, respectively.  We may as well suppose that
$x_j \ne y_j$ for each $j$, and indeed that
\begin{equation}
\label{d(x_j, y_j) ge d(x', y')/2}
        d(x_j, y_j) \ge d(x', y')/2
\end{equation}
for each $j$, since these conditions will hold for all but finitely
many $j$ anyway.  Also let $Z$ be the set of these pairs $(x_j, y_j)$.
Suppose that $E$ is uniformly totally separated along $Z$, and let
$\eta$ be as in the definition of that property.  Thus for each
positive integer $j$ there are $\eta$-separated sets $A_j, B_j
\subseteq E$ such that $x_j \in A_j$, $y_j \in B_j$, and $A_j \cup B_j
= E$.  As before, this implies that $\overline{A_j}$ and
$\overline{B_j}$ are $\eta$-separated and satisfy
\begin{equation}
\label{overline{A_j} cup overline{B_j} = overline{E}}
        \overline{A_j} \cup \overline{B_j} = \overline{E}
\end{equation}
for each $j$.  We also have that
\begin{equation}
\label{B(x_j, eta) cap E subseteq A_j, B(y_j, eta) cap E subseteq B_j}
 B(x_j, \eta) \cap E \subseteq A_j, \quad B(y_j, \eta) \cap E \subseteq B_j
\end{equation}
for each $j$, as in (\ref{B(x, t_1) cap E subseteq A, B(y, t_2) cap E
  subseteq B}), and hence
\begin{equation}
\label{B(x_j, eta) cap overline{E} subseteq overline{A_j}, ...}
        B(x_j, \eta) \cap \overline{E} \subseteq \overline{A_j}, \quad
         B(y_j, \eta) \cap \overline{E} \subseteq \overline{B_j}
\end{equation}
for each $j$, as in (\ref{B(x, t_1) cap overline{E} subseteq
  overline{A}, ...}).  If $j$ is sufficiently large so that
\begin{equation}
\label{d(x', x_j), d(y', y_j) < eta}
        d(x', x_j), \ d(y', y_j) < \eta,
\end{equation}
then (\ref{B(x_j, eta) cap overline{E} subseteq overline{A_j}, ...})
implies that $x' \in \overline{A_j}$, $y' \in \overline{B_j}$.

        If $A, B \subseteq E$ are $\eta$-separated for some $\eta > 0$
and satisfy $A \cup B = E$, then $E$ is automatically uniformly
totally separated along $A \times B$, since one can use $A$, $B$ in
the definition of being uniformly totally separated along $A \times B$
for every $(x, y) \in A \times B$.  Similarly, if $A, B \subseteq
\overline{E}$ are $\eta$-separated for some $\eta > 0$ and satisfy $A
\cup B = \overline{E}$, then $\overline{E}$ is uniformly totally
separated along $A \times B$.  In this case, $A \cap E$ and $B \cap E$
are $\eta$-separated subsets of $E$ whose union is equal to $E$, and
$E$ is uniformly totally separated along $(A \cap E) \times (B \cap
E)$.  If $x' \in A$ and $y' \in B$, then there are sequences
$\{x_j\}_{j = 1}^\infty$ and $\{y_j\}_{j = 1}^\infty$ of elements of
$A \cap E$ and $B \cap E$ that converge to $x'$ and $y'$,
respectively.  In particular, $E$ is uniformly totally separated along
the set $Z$ of pairs $(x_j, y_j)$ under these conditions.

        Let us say that $E \subseteq M$ is \emph{uniformly $0$-dimensional
along $K \subseteq E$}\index{uniformly 0-dimensional spaces@uniformly
$0$-dimensional spaces} if for each $r > 0$ there is an $\eta = 
\eta(K, r) > 0$ such that for every $x \in K$ there are $\eta$-separated
sets $A, B \subseteq E$ with $x \in A$, $A \cup B = E$, and $A
\subseteq B(x, r)$.  Of course, this implies that $E$ is strongly
$0$-dimensional at each point in $K$.  Note that $E$ is uniformly
$0$-dimensional with respect to the restriction of $d(\cdot, \cdot)$
to $E$, as in Section \ref{uniform conditions}, if and only if $E$ is
uniformly $0$-dimensional along $K = E$.  If $K \subseteq E$ is
compact, and $E$ is strongly $0$-dimensional at each element of $K$,
then $E$ is uniformly $0$-dimensional along $K$.  This was already shown
in Section \ref{uniform conditions}, with slightly different terminology
and notation.

        Let $x' \in \overline{E}$ be given, and let $\{x_j\}_{j = 1}^\infty$
be a sequence of elements of $E$ that converges to $x'$.  Also let $K$
be the subset of $E$ consisting of the $x_j$'s, and suppose that $E$
is uniformly $0$-dimensional along $K$.  Under these conditions, one
can check that $\overline{E}$ is strongly $0$-dimensional at $x'$.
This is analogous to the arguments used earlier in this and the
previous section.  In particular, if $E$ is uniformly $0$-dimensional
with respect to the restriction of $d(\cdot, \cdot)$ to $E$, then
essentially the same type of argument shows that $\overline{E}$ is
uniformly $0$-dimensional as well.

        In the other direction, suppose that $\overline{E}$ is strongly
$0$-dimensional at a point $x' \in \overline{E}$.  Thus for each $r > 0$
there an $\eta > 0$ and $\eta$-separated sets $A, B \subseteq \overline{E}$
such that $x' \in A$, $A \cup B = \overline{E}$, and
\begin{equation}
\label{A subseteq B(x', r)}
        A \subseteq B(x', r).
\end{equation}
It follows that
\begin{equation}
\label{A subseteq B(x, 2^{1/q} r)}
        A \subseteq B(x, 2^{1/q} \, r)
\end{equation}
for every $x \in B(x', r)$, by the $q$-metric version of the triangle
inequality.  Of course, $A \cap E$ and $B \cap E$ are $\eta$-separated
subsets of $E$ whose union is equal to $E$.  Let $\{x_j\}_{j =
  1}^\infty$ be a sequence of elements of $E$ that converges to $x'$,
and let $K$ be the subset of $E$ consisting of the $x_j$'s, as before.
If $E$ is strongly $0$-dimensional at $x_j$ for each $j$, and
$\overline{E}$ is strongly $0$-dimensional at $x'$, then it is easy to
see that $E$ is uniformly $0$-dimensional along $K$.  One can verify
this directly, or using the fact that $\overline{E}$ is strongly
$0$-dimensional at $x_j$ for each $j$ too, as in Section \ref{some
  examples, remarks}.  This means that $\overline{E}$ is strongly
$0$-dimensional at each point in $K \cup \{x'\}$, which is a compact
set.  Hence $\overline{E}$ is uniformly $0$-dimensional along $K \cup
\{x'\}$, as before, which implies that $E$ is uniformly
$0$-dimensional along $K$.

\part{Simple functions}
\label{simple functions}

\section{Basic notions}
\label{basic notions}
\setcounter{equation}{0}

        Let $k$ be a field, and let $X$ be a nonempty set.  If $E$
is a subset of $X$, then we let ${\bf 1}_E(x)$ be the
\emph{characteristic}\index{characteristic functions} or
\emph{indicator function}\index{indicator functions} associated to $E$
on $X$, equal to $1$ when $x \in E$ and to $0$ when $x \in X \setminus
E$.  More precisely, ${\bf 1}_E(x)$ is considered here as a $k$-valued
function on $X$, so that $0$ and $1$ refer to the additive and
multiplicative identity elements in $k$.

        Let $V$ be a vector space over $k$, and let $f$ be a $V$-valued
\emph{simple function}\index{simple functions} on $X$, which is to say
a function on $X$ that takes only finitely many values in $V$.  Thus
$f$ can be expressed as
\begin{equation}
\label{f(x) = sum_{j = 1}^n v_j {bf 1}_{E_j}(x)}
        f(x) = \sum_{j = 1}^n v_j \, {\bf 1}_{E_j}(x),
\end{equation}
where $v_1, \ldots, v_n$ are the nonzero values of $f$, without
repetitions, and
\begin{equation}
\label{E_j = f^{-1}({v_j}) = {x in X : f(x) = v_j}}
        E_j = f^{-1}(\{v_j\}) = \{x \in X : f(x) = v_j\}
\end{equation}
for each $j = 1, \ldots, n$.  Note that the $E_j$'s are nonempty and
pairwise disjoint in this representation of $f$.  Conversely, if $v_1,
\ldots, v_n$ are finitely many vectors in $V$, and if $E_1, \ldots,
E_n$ are finitely many subsets of $X$, then (\ref{f(x) = sum_{j = 1}^n
  v_j {bf 1}_{E_j}(x)}) defines a $V$-valued simple function on $X$.
As usual, it is easy to reduce to the case where the $E_j$'s are
pairwise disjoint, using the various intersections of the $E_j$'s and
their complements in $X$.  One can also reduce to the case where the
$v_j$'s are nonzero and distinct, by combining the $E_j$'s as needed.
In particular, the space of $V$-valued simple functions on $X$ is a
vector space over $k$ with respect to pointwise addition and scalar
multiplication.

        Let $\mathcal{A}$ be an algebra of measurable subsets of $X$,
so that $\mathcal{A}$ is a collection of subsets of $X$ that contains
$\emptyset$, $X$ and is closed under finite unions, finite
intersections, and complementation.  A $V$-valued simple function $f$
on $X$ is said to be \emph{measurable}\index{measurable simple
  functions} with respect to $\mathcal{A}$ if
\begin{equation}
\label{f^{-1}({v}) in mathcal{A}}
        f^{-1}(\{v\}) \in \mathcal{A}
\end{equation}
for each $v \in V$.  Of course, $f^{-1}(\{v\})$ is the empty set for
all but finitely many $v \in V$ when $f$ is simple.  If $f$ is a
measurable $V$-valued simple function on $X$, then $f$ can be
expressed as in (\ref{f(x) = sum_{j = 1}^n v_j {bf 1}_{E_j}(x)}),
where $v_1, \ldots, v_n$ are the nonzero values of $f$ on $X$, without
repetitions, and (\ref{E_j = f^{-1}({v_j}) = {x in X : f(x) = v_j}})
is in $\mathcal{A}$ for each $j$.  Conversely, if $v_1, \ldots, v_n$
are finitely many vectors in $V$, and if $E_1, \ldots, E_n$ are
finitely many elements of $\mathcal{A}$, then (\ref{f(x) = sum_{j =
    1}^n v_j {bf 1}_{E_j}(x)}) defines a measurable $V$-valued simple
function on $V$.  This is clear when the $E_j$'s are pairwise
disjoint, and otherwise one can reduce to this case as in the
preceding paragraph.  It follows that the space of $V$-valued
measurable simple functions on $X$ is a linear subspace of the vector
space of all $V$-valued simple functions on $X$.  Let $S(X,
V)$\index{S(X, V)@$S(X, V)$} be the space of $V$-valued measurable
simple functions on $X$, which implicitly depends on the algebra
$\mathcal{A}$ too.

        It is easy to see that the product of two $k$-valued simple
functions on $X$ is a $k$-valued simple function on $X$, which is
measurable when the first two functions are measurable.  One can also
multiply a $k$-valued simple function on $X$ and a $V$-valued simple
function on $X$ to get another $V$-valued simple function on $X$,
which is measurable when the first two functions are measurable.  Note
that a real-valued simple function on $X$ is nonnegative at every
point in $X$ if and only if it can be expressed as a linear
combination of indicator functions with nonnegative coefficients.
Similarly, a real-valued measurable simple function on $X$ is
nonnegative on $X$ if and only if it can be expressed as a linear
combination of indicator functions of measurable sets with nonnegative
coefficients.  As usual, one can add and multiply such expressions, to
get another expression of the same type.

\section{Finitely-additive nonnegative measures}
\label{finitely-additive nonnegative measures}
\setcounter{equation}{0}

        Let $\mathcal{A}$ be an algebra of measurable subsets of a
nonempty set $X$ again, and let $\mu$ be a finitely-additive nonnegative
measure on $(X, \mathcal{A})$.  This means that $\mu$ is a nonnegative
extended-real-valued function on $\mathcal{A}$ which is finitely
additive on pairwise-disjoint measurable sets and satisfies
$\mu(\emptyset) = 0$.  If $f$ is a nonnegative real-valued measurable
simple function on $X$, then $f$ can be expressed as
\begin{equation}
\label{f(x) = sum_{j = 1}^n t_j {bf 1}_{E_j}(x)}
        f(x) = \sum_{j = 1}^n t_j \, {\bf 1}_{E_j}(x)
\end{equation}
for some nonnegative real numbers $t_1, \ldots, t_n$ and measurable
subsets $E_1, \ldots, E_n$ of $X$.  In this case, the integral of $f$
over $X$ with respect to $\mu$ is defined as a nonnegative extended
real number by
\begin{equation}
\label{int_X f d mu = sum_{j = 1}^n t_j mu(E_j)}
        \int_X f \, d\mu = \sum_{j = 1}^n t_i \, \mu(E_j),
\end{equation}
with the standard convention that $0 \cdot \infty = 0$.  As usual, one
can check that the value of the integral does not depend on the
particular representation (\ref{f(x) = sum_{j = 1}^n t_j {bf
    1}_{E_j}(x)}) of $\phi$.  If $a$ is a nonnegative real number,
then $a \, f(x)$ is also a nonnegative real-valued measurable simple
function on $X$, and
\begin{equation}
\label{int_X a f d mu = a int_X f d mu}
        \int_X a \, f \, d\mu = a \, \int_X f \, d\mu,
\end{equation}
using the convention $0 \cdot \infty = 0$ again, when necessary.
Similarly, if $g$ is another nonnegative real-valued measurable simple
function on $X$, then
\begin{equation}
\label{int_X (f + g) d mu = int_X f d mu + int_X g d mu}
 \int_X (f + g) \, d\mu = \int_X f \, d\mu + \int_X g \, d\mu.
\end{equation}
If $f(x) \le g(x)$ for every $x \in X$, then
\begin{equation}
\label{int_X f d mu le int_X g d mu}
        \int_X f \, d\mu \le \int_X g \, d\mu.
\end{equation}

        If $f$ is a nonnegative real-valued measurable simple function
on $X$, then $f(x)^r$ is a measurable simple function on $X$ for each
positive real number $r$, and we put
\begin{equation}
\label{||f||_r = ||f||_{L^r(X)} = (int_X f(x)^r d mu(x))^{1/r}}
        \|f\|_r = \|f\|_{L^r(X)} = \Big(\int_X f(x)^r \, d\mu(x)\Big)^{1/r}.
\end{equation}
Also put
\begin{equation}
\label{||f||_{L^infty(X)} = max {t ge 0 : mu(f^{-1}({t})) > 0}}
        \|f\|_{L^\infty(X)} = \max \{t \ge 0 : \mu(f^{-1}(\{t\})) > 0\},
\end{equation}
which may be described as the \emph{essential maximum}\index{essential
  maximum} of $f$ on $X$.  More precisely, the maximum on the right
side of (\ref{||f||_{L^infty(X)} = max {t ge 0 : mu(f^{-1}({t})) >
    0}}) is taken over all nonnegative real numbers $t$ such that
$\mu(f^{-1}(\{t\})) > 0$.  Because $f$ is a simple function,
$f^{-1}(\{t\}) = \emptyset$ for all but finitely many $t$, so that
$\mu(f^{-1}(\{t\})) = 0$ for all but finitely many $t$.  Thus the
right side of (\ref{||f||_{L^infty(X)} = max {t ge 0 : mu(f^{-1}({t}))
    > 0}}) reduces to the maximum of a finite set of nonnegative real
numbers.  This set is empty in the trivial situation where $\mu(X) =
0$, in which case we interpret (\ref{||f||_{L^infty(X)} = max {t ge 0
    : mu(f^{-1}({t})) > 0}}) as being equal to $0$.  Equivalently,
$\|f\|_{L^\infty(X)}$ is the smallest nonnegative real number such
that
\begin{equation}
\label{f(x) le ||f||_{L^infty(X)}}
        f(x) \le \|f\|_{L^\infty(X)}
\end{equation}
for almost every $x \in X$ with respect to $\mu$.

        Observe that
\begin{equation}
\label{||a f||_{L^r(X)} = a ||f||_{L^r(X)}}
        \|a \, f\|_{L^r(X)} = a \, \|f\|_{L^r(X)}
\end{equation}
for every nonnegative real number $a$ and $0 < r \le \infty$, using
the convention $0 \cdot \infty = 0$ again when needed.  If $g$ is
another nonnegative real-valued measurable simple function on $X$, then
\begin{equation}
\label{||f + g||_{L^r(X)} le ||f||_{L^r(X)} + ||g||_{L^r(X)}}
        \|f + g\|_{L^r(X)} \le \|f\|_{L^r(X)} + \|g\|_{L^r(X)}
\end{equation}
when $1 \le r \le \infty$.  This is version of Minkowski's
inequality,\index{Minkowski's inequality} which is straightforward
when $r = 1$ and $r = \infty$.  If $0 < r \le 1$, then
\begin{equation}
\label{(f(x) + g(x))^r le f(x)^r + g(x)^r, 2}
        (f(x) + g(x))^r \le f(x)^r + g(x)^r
\end{equation}
for every $x \in X$, as in (\ref{a^{q_2} + b^{q_2} le ... = (a^{q_1} +
  b^{q_1})^{q_2/q_1}}), with $q_1 = r$ and $q_2 = 1$.  This implies
that
\begin{equation}
\label{||f + g||_{L^r(X)}^r le ||f||_{L^r(X)}^r + ||g||_{L^r(X)}^r}
        \|f + g\|_{L^r(X)}^r \le \|f\|_{L^r(X)}^r + \|g\|_{L^r(X)}^r
\end{equation}
when $0 < r \le 1$, by integrating both sides of (\ref{(f(x) + g(x))^r
  le f(x)^r + g(x)^r, 2}) over $X$ with respect to $\mu$.

        Let us briefly consider the case where $\mathcal{A}$ is the
algebra of all subsets of a nonempty set $X$, and $\mu$ is
\emph{counting measure}\index{counting measure} on $X$.  Thus
$\mu(E)$ is defined to be the number of elements of $E \subseteq X$,
which is a nonnegative integer when $E$ has only finitely many elements,
and which is interpreted as being $+\infty$ when $E$ has infinitely
many elements.  If $f$ is any nonnegative real-valued function on $X$,
then $f$ is automatically measurable on $X$, and the Lebesgue integral
of $f$ with respect to counting measure on $X$ is the same as the sum
\begin{equation}
\label{sum_{x in X} f(x), 2}
        \sum_{x \in X} f(x),
\end{equation}
defined as the supremum of the corresponding finite subsums, as in
Section \ref{summable functions}.  If $f$ is a nonnegative real-valued
simple function on $X$, then this is consistent with (\ref{int_X f d
  mu = sum_{j = 1}^n t_j mu(E_j)}).  Note that (\ref{sum_{x in X}
  f(x), 2}) reduces to a finite sum when $f$ has finite support in
$X$, and that (\ref{sum_{x in X} f(x), 2}) is infinite when $f$ is a
nonnegative real-valued simple function on $X$ whose support has
infinitely many elements.

\section{Vector-valued functions}
\label{vector-valued functions}
\setcounter{equation}{0}

        Let $\mathcal{A}$ be an algebra of subsets of a nonempty set
$X$ again, and let $\mu$ be a finitely-additive nonnegative measure
on $(X, \mathcal{A})$.  Also let $|\cdot|$ be a $q$-absolute value
function on a field $k$ for some $q > 0$, and let $N$ be a $q$-norm on
a vector space $V$ over $k$ with respect to $|\cdot|$ on $k$.  If
$f(x)$ is a $V$-valued measurable simple function on $X$, then
$N(f(x))$ is a nonnegative real-valued measurable simple function on
$X$.  Thus we put
\begin{equation}
\label{||f||_r = ||f||_{L^r(X, V)} = (int_X N(f(x))^r d mu(x))^{1/r}}
 \|f\|_r = \|f\|_{L^r(X, V)} = \Big(\int_X N(f(x))^r \, d\mu(x)\Big)^{1/r}
\end{equation}
for every positive real number $r$, which is defined as an extended
real number, as in the previous section.  Similarly, we take
$\|f\|_{L^\infty(X, V)}$ to be the essential maximum of $N(f(x))$ on
$X$, as in (\ref{||f||_{L^infty(X)} = max {t ge 0 : mu(f^{-1}({t})) >
    0}}).  This is the same as saying that $\|f\|_{L^\infty(X, V)}$ is
the smallest nonnegative real number such that
\begin{equation}
\label{N(f(x)) le ||f||_{L^infty(X, V)}}
        N(f(x)) \le \|f\|_{L^\infty(X, V)}
\end{equation}
for almost every $x \in X$ with respect to $\mu$.  Note that
\begin{equation}
\label{||a f||_{L^r(X, V)} = |a| ||f||_{L^r(X, V)}}
        \|a \, f\|_{L^r(X, V)} = |a| \, \|f\|_{L^r(X, V)}
\end{equation}
for every $a \in k$ and $0 < r \le \infty$, with the usual convention
that $0 \cdot \infty = 0$.

        If $g(x)$ is another $V$-valued measurable simple function on
$X$, then we have that
\begin{equation}
\label{N(f(x) + g(x))^r le (N(f(x))^q + N(g(x))^q)^{r/q}, 2}
        N(f(x) + g(x))^r \le (N(f(x))^q + N(g(x))^q)^{r/q}
\end{equation}
for every $x \in X$ and positive real number $r$, by the $q$-norm
version of the triangle inequality for $N$ on $V$.  If $r \le q$, then
it follows that
\begin{equation}
\label{N(f(x) + g(x))^r le N(f(x))^r + N(g(x))^r, 2}
        N(f(x) + g(x))^r \le N(f(x))^r + N(g(x))^r
\end{equation}
for every $x \in X$, as in (\ref{(f(x) + g(x))^r le f(x)^r + g(x)^r, 2}),
with $r$ replaced by $r/q$.  This implies that
\begin{equation}
\label{||f + g||_{L^r(X, V)}^r le ||f||_{L^r(X, V)}^r + ||g||_{L^r(X, V)}^r}
 \|f + g\|_{L^r(X, V)}^r \le \|f\|_{L^r(X, V)}^r + \|g\|_{L^r(X, V)}^r
\end{equation}
when $0 < r \le q$, by integrating both sides of (\ref{N(f(x) +
  g(x))^r le N(f(x))^r + N(g(x))^r}) with respect to $\mu$ on $X$.  If
$q \le r < \infty$, then we get that
\begin{equation}
\label{||f + g||_{L^r(X, V)}^q le ||f||_{L^r(X, V)}^q + ||g||_{L^r(X, V)}^q}
 \|f + g\|_{L^r(X, V)}^q \le \|f\|_{L^r(X, V)}^q + \|g\|_{L^r(X, V)}^q,
\end{equation}
using (\ref{N(f(x) + g(x))^r le (N(f(x))^q + N(g(x))^q)^{r/q}, 2}) and
(\ref{||f + g||_{L^r(X)} le ||f||_{L^r(X)} + ||g||_{L^r(X)}}), with
$r$ replaced by $r/q$ in the latter.  It is easy to check directly
that (\ref{||f + g||_{L^r(X, V)}^q le ||f||_{L^r(X, V)}^q +
  ||g||_{L^r(X, V)}^q}) holds when $r = \infty$, using (\ref{N(f(x) +
  g(x))^r le (N(f(x))^q + N(g(x))^q)^{r/q}, 2}) with $r = q$.  If $N$
is an ultranorm on $V$, then (\ref{N(f(x) + g(x))^r le N(f(x))^r +
  N(g(x))^r, 2}) holds for every positive real number $r$, which
implies that (\ref{||f + g||_{L^r(X, V)}^r le ||f||_{L^r(X, V)}^r +
  ||g||_{L^r(X, V)}^r}) holds when $0 < r < \infty$, as before.  This
corresponds to $q = \infty$, in which case we also have that
\begin{equation}
\label{||f + g||_{L^infty(X, V)} le ...}
 \|f + g\|_{L^\infty(X, V)} \le \max(\|f\|_{L^\infty(X, V)}, \|g\|_{L^\infty(X, V)}),
\end{equation}
as one can easily verify.

        Note that
\begin{equation}
\label{{x in X : f(x) ne 0} = X setminus f^{-1}({0})}
        \{x \in X : f(x) \ne 0\} = X \setminus f^{-1}(\{0\})
\end{equation}
is a measurable subset of $X$ when $f$ is a $V$-valued measurable
simple function on $X$.  If $r$ is any positive real number, then
\begin{equation}
\label{||f||_{L^r(X, V)} < infty}
        \|f\|_{L^r(X, V)} < \infty
\end{equation}
if and only if
\begin{equation}
\label{mu({x in X : f(x) ne 0}) < infty}
        \mu(\{x \in X : f(x) \ne 0\}) < \infty.
\end{equation}
Let $S_0(X, V)$\index{S_0(X, V)@$S_0(X, V)$} be the space of
$V$-valued measurable simple functions on $X$ that satisfy (\ref{mu({x
    in X : f(x) ne 0}) < infty}), which is a linear subspace of the
space $S(X, V)$ of all $V$-valued measurable simple functions on $X$.
Of course, (\ref{||f||_{L^r(X, V)} < infty}) holds for every
$V$-valued measurable simple function on $X$ when $r = \infty$.
Similarly, if $0 < r \le \infty$, then
\begin{equation}
\label{||f||_{L^r(X, V)} = 0}
        \|f\|_{L^r(X, V)} = 0
\end{equation}
if and only if
\begin{equation}
\label{mu({x in X : f(x) ne 0}) = 0}
        \mu(\{x \in X : f(x) \ne 0\}) = 0.
\end{equation}
Let us say that $V$-valued measurable simple functions $f$, $g$ on $X$
are \emph{equivalent}\index{equivalent measurable functions} when
\begin{equation}
\label{mu({x in X : f(x) ne g(x)}) = 0}
        \mu(\{x \in X : f(x) \ne g(x)\}) = 0,
\end{equation}
in which case
\begin{equation}
\label{||f||_{L^r(X, V)} = ||g||_{L^r(X, V)}}
        \|f\|_{L^r(X, V)} = \|g\|_{L^r(X, V)}
\end{equation}
for every $0 < r \le \infty$.  This defines an equivalence relation on
$S(X, V)$, and we let $\widetilde{S}(X, V)$\index{S^~(X,
  V)@$\widetilde{S}(X, V)$} be the corresponding space of equivalence
classes.  This is the same as taking the quotient of $S(X, V)$ by the
linear subspace consisting of functions equal to $0$ almost everywhere
on $X$ with respect to $\mu$.  In particular, $\widetilde{S}(X, V)$ is
also a vector space over $k$ in a natural way, and it is easy to see
that $\|f\|_{L^\infty(X, V)}$ determines a well-defined $q$-norm on
$\widetilde{S}(X, V)$.  Let $\widetilde{S}_0(X, V)$\index{S_0^~(X,
  V)@$\widetilde{S}_0(X, V)$} be the image of $S_0(X, V)$ in
$\widetilde{S}(X, V)$ under this quotient mapping, which consists of
equivalence classes of $V$-valued measurable simple functions on $X$
that satisfy (\ref{mu({x in X : f(x) ne 0}) < infty}).  Thus
$\widetilde{S}_0(X, V)$ is a linear subspace of $\widetilde{S}(X, V)$,
since $S_0(X, V)$ is a linear subspace of $S(X, V)$.  One can check
that $\|f\|_{L^r(X, V)}$ determines a well-defined $q$-norm on
$\widetilde{S}_0(X, V)$ when $q \le r < \infty$, and that
$\|f\|_{L^r(X, V)}$ determines a well-defined $r$-norm on
$\widetilde{S}_0(X, V)$ when $0 < r \le q$, by (\ref{||a f||_{L^r(X,
    V)} = |a| ||f||_{L^r(X, V)}}), (\ref{||f + g||_{L^r(X, V)}^r le
  ||f||_{L^r(X, V)}^r + ||g||_{L^r(X, V)}^r}), and (\ref{||f +
  g||_{L^r(X, V)}^q le ||f||_{L^r(X, V)}^q + ||g||_{L^r(X, V)}^q}).

        Suppose for the moment that $\mathcal{A}$ is the algebra of
all subsets of a nonempty set $X$, and that $\mu$ is counting measure
on $X$.  In this case, every $V$-valued simple function on $X$ is
automatically measurable, and $S_0(X, V)$ is the same as the space
$c_{00}(X, V)$ of $V$-valued functions on $X$ with finite support.
Any two functions on $X$ that are equal almost everywhere with respect
to counting measure are in fact equal everywhere on $X$, so that
$\widetilde{S}(X, V)$ is the same as $S(X, V)$, and
$\widetilde{S}_0(X, V)$ is the same as $S_0(X, V) = c_{00}(X, V)$.  If
$f$ is a $V$-valued simple function on $X$, then $\|f\|_{L^\infty(X,
  V)}$ is the same as the supremum norm of $f$, as in Section
\ref{supremum metrics, norms}.  If $r$ is a positive real number and
$f \in S_0(X, V)$, then $\|f\|_{L^r(X, V)}$ is the same as
$\|f\|_{\ell^r(X, V)}$, as in Section \ref{l^r norms}.  If $f$ is a
$V$-valued simple function on $X$ not in $S_0(X, V)$, then
$\|f\|_{L^r(X, V)} = +\infty$, but $\|f\|_{\ell^r(X, V)}$ was not
defined in Section \ref{l^r norms}, strictly speaking.  However, if
$\|f\|_{\ell^r(X, V)}$ were defined in the same way as in Section
\ref{l^r norms}, then it would also be infinite under these
conditions.

\section{The unit interval}
\label{unit interval}
\setcounter{equation}{0}

        Let us consider the case where $X$ is the closed unit interval
$[0, 1]$, and $\mathcal{A}$ is an algebra of subsets of $[0, 1]$ that
includes all closed subintervals of $[0, 1]$.  In particular, a subset
of $[0, 1]$ with only one element is considered as a closed interval
of length $0$, and hence should be in $\mathcal{A}$.  Also let $\mu$
be a finitely-additive nonnegative measure on $(X, \mathcal{A})$ such
that $\mu([0, 1]) < \infty$,
\begin{equation}
\label{lim_{t to a+} mu([a, t]) = 0}
        \lim_{t \to a+} \mu([a, t]) = 0
\end{equation}
when $0 \le a < 1$, and
\begin{equation}
\label{lim_{t to b-} mu([t, b]) = 0}
        \lim_{t \to b-} \mu([t, b]) = 0
\end{equation}
when $0 < b \le 1$.  This implies that
\begin{equation}
\label{mu({a}) = 0}
        \mu(\{a\}) = 0
\end{equation}
for every $a \in [0, 1]$, and in some situations (\ref{lim_{t to
    a+} mu([a, t]) = 0}) and (\ref{lim_{t to b-} mu([t, b]) =
  0}) can be derived from (\ref{mu({a}) = 0}).  Of course, these
conditions hold when
\begin{equation}
\label{mu([a, b]) = b - a}
        \mu([a, b]) = b - a
\end{equation}
for $0 \le a \le b \le 1$.

        Let $k$ be a field with a $q$-absolute value function for some
$q > 0$, and $V$ be a vector space over $k$ with a $q$-norm $N$, as
before.  If $f$ is a $V$-valued measurable simple function on $[0,
  1]$, then put
\begin{equation}
\label{f_t(x) = {bf 1}_{[0, t]}(x) f(x)}
        f_t(x) = {\bf 1}_{[0, t]}(x) \, f(x)
\end{equation}
for every $t, x \in [0, 1]$, where ${\bf 1}_{[0, t]}$ is considered as
a $k$-valued indicator function on $[0, 1]$.  Thus $f_t(x)$ is a
$V$-valued measurable simple function on $[0, 1]$ as a function of $x$
for each $t \in [0, 1]$.  By construction, $f_1 = f$, and $f_0(x) = 0$
for every $x \in (0, 1]$.  This means that $f_0 = 0$ almost everywhere
on $[0, 1]$ with respect to $\mu$, by (\ref{mu({a}) = 0}) with $a =
0$, and one could also change the definitions slightly to get $f_0 =
0$ everywhere on $X$.  If $0 \le t_1 \le t_2 \le 1$, then
\begin{equation}
\label{f_{t_2}(x) - f_{t_1}(x) = {bf 1}_{(t_1, t_2]}(x) f(x)}
        f_{t_2}(x) - f_{t_1}(x) = {\bf 1}_{(t_1, t_2]}(x) \, f(x)
\end{equation}
for every $x \in [0, 1]$, so that
\begin{equation}
\label{||f_{t_2} - f_{t_1}||_{L^r([0, 1], V)} = ...}
        \|f_{t_2} - f_{t_1}\|_{L^r([0, 1], V)}
                  = \Big(\int_{(t_1, t_2]} N(f(x))^r \, d\mu(x)\Big)^{1/r}
\end{equation}
for every positive real number $r$.  In particular,
\begin{equation}
\label{||f_{t_2} - f_{t_1}||_{L^r([0, 1], V)} le ...}
        \|f_{t_2} - f_{t_1}\|_{L^r([0, 1], V)}
                      \le \mu((t_1, t_2])^{1/r} \, \|f\|_{L^\infty([0, 1], V)}
\end{equation}
for every $0 \le t_1 \le t_2 \le 1$ and $r > 0$.

        As in the preceding section, we can identify $V$-valued measurable
simple functions on $[0, 1]$ that are equal almost everywhere with
respect to $\mu$, to get a vector space $\widetilde{S}([0, 1], V)$
over $k$.  We have seen that $\|\cdot\|_{L^r([0, 1], V)}$ determines a
$q$-norm on $\widetilde{S}([0, 1], V)$ when $r \ge q$, and an $r$-norm
on $\widetilde{S}([0, 1], V)$ when $0 < r \le q$.  This leads to a
$q$-metric on $\widetilde{S}([0, 1], V)$ when $r \ge q$, and to an
$r$-metric on $\widetilde{S}([0, 1], V)$ when $0 < r \le q$, as usual.
In both cases, we get a topology on $\widetilde{S}([0, 1], V)$
corresponding to $\|\cdot\|_{L^r(X, V)}$.  It follows from
(\ref{lim_{t to a+} mu([a, t]) = 0}), (\ref{lim_{t to b-} mu([t, b]) =
  0}), and (\ref{||f_{t_2} - f_{t_1}||_{L^r([0, 1], V)} le ...})  that
\begin{equation}
\label{t mapsto f_t}
        t \mapsto f_t
\end{equation}
leads to a continuous mapping from $[0, 1]$ into $\widetilde{S}([0,
  1], V)$ with respect to the topology on $\widetilde{S}([0, 1], V)$
corresponding to $\|\cdot\|_{L^r([0, 1], V)}$ when $0 < r < \infty$.
This shows that $\widetilde{S}([0, 1], V)$ is pathwise connected with
respect to this topology when $0 < r < \infty$.  Essentially the same
argument shows that $\widetilde{S}([0, 1], V)$ is contractible with
respect to this topology when $0 < r < \infty$.  If $g$ is any other
$V$-valued measurable simple function on $[0, 1]$, then one can
consider families of the form $f_t + g$, in order to get
contractibility centered at $g$ instead of $0$.  Note that (\ref{t
  mapsto f_t}) is not normally continuous with respect to the topology
on $\widetilde{S}([0, 1], V)$ corresponding to
$\|\cdot\|_{L^\infty([0, 1], V)}$.

        Of course, one can get contractibility of vector spaces over
the real or complex numbers using scalar multiplication.  In
particular, if $k = {\bf R}$ or ${\bf C}$ equipped with the standard
absolute value function, then one can use this to get contractibility
of $\widetilde{S}([0, 1], V)$ with respect to the topology
corresponding to $\|\cdot\|_{L^\infty([0, 1], V)}$.  Otherwise, if $k$
is a field equipped with an ultrametric absolute value function, and
$N$ is an ultranorm on $V$, then $\|\cdot\|_{L^\infty([0, 1], V)}$
satisfies the ultrametric version of the triangle inequality, as in
(\ref{||f + g||_{L^infty(X, V)} le ...}).  This means that
$\|\cdot\|_{L^\infty([0, 1], V)}$ determines an ultranorm on
$\widetilde{S}([0, 1], V)$, so that the $\widetilde{S}([0, 1], V)$ is
uniformly $0$-dimensional with respect to the corresponding
ultrametric, as in Section \ref{uniform conditions}.

\section{Pushing measures forward}
\label{pushing measures forward}
\setcounter{equation}{0}

        Let $X$, $Y$ be nonempty sets, and let $\mathcal{A}$, $\mathcal{B}$
be algebras of subsets of $X$, $Y$, respectively.  Suppose that a
mapping $\phi : X \to Y$ is \emph{measurable}\index{measurable
  mappings} in the sense that $\phi^{-1}(E) \in \mathcal{A}$ for every
$E \in \mathcal{B}$.  If $\mu$ is a finitely-additive nonnegative measure
on $(X, \mathcal{A})$, then it is easy to see that
\begin{equation}
\label{nu(E) = mu(phi^{-1}(E))}
        \nu(E) = \mu(\phi^{-1}(E))
\end{equation}
defines a finitely-additive nonnegative measure on $(Y, \mathcal{B})$.
This is the measure on $Y$ obtained by \emph{pushing $\mu$
  forward}\index{pushing measures forward} using $\phi$.  If $f$ is a
nonnegative real-valued simple function on $Y$ that is measurable with
respect to $\mathcal{B}$, then one can check that $f \circ \phi$ is a
nonnegative real-valued simple function on $X$ that is measurable with
respect to $\mathcal{A}$, and that
\begin{equation}
\label{int_X f circ phi d mu = int_Y f d nu}
        \int_X f \circ \phi \, d\mu = \int_Y f \, d\nu
\end{equation}
under these conditions.

        Let $V$ be a vector space over a field $k$ again.  If $f$ is a
$V$-valued simple function on $Y$ that is measurable with respect to
$\mathcal{B}$, then it is easy to see that $f \circ \phi$ is a $V$-valued
simple function on $X$ that is measurable with respect to $\mathcal{A}$.
Suppose that $|\cdot|$ is a $q$-absolute value function on $k$ for some
$q > 0$, and that $N$ is a $q$-norm on $V$ with respect to $|\cdot|$ on $k$.
If $f$ is as before, then $N(f(x))$ is a nonnegative real-valued measurable
simple function on $Y$, and $N(f(\phi(x)))$ is a nonnegative real-valued
measurable simple function on $X$.  Observe that
\begin{equation}
\label{||f circ phi||_{L^r(X, V)} = ||f||_{L^r(Y, V)}}
        \|f \circ \phi\|_{L^r(X, V)} = \|f\|_{L^r(Y, V)}
\end{equation}
for every $0 < r \le \infty$, using (\ref{int_X f circ phi d mu =
  int_Y f d nu}) when $r < \infty$, and going back to the definition
of the essential maximum when $r = \infty$.

        Suppose for the moment that $Y = [0, 1]$, and that $\mathcal{B}$
includes all closed subintervals of $[0, 1]$, as in the previous
section.  Suppose also that
\begin{equation}
        \nu([0, 1]) = \mu(X) < \infty,
\end{equation}
and that $\nu$ satisfies the analogues of (\ref{lim_{t to a+} mu([a,
    t]) = 0}) and (\ref{lim_{t to b-} mu([t, b]) = 0}) in this
situation, and hence (\ref{mu({a}) = 0}).  Let $f$ be a measurable
$V$-valued simple function on $X$, and put
\begin{equation}
\label{f_t(x) = {bf 1}_{[0, t]}(phi(x)) f(x)}
        f_t(x) = {\bf 1}_{[0, t]}(\phi(x)) \, f(x)
\end{equation}
for every $x \in X$ and $0 \le t \le 1$, where ${\bf 1}_{[0, t]}$ is
the indicator function on $[0, 1]$ associated to $[0, t]$.
Equivalently,
\begin{equation}
\label{f_t(x) = {bf 1}_{phi^{-1}([0, t])}(x) f(x)}
        f_t(x) = {\bf 1}_{\phi^{-1}([0, t])}(x) \, f(x)
\end{equation}
for every $x \in X$ and $0 \le t \le 1$, where ${\bf 1}_{\phi^{-1}([0,
    t])}$ is the indicator function on $X$ associated to
$\phi^{-1}([0, t])$.  Thus $f_t(x)$ is a $V$-valued measurable simple
function of $x$ on $X$ for every $t \in [0, 1]$, $f_1 = f$, and
$f_0(x) = 0$ when $\phi(x) \ne 0$.  By hypothesis, $\{0\}$ is a
measurable subset of $Y$ with respect to $\mathcal{B}$, so that
$\phi^{-1}(\{0\})$ is a measurable subset of $X$ with respect to
$\mathcal{A}$, and
\begin{equation}
\label{mu(phi^{-1}({0})) = nu({0}) = 0}
        \mu(\phi^{-1}(\{0\})) = \nu(\{0\}) = 0.
\end{equation}
This implies that $f_0 = 0$ almost everywhere on $X$ with respect to
$\mu$, although one could again make some changes to get $f_0 = 0$
everywhere on $X$, if desired.  If $0 \le t_1 \le t_2 \le 1$, then we
have that
\begin{equation}
\label{f_{t_2}(x) - f_{t_1}(x) = {bf 1}_{(t_1, t_2]}(phi(x)) f(x)}
        f_{t_2}(x) - f_{t_1}(x) = {\bf 1}_{(t_1, t_2]}(\phi(x)) \, f(x)
\end{equation}
for every $x \in X$, and hence
\begin{equation}
\label{||f_{t_2} - f_{t_1}||_{L^r(X, V)} = ...}
        \|f_{t_2} - f_{t_1}\|_{L^r(X, V)}
            = \Big(\int_{\phi^{-1}((t_1, t_2])} N(f(x))^r \, d\mu(x)\Big)^{1/r}
\end{equation}
for every positive real number $r$.  It follows that
\begin{eqnarray}
\label{||f_{t_2} - f_{t_1}||_{L^r(X, V)} le ...}
 \|f_{t_2} - f_{t_1}\|_{L^r(X, V)} & \le & \mu(\phi^{-1}((t_1, t_2])^{1/r} \,
                                         \|f\|_{L^\infty(X, V)} \\
 & = & \nu((t_1, t_2])^{1/r} \, \|f\|_{L^\infty(X, V)} \nonumber
\end{eqnarray}
when $0 \le t_1 \le t_2 \le 1$ and $0 < r < \infty$.  As in Section
\ref{vector-valued functions}, we can identify $V$-valued measurable
simple functions on $X$ that are equal almost everywhere with respect
to $\mu$ to get a vector space $\widetilde{S}(X, V)$ over $k$, and
$\|\cdot\|_{L^r(X, V)}$ determines a $q$-norm on $\widetilde{S}(X, V)$
when $r \ge q$, and an $r$-norm on $\widetilde{S}(X, V)$ when $0 < r
\le q$.  Using (\ref{||f_{t_2} - f_{t_1}||_{L^r(X, V)} le ...}) and
the analogues of (\ref{lim_{t to a+} mu([a, t]) = 0}) and (\ref{lim_{t
    to b-} mu([t, b]) = 0}) for $\nu$, we get that
\begin{equation}
\label{t mapsto f_t, 2}
        t \mapsto f_t
\end{equation}
leads to a continuous mapping from $[0, 1]$ into $\widetilde{S}(X, V)$
with respect to the topology on $\widetilde{S}(X, V)$ corresponding to
$\|\cdot\|_{L^r(X, V)}$ when $0 < r < \infty$.  As before, this
implies that $\widetilde{S}(X, V)$ is pathwise connected with respect
to this topology when $0 < r < \infty$, and in fact contractible.  One
can also consider $f_t + g$ for any other $V$-valued measurable simple
function $g$ on $X$, to get contractibility centered at $g$ instead of
$0$.

        Suppose now that $X$, $Y$ are topological spaces, and that
$\mathcal{A}$, $\mathcal{B}$ are the $\sigma$-aglebras of Borel subsets
of $X$, $Y$, respectively.  If $\phi$ is a continuous mapping from $X$
into $Y$, then $\phi$ is automatically measurable with respect to the
Borel sets.  In particular, we can apply the discussion in the
preceding paragraph to the case where $Y = [0, 1]$ with the standard
topology.  It is well known that there are continuous mappings $\phi$
from topological Cantor sets $X$ onto $[0, 1]$, and that one can do
this in such a way that Lebesgue measure on $[0, 1]$ corresponds to
pushing forward a finite nonnegative Borel measure $\mu$ on $X$.  Thus
connectedness of $X$ as a topological space is not really needed here.

\section{Countability conditions}
\label{countability conditions}
\setcounter{equation}{0}

        Let $X$ be a nonempty set, let $\mathcal{A}$ be an algebra of
subsets of $X$, and let $\mu$ be a finitely-additive nonnegative
measure on $(X, \mathcal{A})$.  Also let $k$ be a field with a
$q$-absolute value function $|\cdot|$ for some $q > 0$, and let $N$ be
a $q$-norm on $V$ with respect to $|\cdot|$ on $k$.  Suppose for the
moment that $\mu$ takes values in a set of finitely or countably many
nonnegative extended real numbers, and that $N$ takes values in a set
of finitely or countably many nonnegative real numbers too.  In
particular, these conditions hold when $\mathcal{A}$ has only finitely
or countably many elements, and $V$ has only finitely or countably
many elements.  This implies that $\|f\|_{L^\infty(X, V)}$ takes
values in a set of finitely or countably many nonnegative real numbers
when $f$ is a $V$-valued measurable simple function on $X$.  As in
Section \ref{vector-valued functions}, we can identify $V$-valued
measurable simple functions on $X$ that are equal almost everywhere
with respect to $\mu$, to get a vector space $\widetilde{S}(X, V)$
over $k$, and $\|\cdot\|_{L^\infty(X, V)}$ determines a $q$-norm on
$\widetilde{S}(X, V)$.  Under the conditions just described, this
$q$-norm takes values in a set of only finitely or countably many
nonnegative real numbers, which implies that $\widetilde{S}(X, V)$ has
topological dimension $0$ with respect to topology determined by the
corresponding $q$-metric.  Of course, if $N$ is an ultranorm on $V$,
then $\|\cdot\|_{L^\infty(X, V)}$ determines an ultranorm on
$\widetilde{S}(X, V)$, and one does not need any countability
conditions to get that $\widetilde{S}(X, V)$ has topological dimension
$0$.

        If $\mu$ and $N$ take values in sets with only finitely or
countably many values again, then for each positive real number $r$,
$\|f\|_{L^r(X, V)}$ takes values in a set of only finitely or
countably many nonnegative extended real numbers when $f$ is a
$V$-valued measurable simple function on $X$.  Let us now restrict our
attention to $V$-valued measurable simple functions $f$ on $X$ that
satisfy (\ref{mu({x in X : f(x) ne 0}) < infty}), and hence
(\ref{||f||_{L^r(X, V)} < infty}), which corresponds to the linear
subspace $S_0(X, V)$ of $S(X, V)$.  This leads to a linear subspace
$\widetilde{S}_0(X, V)$ of $\widetilde{S}(X, V)$, after identifying
such functions that are equal almost everywhere with respect to $\mu$,
as in Section \ref{vector-valued functions}.  Remember that
$\|\cdot\|_{L^r(X, V)}$ determines a $q$-norm on $\widetilde{S}_0(X,
V)$ when $r \ge q$, and an $r$-norm when $0 < r \le q$.  If $\mu$ and
$N$ take values in sets with only finitely or countably many elements,
then $\widetilde{S}_0(X, V)$ has topological dimension $0$ with
respect to the $q$ or $r$-metric associated to $\|\cdot\|_{L^r(X, V)}$
for each $r > 0$, for the same reasons as before.

        Suppose for the rest of the section that $|\cdot|$ is an
ultrametric absolute value function on $k$, and that $N$ is an
ultranorm on $V$ with respect to $|\cdot|$ on $k$.  If $N$ does not
already take values in a set of finitely or countably many nonnegative
real numbers, then we can basically reduce to that case by modifying
$N$, as in Section \ref{some variants}.  As before, let $h(t)$ be a
monotonically increasing real-valued function defined on the set of
nonnegative real numbers such that $h(0) = 0$ and $h(t) > 0$ when $t >
0$.  Thus
\begin{equation}
\label{h(N(v - w)), 2}
        h(N(v - w))
\end{equation}
defines an ultrametric on $V$ which determines the same topology on
$V$ as the ultrametric $N(v - w)$ associated to $N$.  If $f$, $g$ are
$V$-valued measurable simple functions on $X$, then $h(N(f(x) -
g(x)))^r$ is a nonnegative real-valued measurable simple function on
$X$ for every positive real number $r$, and we put
\begin{equation}
\label{d_r(f, g) = (int_X h(N(f(x) - g(x)))^r d mu(x))^{1/r}}
        d_r(f, g) = \Big(\int_X h(N(f(x) - g(x)))^r \, d\mu(x)\Big)^{1/r}.
\end{equation}
It is easy to see that (\ref{d_r(f, g) = (int_X h(N(f(x) - g(x)))^r d
  mu(x))^{1/r}}) satisfies the $r$-metric version of the triangle
inequality for each $r > 0$ under these conditions, for essentially
the same reasons as for
\begin{equation}
\label{||f - g||_{L^r(X, V)}}
        \|f - g\|_{L^r(X, V)},
\end{equation}
as in Section \ref{vector-valued functions}.  If $f$, $g$ also satisfy
(\ref{mu({x in X : f(x) ne 0}) < infty}), which is to say that $f, g
\in S_0(X, V)$, then $f - g \in S_0(X, V)$ too, which implies that
(\ref{d_r(f, g) = (int_X h(N(f(x) - g(x)))^r d mu(x))^{1/r}}) is
finite.  As usual, we get a vector space $\widetilde{S}_0(X, V)$ by
identifying $V$-valued measurable simple functions on $X$ that satisfy
(\ref{mu({x in X : f(x) ne 0}) < infty}) and which are equal almost
everywhere with respect to $\mu$, and (\ref{d_r(f, g) = (int_X
  h(N(f(x) - g(x)))^r d mu(x))^{1/r}}), (\ref{||f - g||_{L^r(X, V)}})
define $r$-metrics on $\widetilde{S}_0(X, V)$.  If $h(t)$ and $t$ are
each bounded by positive constant multiples of the other on $[0,
  +\infty)$, then (\ref{d_r(f, g) = (int_X h(N(f(x) - g(x)))^r d
    mu(x))^{1/r}}) and (\ref{||f - g||_{L^r(X, V)}}) are each bounded
  by the same constant multiples of the other on $S_0(X, V)$, and so
  the corresponding $r$-metrics on $\widetilde{S}_0(X, V)$ are bounded
  by the same constant multiples of each other.  This implies that the
  corresponding $r$-metrics determine the same topology on
  $\widetilde{S}_0(X, V)$.

        As in Section \ref{some variants}, we can choose $h(t)$ so that
it takes values in a countable subset of ${\bf R}$, in addition to the
properties already mentioned.  If $\mu$ takes values in a set of
finitely or countably many nonnegative extended real numbers, then it
follows that for each positive real number $r$, (\ref{d_r(f, g) =
  (int_X h(N(f(x) - g(x)))^r d mu(x))^{1/r}}) takes values in a set of
finitely or countably many nonnegative real numbers when $f, g \in
S_0(X, V)$.  This means that for each $0 < r < \infty$, the $r$-metric
on $\widetilde{S}_0(X, V)$ corresponding to (\ref{d_r(f, g) = (int_X
  h(N(f(x) - g(x)))^r d mu(x))^{1/r}}) takes values in the same set of
finitely or countably many nonnegative real numbers.  Under these
conditions, we get that $\widetilde{S}_0(X, V)$ has topological
dimension $0$ with respect to the topology determined by the
$r$-metrics corresponding to (\ref{d_r(f, g) = (int_X h(N(f(x) -
  g(x)))^r d mu(x))^{1/r}}) or (\ref{||f - g||_{L^r(X, V)}}) when $0 <
r < \infty$.

\section{Measurable sets}
\label{measurable sets}
\setcounter{equation}{0}

        The \emph{symmetric difference}\index{symmetric difference}
of two sets $A$, $B$ is the set
\begin{equation}
\label{A bigtriangleup B = (A setminus B) cup (B setminus A)}
 A \bigtriangleup B = (A \setminus B) \cup (B \setminus A).
\end{equation}
If $C$ is another set, then
\begin{eqnarray}
\label{A bigtriangleup C = ...}
 A \bigtriangleup C & = & (A \setminus C) \cup (C \setminus A) \\
 & \subseteq & ((A \setminus B) \cup (B \setminus C))
           \cup ((C \setminus B) \cup (B \setminus C)) \nonumber \\
 & = & (A \bigtriangleup B) \cup (B \bigtriangleup C). \nonumber
\end{eqnarray}
Similarly, observe that
\begin{equation}
\label{(A cap C) bigtriangleup (B cap C) = (A bigtriangleup B) cap C}
        (A \cap C) \bigtriangleup (B \cap C) = (A \bigtriangleup B) \cap C
\end{equation}
and
\begin{equation}
\label{(A cup C) bigtriangleup (B cup C) = ...}
        (A \cup C) \bigtriangleup (B \cup C)
               = (A \setminus C) \bigtriangleup (B \setminus C)
               = (A \bigtriangleup B) \setminus C.
\end{equation}
Let $X$ be a nonempty set, let $\mathcal{A}$ be an algebra of subsets
of $X$, and let $\mu$ be a finitely-additive nonnegative measure on
$(X, \mathcal{A})$.  If $A, B \in \mathcal{A}$, then $A \bigtriangleup
B \in \mathcal{A}$, and we put
\begin{equation}
\label{d_mu(A, B) = mu(A bigtriangleup B)}
        d_\mu(A, B) = \mu(A \bigtriangleup B),
\end{equation}
which is defined as a nonnegative extended real number.  If $C \in
\mathcal{A}$ too, then we get that
\begin{equation}
\label{d_mu(A, C) le d_mu(A, B) + d_mu(B, C)}
        d_\mu(A, C) \le d_\mu(A, B) + d_\mu(B, C),
\end{equation}
by (\ref{A bigtriangleup C = ...}).  Note that
\begin{equation}
\label{d_mu(A cap C, B cap C) le d_mu(A, B)}
        d_\mu(A \cap C, B \cap C) \le d_\mu(A, B)
\end{equation}
and
\begin{equation}
\label{d_mu(A cup C, B cup C) le d_mu(A, B)}
        d_\mu(A \cup C, B \cup C) \le d_\mu(A, B),
\end{equation}
by (\ref{(A cap C) bigtriangleup (B cap C) = (A bigtriangleup B) cap
  C}) and (\ref{(A cup C) bigtriangleup (B cup C) = ...}).

        Let us say that $A, B \in \mathcal{A}$ are
\emph{equivalent}\index{equivalent measurable sets} if
(\ref{d_mu(A, B) = mu(A bigtriangleup B)}) is equal to $0$.
This defines an equivalence relation on $\mathcal{A}$, and we let
$\widetilde{\mathcal{A}}$ denote the corresponding collection of
equivalence classes.  Put
\begin{equation}
\label{mathcal{A}_0 = {A in mathcal{A} : mu(A) < infty}}
        \mathcal{A}_0 = \{A \in \mathcal{A} : \mu(A) < \infty\},
\end{equation}
and observe that $A \bigtriangleup B \in \mathcal{A}_0$ when $A, B \in
\mathcal{A}_0$, so that (\ref{d_mu(A, B) = mu(A bigtriangleup B)}) is
finite.  Also let $\widetilde{\mathcal{A}}_0$ be the subset of
$\widetilde{\mathcal{A}}$ consisting of equivalence classes of
elements of $\mathcal{A}_0$.  It is easy to see that (\ref{d_mu(A, B)
  = mu(A bigtriangleup B)}) leads to a well-defined metric on
$\widetilde{\mathcal{A}}_0$, by standard arguments.

        Suppose for the moment that $X$, $\mathcal{A}$, and $\mu$
are as in Section \ref{unit interval}, so that $X = [0, 1]$,
$\mathcal{A}$ contains all closed subintervals of $[0, 1]$, $\mu([0,
  1]) < \infty$, and $\mu$ satisfies (\ref{lim_{t to a+} mu([a, t]) =
  0}) and (\ref{lim_{t to b-} mu([t, b]) = 0}).  In this case, one can
use finite unions of subintervals of $[0, 1]$ to get families of
elements of $\mathcal{A}$ that depend on arbitrarily many parameters.
More precisely, this leads to continuous families of elements of
$\widetilde{\mathcal{A}}$ with respect to the metric on
$\widetilde{\mathcal{A}}$ corresponding to (\ref{d_mu(A, B) = mu(A
  bigtriangleup B)}), because of (\ref{lim_{t to a+} mu([a, t]) = 0})
and (\ref{lim_{t to b-} mu([t, b]) = 0}).  Let us take $\mu$ to be as
in (\ref{mu([a, b]) = b - a}), for simplicity.  Under these conditions,
it is easy to see that $\widetilde{\mathcal{A}}$ has infinite
topological dimension with respect to the topology determined by the
metric corresponding to (\ref{d_mu(A, B) = mu(A bigtriangleup B)}),
using families of elements of $\mathcal{A}$ like these.

        Let $X$, $\mathcal{A}$, and $\mu$ be arbitrary again, let $k$
be a field with a $q$-absolute value function $|\cdot|$ for some $q >
0$, and let $V$ be a vector space over $k$ with a $q$-norm $N$ with
respect to $|\cdot|$ on $k$.  If $A \subseteq X$, then let ${\bf
  1}_A^k(x)$ be the $k$-valued indicator function on $X$ associated to
$A$, so that ${\bf 1}_A^{\bf R}$ is the real-valued indicator
function.  Let $v_0$ be an element of $V$, so that
\begin{equation}
\label{{bf 1}_A^k(x) v_0}
        {\bf 1}_A^k(x) \, v_0
\end{equation}
is a $V$-valued simple function on $X$, which is measurable when $A
\in \mathcal{A}$.  If $A, B \subseteq X$, then
\begin{eqnarray}
\label{N({bf 1}_A^k(x) v_0 - {bf 1}_B^k(x) v_0) = ...}
        N({\bf 1}_A^k(x) \, v_0 - {\bf 1}_B^k(x) \, v_0)
           & = & |{\bf 1}_A^k(x) - {\bf 1}_B^k(x)| \, N(v_0) \\
           & = & {\bf 1}_{A \bigtriangleup B}^{\bf R}(x) \, N(v_0) \nonumber
\end{eqnarray}
for every $x \in X$.  Note that (\ref{(A cap C) bigtriangleup (B cap
  C) = (A bigtriangleup B) cap C}) corresponds to multiplying this by
$|{\bf 1}_C^k(x)| = {\bf 1}_C^{\bf R}(x)$.  If $A, B \in \mathcal{A}$,
then it follows that
\begin{equation}
\label{||{bf 1}_A^k v_0 - {bf 1}_B^k v_0||_{L^r(X, V)} = ...}
        \|{\bf 1}_A^k \, v_0 - {\bf 1}_B^k \, v_0\|_{L^r(X, V)}
                     = \mu(A \bigtriangleup B)^{1/r} \, N(v_0)
\end{equation}
for every positive real number $r$.  Similarly,
\begin{equation}
\label{||{bf 1}_A^k v_0 - {bf 1}_B^k v_0||_{L^infty(X, V)} = ...}
        \|{\bf 1}_A^k \, v_0 - {\bf 1}_B^k \, v_0\|_{L^\infty(X, V)}
 = \|{\bf 1}_{A \bigtriangleup B}^{\bf R}\|_{L^\infty(X, {\bf R})} \, N(v_0)
\end{equation}
is equal to $N(v_0)$ when $\mu(A \bigtriangleup B) > 0$, and to $0$
otherwise.

\section{Measurable sets, continued}
\label{measurable sets, continued}
\setcounter{equation}{0}

        Suppose for the moment that $X$, $\mathcal{A}$, and $\mu$ are
as in Section \ref{unit interval} again, so that $X = [0, 1]$,
$\mathcal{A}$ contains all closed subintervals of $[0, 1]$, $\mu([0,
  1]) < \infty$, and $\mu$ satisfies (\ref{lim_{t to a+} mu([a, t]) =
  0}) and (\ref{lim_{t to b-} mu([t, b]) = 0}).  Thus for each $x \in
[0, 1]$ and $\epsilon > 0$ there is a $\delta(x) > 0$ such that
\begin{equation}
\label{mu([0, 1] cap (x - delta(x), x + delta(x))) < epsilon}
        \mu([0, 1] \cap (x - \delta(x), x + \delta(x))) < \epsilon.
\end{equation}
Because $[0, 1]$ is compact with respect to the standard topology,
there are finitely many points $x_1, \ldots, x_n$ in $[0, 1]$ such
that
\begin{equation}
\label{[0, 1] subseteq bigcup_j (x_j - delta(x_j)/2, x_j + delta(x_j)/2)}
 [0, 1] \subseteq \bigcup_{j = 1}^n (x_j - \delta(x_j)/2, x_j + \delta(x_j)/2).
\end{equation}
If we put
\begin{equation}
\label{delta = min_{1 le j le n} delta(x_j)/2}
        \delta = \min_{1 \le j \le n} \delta(x_j)/2,
\end{equation}
then we get that
\begin{equation}
\label{mu([0, 1] cap (x - delta, x + delta)) < epsilon}
        \mu([0, 1] \cap (x - \delta, x + \delta)) < \epsilon
\end{equation}
for every $x \in [0, 1]$.  More precisely, for each $x \in [0, 1]$, we
have that $x$ is contained in one of the intervals on the right side
of (\ref{[0, 1] subseteq bigcup_j (x_j - delta(x_j)/2, x_j +
  delta(x_j)/2)}).  This implies that
\begin{equation}
\label{(x - delta, x + delta) subseteq (x_j - delta(x_j), x_j + delta(x_j))}
 (x - \delta, x + \delta) \subseteq (x_j - \delta(x_j), x_j + \delta(x_j))
\end{equation}
for some $j$, by the definition (\ref{delta = min_{1 le j le n}
  delta(x_j)/2}) of $\delta$.  Thus (\ref{mu([0, 1] cap (x - delta, x
  + delta)) < epsilon}) follows from (\ref{(x - delta, x + delta)
  subseteq (x_j - delta(x_j), x_j + delta(x_j))}) and (\ref{mu([0, 1]
  cap (x - delta(x), x + delta(x))) < epsilon}), where the latter is
applied to $x_j$.

        Now let $X$ be a nonempty set, let $\mathcal{A}$ be an
algebra of subsets of $X$, and let $\mu$ be a finitely-additive
nonnegative measure on $(X, \mathcal{A})$.  Suppose that for each
$E \in \mathcal{A}$ with $\mu(E) < \infty$ and $\epsilon > 0$
there are finitely many measurable sets $A_1, \ldots, A_n \subseteq X$
such that
\begin{equation}
\label{mu(A_j) < epsilon}
        \mu(A_j) < \epsilon
\end{equation}
for each $j$, and
\begin{equation}
\label{bigcup_{j = 1}^n A_j = E}
        \bigcup_{j = 1}^n A_j = E.
\end{equation}
Put $E_0 = \emptyset$ and
\begin{equation}
\label{E_l = bigcup_{j = 1}^l A_j}
        E_l = \bigcup_{j = 1}^l A_j
\end{equation}
for $l = 1, \ldots, n$, so that $E_{l - 1} \subseteq E_l$ for $l = 1,
\ldots, n$ and $E_n = E$.  If $d_\mu(\cdot, \cdot)$ is as in
(\ref{d_mu(A, B) = mu(A bigtriangleup B)}), then we have that
\begin{equation}
\label{d_mu(E_{l - 1}, E_l) = mu(E_l setminus E_{l - 1}) le mu(A_l) < epsilon}
 d_\mu(E_{l - 1}, E_l) = \mu(E_l \setminus E_{l - 1}) \le \mu(A_l) < \epsilon
\end{equation}
for $l = 1, \ldots, n$, by (\ref{mu(A_j) < epsilon}).  Let
$\widetilde{\mathcal{A}}_0$ be as in the previous section, so that
$d_\mu(\cdot, \cdot)$ determines a metric on
$\widetilde{\mathcal{A}}_0$, as before.  It follows from this
discussion that $\widetilde{\mathcal{A}}_0$ is chain connected with
respect to this metric under these conditions.

        Similarly, let $k$ be a field with a $q$-absolute value
function $|\cdot|$ for some $q > 0$, and let $V$ be a vector space
over $k$ with a $q$-norm $N$.  Also let $\widetilde{S}_0(X, V)$ be as
in Section \ref{vector-valued functions}, and remember that
$\|\cdot\|_{L^r(X, V)}$ determines a $q$-norm on $\widetilde{S}(X, V)$
when $q \le r$, and an $r$-norm on $\widetilde{S}_0(X, V)$ when $0 < r
\le q$.  Under the conditions considered in the preceding paragraph,
one can check that $\widetilde{S}_0(X, V)$ is chain connected with
respect to the $q$ or $r$-metric corresponding to $\|\cdot\|_{L^r(X,
  V)}$ when $0 < r < \infty$.  More precisely, this uses the chain
connectedness of $\widetilde{\mathcal{A}}_0$ with respect to the
metric corresponding to $d_\mu(\cdot, \cdot)$.

        Let $X$ be any nonempty set again, let $\mathcal{A}$ be an
algebra of subsets of $X$, and let $\mu$ be a finitely-additive
nonnegative measure on $(X, \mathcal{A})$.  Suppose that
$\widetilde{\mathcal{A}}_0$ is chain connected with respect to the
metric associated to $d_\mu(\cdot, \cdot)$, and let us show that we
get the same type of conditions on $X$, $\mathcal{A}$, and $\mu$ as
before, in (\ref{mu(A_j) < epsilon}) and (\ref{bigcup_{j = 1}^n A_j =
  E}).  Let $E \subseteq X$ be a measurable set with $\mu(E) <
\infty$, and let $\epsilon > 0$ be given.  If
$\widetilde{\mathcal{A}}_0$ is chain connected, then there is an
$\epsilon$-chain in $\widetilde{\mathcal{A}}_0$ that connects the
elements of $\widetilde{\mathcal{A}}_0$ corresponding to $\emptyset$
and $E$.  Equivalently, this means that there are finitely many
measurable subsets $E_0, E_1, \ldots, E_n$ of $X$ such that $E_0 =
\emptyset$, $E_n = E$, $\mu(E_l) < \infty$ for each $l$, and
\begin{equation}
\label{d_mu(E_{l - 1}, E_l) < epsilon}
        d_\mu(E_{l - 1}, E_l) < \epsilon
\end{equation}
for $l = 1, \ldots, n$.  If we put
\begin{equation}
\label{E_l' = E_l cap E}
        E_l' = E_l \cap E
\end{equation}
for each $l = 0, 1, \ldots, n$, then we have that $E_0' = \emptyset$,
$E_n' = E$, and $E_l' \subseteq E$ for each $l$.  We also have that
\begin{equation}
\label{E_l' setminus E_{l - 1}' = ...}
        E_l' \setminus E_{l - 1}' = (E_l \setminus E_{l - 1}) \cap E
                                         \subseteq E_l \setminus E_{l - 1}
\end{equation}
for $l = 1, \ldots, n$, and hence
\begin{equation}
\label{mu(E_l' setminus E_{l - 1}') le ...}
        \mu(E_l' \setminus E_{l - 1}') \le \mu(E_l \setminus E_{l - 1})
                                       \le d_\mu(E_{l - 1}, E_l)
\end{equation}
Similarly, put
\begin{equation}
\label{E_l'' = bigcup_{j = 1}^l E_j'}
        E_l'' = \bigcup_{j = 1}^l E_j'
\end{equation}
for $l = 1, \ldots, n$, and $E_0'' = \emptyset$, so that $E_{l - 1}''
\subseteq E_l''$ for $l = 1, \ldots, n$ by construction.  Note that
$E_l'' \subseteq E$ for each $l$, because $E_l' \subseteq E$ for
each $l$, and that $E_n'' = E$.  It is easy to see that
\begin{equation}
\label{E_l'' setminus E_{l - 1}'' = ...}
        E_l'' \setminus E_{l - 1}'' = E_l' \setminus E_{l - 1}''
                             \subseteq E_l' \setminus E_{l - 1}'
\end{equation}
for $l = 1, \ldots, n$, which implies that
\begin{equation}
\label{mu(E_l'' setminus E_{l - 1}'') le mu(E_l' setminus E_{l - 1}')}
        \mu(E_l'' \setminus E_{l - 1}'') \le \mu(E_l' \setminus E_{l - 1}').
\end{equation}
Put $A_l = E_l'' \setminus E_{l - 1}''$ for $l = 1, \ldots, n$, so
that
\begin{equation}
\label{mu(A_l) < epsilon}
        \mu(A_l) < \epsilon
\end{equation}
for each $l$, by (\ref{d_mu(E_{l - 1}, E_l) < epsilon}), (\ref{mu(E_l'
  setminus E_{l - 1}') le ...}), and (\ref{mu(E_l'' setminus E_{l -
    1}'') le mu(E_l' setminus E_{l - 1}')}).  By construction, the
$A_l$'s are pairwise-disjoint measurable subsets of $X$ such that
\begin{equation}
\label{E_l'' = bigcup_{j = 1}^l A_j}
        E_l'' = \bigcup_{j = 1}^l A_j
\end{equation}
for $l = 1, \ldots, n$.  In particular, (\ref{E_l'' = bigcup_{j = 1}^l
  A_j}) is equal to $E_n'' = E$ when $l = n$, as desired.

\section{Continuous simple functions}
\label{continuous simple functions}
\setcounter{equation}{0}

        Let $X$ be a nonempty topological space, and let $Z$ be a
nonempty set.  A $Z$-valued function $f$ on $X$ is said to be
\emph{locally constant}\index{locally constant functions} at a point
$x \in X$ if $f$ is constant on an open subset of $X$ that contains
$x$.  Similarly, $f$ is said to be locally constant on $X$ if $f$ is
locally constant at every element of $X$.  It is easy to see that this
happens if and only if $f^{-1}(\{z\})$ is an open set in $X$ for every
$z \in Z$.  This implies that the inverse image of every subset of $Z$
under $f$ is an open set in $X$, since the union of any collection of
open subsets of $X$ is also an open set in $X$.  It follows that the
inverse image of every subset of $Z$ under $f$ is a closed set in $X$
too, because its complement is an open set.  In particular,
$f^{-1}(\{z\})$ is a closed set in $X$ for every $z \in Z$ when $f$ is
locally constant.

        Of course, if $f$ is locally constant on $X$, then $f$ is
continuous with respect to any topology on $Z$.  If $Z$ is equipped
with the discrete topology, then every continuous mapping into $Z$ is
locally constant.  If $X$ is connected, then every locally constant
function on $X$ is constant.  Conversely, if $X$ is not connected, and
if $Z$ has at least two elements, then there is a locally constant
mapping from $X$ into $Z$ which is not constant.  Similarly, $X$ is
totally separated if and only if the collection of locally constant
mappings from $X$ into any set $Z$ with at least two elements separates
points in $X$.

        Let $\mathcal{A}$ be the collection of subsets of $X$ that are
both open and closed.  It is easy to see that this defines an algebra
of subsets of $X$.  Also let $k$ be a field, and let $V$ be a vector
space over $k$.  In this case, a $V$-valued function $f$ on $X$ is
locally constant on $X$ if and only if
\begin{equation}
\label{f^{-1}({v}) in mathcal{A}, 2}
        f^{-1}(\{v\}) \in \mathcal{A}
\end{equation}
for every $v \in V$, by the earlier remarks.  Thus a $V$-valued simple
function $f$ on $X$ is measurable with respect to $\mathcal{A}$, as in
Section \ref{basic notions}, if and only if $f$ is locally constant.

        Let $CS(X, V)$\index{CS(X, V)@$CS(X, V)$} be the space of
$V$-valued simple functions on $X$ that are locally constant.  This is
the same as the space $S(X, V)$ defined in Section \ref{basic notions}
when $\mathcal{A}$ is as in the preceding paragraph.  In particular,
$CS(X, V)$ is a vector space over $k$ with respect to pointwise addition
and scalar multiplication, which can easily be verified directly as well.
Every element of $CS(X, V)$ is continuous with respect to any topology
on $V$, as before.  If $V$ is equipped with any topology that satisfies
the first separation condition, and if a $V$-valued simple function $f$
on $X$ is continuous with respect to this topology on $V$, then one
can check that $f$ is locally constant on $X$.

        If $f$ is a locally constant function on $X$ and $K \subseteq X$
is compact, then $f$ takes only finitely many values on $K$.  Suppose
that $f$ is a locally constant $V$-valued function on $K$, so that
\begin{equation}
\label{{x in X : f(x) ne 0}}
        \{x \in X : f(x) \ne 0\}
\end{equation}
is a closed set in $X$, which is the same as the support of $f$ in
this case.  If (\ref{{x in X : f(x) ne 0}}) is compact, then $f$ takes
only finitely many values on (\ref{{x in X : f(x) ne 0}}), and hence
$f$ takes only finitely many values on $X$ too.  Let $CS_{com}(X,
V)$\index{CS_com(X, V)@$CS_{com}(X, V)$} be the collection of locally
constant $V$-valued functions on $X$ with compact support.  Thus
$CS_{com}(X, V)$ is contained in $CS(X, V)$, by the previous remarks,
and in fact $CS_{com}(X, V)$ is a linear subspace of $CS(X, V)$.

        Let $\mathcal{A}_{com}$ be the collection of subsets of $X$
that are open, closed, and compact.  Of course, compact subsets of $X$
are automatically closed when $X$ is Hausdorff, and closed sets in $X$
are compact when $X$ is compact.  If $X$ is any topological space,
then $\emptyset \in \mathcal{A}_{com}$, and the union of any two
elements of $\mathcal{A}$ is an element of $\mathcal{A}_{com}$ as
well.  Similarly, if $A \in \mathcal{A}_{com}$, and $B \subseteq X$ is
both open and closed, then $A \cap B$ and $A \setminus B$ are elements
of $\mathcal{A}_{com}$.  If we put
\begin{equation}
\label{mathcal{A}_1 = ...}
        \mathcal{A}_1 = \{A \subseteq X : A \in \mathcal{A}_{com}
                             \hbox{ or } X \setminus A \in \mathcal{A}_{com}\},
\end{equation}
then one can check that $\mathcal{A}_1$ is an algebra of subsets of $X$,
which is obviously contained in the algebra $\mathcal{A}$ defined earlier.

        If $A \subseteq X$ is an open set and $\mathcal{B}$ is a base
for the topology of $X$, then $A$ can be expressed as a union of
elements of $\mathcal{B}$.  If $A \subseteq X$ is compact and open,
then it follows that $A$ can be expressed as the union of finitely
many elements of $\mathcal{B}$.  If $\mathcal{B}$ has only finitely or
countably many elements, then there can only be finitely or countably
many subsets of $X$ that are compact and open.  This implies that
$\mathcal{A}_{com}$ has only finitely or countably many elements, and
hence that $\mathcal{A}_1$ has only finitely or countably many
elements.

        If $X$ is compact, then $CS_{com}(X, V) = CS(X, V)$, and
$\mathcal{A}_1 = \mathcal{A}$.  Otherwise, suppose for the moment that
$X$ is not compact.  Let $CS_1(X, V)$ be the collection of locally
constant $V$-valued simple functions $f$ on $X$ such that
\begin{equation}
\label{X setminus f^{-1}({v_0}) is a compact subset of X}
        X \setminus f^{-1}(\{v_0\}) \hbox{ is a compact subset of } X
\end{equation}
for some $v_0 \in V$.  This is the same as saying that
\begin{equation}
\label{f_0(x) = f(x) - v_0 in CS_{com}(X, V)}
        f_0(x) = f(x) - v_0 \in CS_{com}(X, V),
\end{equation}
so that $CS_1(X, V)$ is the same as the linear span in $CS(X, V)$ of
$CS_{com}(X, V)$ and the space of $V$-valued constant functions on
$X$.  If $f$ is an element of $CS_1(X, V)$, $v_0 \in V$ is as in
(\ref{X setminus f^{-1}({v_0}) is a compact subset of X}), and $v \in
V \setminus \{v_0\}$, then
\begin{equation}
\label{f^{-1}({v}) subseteq X setminus f^{-1}({v_0})}
        f^{-1}(\{v\}) \subseteq X \setminus f^{-1}(\{v_0\}),
\end{equation}
which implies that
\begin{equation}
\label{f^{-1}({v}) is a compact subset of X}
        f^{-1}(\{v\}) \hbox{ is a compact subset of } X.
\end{equation}
In particular, (\ref{X setminus f^{-1}({v_0}) is a compact subset of
  X}) can hold for at most one element $v_0$ of $V$, because $X$ is
not compact.  It follows from (\ref{X setminus f^{-1}({v_0}) is a
  compact subset of X}) and (\ref{f^{-1}({v}) is a compact subset of
  X}) that every element of $CS_1(X, V)$ is measurable with respect to
$\mathcal{A}_1$.

        Conversely, suppose that $f$ is a $V$-valued simple function on
$X$ that is measurable with respect to $\mathcal{A}_1$.  In particular,
this means that $f$ is measurable with respect to $\mathcal{A}$, since
$\mathcal{A}_1 \subseteq \mathcal{A}$, so that $f$ is locally constant
on $X$.  If $v \in V$, then $f^{-1}(\{v\}) \in \mathcal{A}_1$, so that
either $f^{-1}(\{v\})$ is compact, or its complement in $X$ is
compact.  Because $f$ is a simple function on $X$, $f^{-1}(\{v\}) \ne
\emptyset$ for only finitely many $v \in V$.  If $f^{-1}(\{v\})$ is
compact for every $v \in V$, then $X$ can be expressed as the union of
finitely many compact sets, which implies that $X$ is compact.  Thus
$f^{-1}(\{v_0\})$ is not compact for at least one $v_0 \in V$ when $X$
is not compact.  This implies that $X \setminus f^{-1}(\{v_0\})$ is
compact for at least one $v_0 \in V$, because $f^{-1}(\{v_0\}) \in
\mathcal{A}_1$.  It follows that $f$ is an element of $CS_1(X, V)$,
so that $CS_1(X, V)$ consists of exactly the $V$-valued simple functions
on $X$ that are measurable with respect to $\mathcal{A}_1$ when $X$
is not compact.

        Suppose now that $X$ is a locally compact Hausdorff topological
space with topological dimension $0$.  If $K \subseteq X$ is compact,
$W \subseteq X$ is an open set, and $K \subseteq W$, then there is an
open set $U \subseteq X$ such that $K \subseteq U$ and $U$ is compact,
as in Section \ref{dimension 0}.  If $X$ is not compact, then one can
apply this with $W = X$ to get that the one-point compactification of
$X$ has topological dimension $0$ too.

        Let $|\cdot|$ be a $q$-absolute value function on $k$ for some
$q > 0$, and let $N$ be a $q$-norm on $V$ with respect to $|\cdot|$ on
$k$.  Also let $f$ be a continuous $V$-valued function on $X$, with respect
to the topology on $V$ determined by the $q$-metric associated to $N$.
If $U \subseteq X$ is compact and open, then one can approximate $f$
uniformly on $U$ by locally constant $V$-valued simple functions on $U$.
More precisely, for each $x \in U$, there is an open set $U(x) \subseteq U$
such that $x \in U$ and $f$ is almost constant on $U(x)$, because $f$
is continuous at $x$.  One can also choose $U(x)$ to be compact, because
$X$ has topological dimension $0$ at $x$.  It follows that there
are finitely many elements $x_1, \ldots, x_n$ of $U$ such that
\begin{equation}
\label{U = bigcup_{j = 1}^n U(x_j)}
        U = \bigcup_{j = 1}^n U(x_j),
\end{equation}
since $U$ is compact, and $U(x_j) \subseteq U$ for each $j$.  If we put
$U_1 = U(x_1)$ and
\begin{equation}
\label{U_l = U(x_l) setminus (bigcup_{j = 1}^{l - 1} U(x_j))}
        U_l = U(x_l) \setminus \Big(\bigcup_{j = 1}^{l - 1} U(x_j)\Big)
\end{equation}
for $l = 2, \ldots, n$, then $U_1, \ldots, U_n$ are pairwise-disjoint
compact open subsets of $U$ such that
\begin{equation}
\label{U = bigcup_{l = 1}^n U_l}
        U = \bigcup_{l = 1}^n U_l.
\end{equation}
By construction, $U_l \subseteq U(x_l)$ for each $l$, which implies
that $f$ is approximately constant on $U_l$ for each $l$.  Thus one
can approximate $f$ by a $V$-valued function that is constant on $U_l$
for each $l$.

        If $f$ has compact support in $X$, then we can first choose
a compact open set $U \subseteq X$ that contains the support of $f$.
The preceding argument permits us to approximate $f$ uniformly on $X$
by elements of $CS_{com}(X, V)$ with support contained in $U$.
Similarly, if $f$ vanishes at infinity on $X$, then there is a compact
set $K \subseteq X$ such that $f$ is small on $X \setminus K$, and we
can choose $U$ so that $K \subseteq U$.  In this case, the preceding
argument implies that $CS_{com}(X, V)$ is dense in $C_0(X, V)$ with
respect to the supremum norm.

\section{Lengths of chains}
\label{lengths of chains}
\setcounter{equation}{0}

        Let $M$ be a nonempty set, and let $d(x, y)$ be a $q$-metric
on $M$ for some $q > 0$.  If
\begin{equation}
\label{x_1, ldots, x_n}
        x_1, \ldots, x_n
\end{equation}
is a finite sequence of elements of $M$ and $a$ is a positive real
number, then let us define the \emph{$a$-length}\index{lengths of
  chains} of (\ref{x_1, ldots, x_n}) to be
\begin{equation}
\label{(sum_{j = 1}^{n - 1} d(x_j, x_{j + 1})^a)^{1/a}}
        \Big(\sum_{j = 1}^{n - 1} d(x_j, x_{j + 1})^a\Big)^{1/a},
\end{equation}
which is interpreted as being equal to $0$ when $n = 1$.  The analogue
of (\ref{(sum_{j = 1}^{n - 1} d(x_j, x_{j + 1})^a)^{1/a}}) with $a =
\infty$ is
\begin{equation}
\label{max_{1 le j < n} d(x_j, x_{j + 1})}
        \max_{1 \le j < n} d(x_j, x_{j + 1}),
\end{equation}
which should also be interpreted as being equal to $0$ when $n = 0$.
As in (\ref{||f||_infty le ||f||_r}) and (\ref{||f||_t le ||f||_r}),
(\ref{max_{1 le j < n} d(x_j, x_{j + 1})}) is less than or equal to
(\ref{(sum_{j = 1}^{n - 1} d(x_j, x_{j + 1})^a)^{1/a}}) for every $a >
0$, and (\ref{(sum_{j = 1}^{n - 1} d(x_j, x_{j + 1})^a)^{1/a}}) is
monotonically decreasing in $a$.  Note that (\ref{max_{1 le j < n}
  d(x_j, x_{j + 1})}) is less than some $\eta > 0$ exactly when
(\ref{x_1, ldots, x_n}) is an $\eta$-chain in $M$, as in Section
\ref{chain connectedness}.

        If $1 \le l \le m \le n$, then
\begin{equation}
\label{d(x_l, x_m)^q le ... le sum_{j = 1}^n d(x_j, x_{j + 1})^q}
        d(x_l, x_m)^q \le \sum_{j = l}^{m - 1} d(x_j, x_{j + 1})^q
                       \le \sum_{j = 1}^n d(x_j, x_{j + 1})^q,
\end{equation}
by the $q$-metric version of the triangle inequality.  This implies that
\begin{equation}
\label{max_{1 le l le m le n} d(x_l, x_m) le ...}
        \max_{1 \le l \le m \le n} d(x_l, x_m)
         \le \Big(\sum_{j = 1}^{n - 1} d(x_j, x_{j + 1})^q\Big)^{1/q},
\end{equation}
so that the diameter of (\ref{x_1, ldots, x_n}) in $M$ is less than or
equal to its $q$-length.  Similarly, if $d(\cdot, \cdot)$ is an
ultrametric on $M$, then
\begin{equation}
\label{max_{1 le l le m le n} d(x_l, x_m) le ..., 2}
 \max_{1 \le l \le m \le n} d(x_l, x_m) \le \max_{1 \le j < n} d(x_j, x_{j + 1}),
\end{equation}
which is the analogue of (\ref{max_{1 le l le m le n} d(x_l, x_m) le
  ...}) with $q = \infty$.

        If $0 < a < b < \infty$, then
\begin{equation}
\label{sum_{j = 1}^{n - 1} d(x_j, x_{j + 1})^b le ...}
        \sum_{j = 1}^{n - 1} d(x_j, x_{j + 1})^b
         \le \Big(\max_{1 \le j < n} d(x_j, x_{j + 1})\Big)^{b - a} \,
                   \sum_{j = 1}^{n - 1} d(x_j, x_{j + 1})^a.
\end{equation}
Equivalently, this means that
\begin{eqnarray}
\label{(sum_{j = 1}^{n - 1} d(x_j, x_{j + 1})^b)^{1/b} le ...}
\lefteqn{\Big(\sum_{j = 1}^{n - 1} d(x_j, x_{j + 1})^b\Big)^{1/b}} \\
 & \le &  \Big(\max_{1 \le j < n} d(x_j, x_{j + 1})\Big)^{1 - (a/b)} \,
 \bigg(\Big(\sum_{j = 1}^{n - 1} d(x_j, x_{j + 1})^a\Big)^{1/a}\bigg)^{a/b}.
                                                     \nonumber
\end{eqnarray}
In particular, if the $a$-length (\ref{(sum_{j = 1}^{n - 1} d(x_j,
  x_{j + 1})^a)^{1/a}}) of (\ref{x_1, ldots, x_n}) is bounded, and if
the maximal step size (\ref{max_{1 le j < n} d(x_j, x_{j + 1})}) is
small, then (\ref{(sum_{j = 1}^{n - 1} d(x_j, x_{j + 1})^b)^{1/b} le
  ...})  says that the $b$-length of (\ref{x_1, ldots, x_n}) should be
small too when $a < b$.



\printindex

\end{document}